\documentstyle[12pt, amstex,  epsf, amscd]{amsart}
\renewcommand{\baselinestretch}{1.35}
\newcommand\bR{{\bf R}}

\newcommand\rp{\bR P}
\newcommand\GL{{\rm GL}}
\newcommand\SL{{\rm SL}}
\newcommand\PGL{{\rm PGL}}
\newcommand\PSO{{\rm PSO}}

\newcommand\SLt{{\rm SL}}

\newcommand\dev{{\bf dev}}
\newcommand\SI{{\bf S}}

\newcommand\inte{{\rm int}}
\newcommand\Bd{{\rm bd}}
\newcommand\clo{{\rm Cl}}
\newcommand\bdd{{\bf d}}
\newcommand\ideal[1]{\tilde #1_\infty}
\newcommand\ra{\rightarrow}
\newcommand\longra{\longrightarrow}
\newcommand\che{\check}
\newcommand\emp{\emptyset}
\newcommand\eps{\epsilon}

\newcommand\vth{\vartheta}
\newcommand\vpi{\varphi}

\newcommand\ovl{\overline}
\newcommand\Aut{{\rm Aut}}

\newcommand\Rar{\Rightarrow}
\newcommand\uar{\uparrow}

\newcommand\awlg{may assume without loss of generality }

\newcommand{\cim}{\check \imath}

\newcommand{\hideal}[1]{#1_{h, \infty}}

\newtheorem{thm}{Theorem}[section]
\newtheorem{lem}{Lemma}[section]
\newtheorem{cor}{Corollary}[section]
\newtheorem{prop}{Proposition}[section]

\theoremstyle{definition} 
\newtheorem{defn}{Definition}[section]
\newtheorem{exmp}{Example}[section]

\theoremstyle{remark} 
\newtheorem{rem}{Remark}[section]

\begin{document}
\setlength{\baselineskip}{16pt}
\title[The $(n-1)$-convexity of real projective manifolds]
{Canonical decomposition of manifolds with flat real projective 
structure into $(n-1)$-convex manifolds and concave affine 
manifolds.}
\author{Suhyoung Choi}
\address{Department of Mathematics \\
College of Natural Sciences \\ 
Seoul National University \\ 
151--742 Seoul, Korea}
\email{shchoi@@math.snu.ac.kr}
\date{\today} 
\thanks{Research partially supported by GARC-KOSEF
and the Ministry of Education of Korea, 1997-001-D00036} 

\small\normalsize

\begin{abstract}
We try  to understand the geometric  properties of  $n$-manifolds 
($n\geq 2$) with geometric structures  modeled  on 
$(\bR P^n,  \PGL(n+1,  \bR))$, i.e., $n$-manifolds with projectively 
flat torsion free affine connections. We define the notion of 
$i$-convexity of such manifolds due to Carri\'ere for integers $i$, 
$1 \leq i \leq n-1$, which are generalization of convexity. Given 
a real  projective  $n$-manifold $M$, we  show that the failure of 
an $(n-1)$-convexity of $M$ implies an existence of a certain geometric 
object, $n$-crescent, in the completion $\che M$ of the universal 
cover $\tilde M$ of $M$.  We show that this  further implies the 
existence of a particular type of affine  submanifold  in $M$ and 
give a natural decomposition of $M$ into simpler real projective 
manifolds, some of which are $(n-1)$-convex and  others are  affine, 
more specifically concave affine. We feel that it is useful to have 
such decomposition particularly in dimension three. Our result will 
later aid us to study the geometric and topological properties of 
radiant affine $3$-manifolds leading to their classification. 
We get a consequence for affine Lie groups. 
\end{abstract}  

\maketitle

\section{Introduction} 

\renewcommand{\baselinestretch}{1.3}

From Ehresmann's definition of geometric structures on manifolds, 
a {\em real projective structure}\/ on a manifold is given by a maximal 
atlas of charts to $\bR P^n$  with  transition functions extending 
to projective transformations. This device lifts  the real projective 
geometry locally and consistently on a manifold. A differentio-geometric 
definition of a real projective structure is a projectively flat 
torsion-free connection. An equivalent way to define a real projective 
structure on a manifold $M$ is to give an immersion 
$\tilde M \ra \rp^n$, a so-called a developing map, 
equivariant with respect to a so-called holonomy homomorphism 
$h:\pi_1(M) \ra \PGL(n+1, \bR)$ where $\pi_1(M)$ is the group of 
deck transformations of the universal cover $M$ and $\PGL(n+1, \bR)$ 
is the group of deck transformations of $\rp^n$. Each of these 
description of a real projective structure gives rise to the unique 
descriptions of the other two kinds. For convenience, we will
assume that the dimension $n$ of manifolds is greater that or
equal to $2$ throughout this paper unless stated otherwise.

The global geometric and topological properties of real projective
manifolds are completely unknown, and are thought to be very complicated.
The study of real projective structure is a fairly obscure area with 
only handful of global results, as it is a very young field with many 
open questions, however seemingly unsovable by traditional methods. 
The complication comes from the fact 
that many compact manifolds are geodesically incomplete, and often 
the holonomy groups are far from being discrete lattices and thought 
to be far from being small such as solvable. There are some early 
indication that this field however offers many challenges for applying 
linear representations of discrete groups (which are not lattices),
group cohomology, classical convex and projective geometry, 
affine and projective differential geometry, real algebraic geometry, 
and analysis. (Since we cannot hope to mention them here appropriately, 
we offer as a reference the Proceedings of Geometry and Topology 
Conference at Seoul National University in 1997 \cite{ProcGT:97}.) 
This area is also an area closely related to the study of affine 
structures, which are more extensively studied with regard
to affine Lie groups. 

Riemannian hyperbolic manifolds admits a canonical real projective 
structure, via the Klein model of hyperbolic geometry as the hyperbolic 
space embeds as the interior of a standard ball in $\bR P^n$ and 
the isometry group $\PSO(1, n)$ as a subgroup of the group 
$\PGL(n+1, \bR)$ of projective automorphisms of $\bR P^n$ 
(see \cite{Gconv:90} and \cite{cdcr1}). 

They belong to the class of particularly understandable real projective 
manifolds which are convex ones. A {\em convex real projective manifold}\/ 
is  a quotient of a convex domain in an affine patch of $\bR  P^n$, 
i.e., the complement of a codimension one subspace with the natural 
affine structure of a complete affine space $\bR^n$,  by a  properly 
discontinuous and free action of a group of real projective 
transformations. It admits a Finsler metric, which has many nice 
geometric properties of a negatively curved Riemmanian manifold 
though the curvature is not bounded in the sense of Alexandrov. 
  
Affine manifolds naturally admit a canonical real projective structure 
since an affine space is canonically identified with the complement 
of codimension one subspace in the real projective space $\rp^n$ and 
affine automorphisms are projective. In particular, euclidean manifolds are 
projective. 

Not all real projective manifolds are convex (see \cite{ST:83} and 
\cite{Gconv:90}). However, in dimension two, we showed that closed 
real projective manifolds are built from convex surfaces. That is, 
a  compact real projective surface of negative Euler characteristic 
with geodesic boundary or empty boundary decomposes  along simple  
closed geodesics into convex surfaces (see  \cite{cdcr1}, \cite{cdcr2}, 
\cite{cdcr3},  and \cite{CG:97}). 

Also, recently, Benoist \cite{Benst:94} classified all real 
projective structures, homogeneous or not, on nilmanifolds some 
of which are not convex. Again, decompositions into parts 
admitting homogeneous structures were the central results. 
His student Dupont \cite{Dupont:98} classifies 
real projective structures on $3$-manifolds modelled on Sol.  

The real projective structures on $3$-manifolds are unexplored area, 
which may give us some insights into the topology of $3$-manifolds
along with hyperbolic or contact structures on $3$-manifolds. 

Let us state an interesting fact: All eight types of $3$-dimensional 
homogeneous Riemannian manifolds, i.e., manifolds with hyperbolic, 
spherical, euclidean, $\SI^2 \times \SI^1$-, $H^2 \times \SI^1$-, Sol, 
Nil, or $\SL(2, \bR)$-structures, admit canonical real projective structures 
since the models of each of the eight $3$-dimensional geometry can 
be realized as pairs of open domains in an affine or real projective 
space or their cover, and subgroups of groups of projective 
automorphisms of such domains. Thus, all Seifert spaces and atoridal 
Haken manifolds admit real projective structures.  

We might ask whether (i) real projective $3$-manifolds 
decomposes into pieces which admit one of the above geometries.
or (ii) conversely pieces with such geometric structures can be made 
into a real projective structures by perturbations. (These are 
questions by Thurston.) 

A question by Goldman (see \cite[p. 336]{Apan:97}) is that 
does all irreducible (Haken) $3$-manifolds admit real projective 
structure? A very exciting development will come from discovering 
ways to put real projective structures on $3$-manifolds other than 
from homogeneous Riemannian structures perhaps starting from 
triangulations of $3$-manifolds. 
(A related question asked by John Nash after his showing that all smooth 
manifolds admit real algebraic structure is that when does a manifold 
admit a rational structure, i.e., an atlas of charts with transition 
functions which are real rational functions. Real projective manifolds 
are rational manifolds with more conditions on the transition functions.)

These questions are at the moment very mysterious. This paper 
initiates some methods to study these with regard to the question 
(i). We will decompose real projective $n$-manifolds into concave 
affine real projective $n$-manifolds and $(n-1)$-convex real 
projective $n$-manifolds.

To do this, we will extend and refine the techniques involved in 
proving the result in dimension two. We work on $n\geq 2$ case 
although $n=2$ case was already done in \cite{cdcr1} and \cite{cdcr2}. 
The point where this paper improves the papers \cite{cdcr1} and 
\cite{cdcr2} even in $n=2$ case is that we will be introducing 
the notion of two-faced submanifolds which makes decomposition 
easier to understand. 

In three-dimensional case, our resulting decompositions into 
$2$-convex $3$-manifolds and concave affine $3$-manifolds seem 
to be along totally geodesic surfaces, which hopefully will be 
essential in $3$-manifold topology terminology. Thus, our remaining
task is to see if $2$-convex real projective $3$-manifolds 
admit nice decompositions or at least nice descriptions.


Our result will be used in the decomposition of radiant affine  
$3$-manifolds, which are $3$-manifolds with flat affine structure 
whose affine holonomy group fixes a point of an affine space 
(see \cite{rdh:97}). In particular, we will be proving there 
the Carri\`ere conjecture (see \cite{Carlet:93}) that every 
radiant affine $3$-manifolds admit a total section to the 
radial flow which exists naturally on radiant affine manifolds with 
the help from Barbot's work \cite{barbot:97a}, \cite{barbot:97b} 
(also see his survey article \cite{barbot:97c}). This will result 
in the classification of radiant affine $3$-manifolds.

Let us begin to state our theorems.
Let $T$ be an $(i+1)$-simplex in an affine space $\bR^n$, $i+1 < n$, 
with sides $F_1, F_2, \dots, F_{i+2}$. A real projective manifold is 
said to be $i$-convex if every real projective immersion 
\[T^o \cup F_2 \cup \dots \cup F_{i+2} \ra M\] 
extends to one from $T$ itself. 

\begin{thm}[Main]\label{thm:main} 
Suppose that $M$ is a compact real projective $n$-manifold with totally
geodesic or empty boundary. If $M$ is not $(n-1)$-convex, then $M$ 
includes a compact concave affine $n$-submanifold $N$ of type I or II 
or $M^o$ includes the two-faced $(n-1)$-submanifold. 
\end{thm} 

We will define the term ``two-faced $(n-1)$-submanifolds of type I 
and II'' in Definitions \ref{defn:twoface1} and \ref{defn:twoface2} 
which arise in separate constructions. But they are totally geodesic 
and are quotients of an open domain in an affine space by a group of 
projective transformations. They are canonically defined.
We will define the term concave affine $n$-submanifold in Definition 
\ref{defn:concaff}: A concave affine  $n$-manifold $M$ is a compact 
real projective manifold with concave boundary such that its 
universal cover is a union of  overlapping  $n$-crescents. 
An $n$-crescent is a convex $n$-ball whose bounding 
sides except one is in the ``infinity'' in the completion of 
the universal or holonomy cover (see Section \ref{sec:trcres}).
Their interiors are projectively diffeomorphic to either a half-space
or an open hemisphere. They are really generalization of half-spaces 
as one of the side is at ``infinity'' or ``missing''.
The manifold-interior $M^o$ of a concave affine 
manifold admits a projectively equivalent affine structure 
of very special nature. We expect them to be very limited.

Let $A$ be   a  properly imbedded $(n-1)$-manifold in $M^o$, 
which may or may not be two-sided and not necessarily connected or 
totally geodesic. The so-called splitting $S$  of $M$ along $A$ is 
obtained by completing $M - N$ by adding boundary which consists of 
either the union of two disjoint copies of components of $A$ or 
a double cover of components of $A$ (see the beginning of Section 
\ref{sec:corspl1}).

A manifold $N$ {\sl decomposes  into}\/ manifolds $N_1, N_2, \dots$ if
there exists a properly imbedded $(n-1)$-submanifold $\Sigma$ so that 
$N_i$ are components of the manifold obtained from splitting $M$ along 
$\Sigma$; $N_1, N_2, \dots$ are said to be the {\em resulting}\/ manifolds 
of the decomposition.

\begin{cor}\label{cor:spl1} Suppose that $M$ is compact but not 
$(n-1)$-convex.~Then 
\begin{enumerate}
\item after splitting  $M$ along  the  two-faced $(n-1)$-manifold  $A_1$
arising from hemispheric  $n$-crescents{\rm ,}  the resulting manifold
$M^{\sf s}$ decomposes properly into concave affine manifolds of type I 
and real projective $n$-manifolds with totally geodesic boundary which
does not include any concave affine manifolds of type I. 
\item We let $N$ be the disjoint union of the resulting manifolds of
the above decomposition other than concave affine ones.  After cutting 
$N$ along the two-faced $(n-1)$-manifold $A_2$ arising from  bihedral
$n$-crescents{\rm ,}  the resulting manifold $N^{\sf s}$ decomposes
into   concave affine  manifolds  of   type  II  and  real  projective
$n$-manifolds with convex boundary which is $(n-1)$-convex and 
includes no concave affine manifold of type II. 
\end{enumerate}
Furthermore, $A_1$ and $A_2$ are canonically defined and the decompositions
are also canonical in the following sense\/{\rm :}
If $M^{\sf s}$ equals $N \cup K$ for $K$ the union of concave affine 
manifolds of type I in $M^{\sf s}$ and $N$ the closure of the complement 
of $K$ includes no concave affine manifolds of type I, then the above 
decomposition agree with the decomposition into components of submanifolds 
in {\rm (1).} If $N^{\sf s}$ equals $S \cup T$ for $T$ the union of concave 
affine manifold of type II in $N^{\sf s}$ and $S$ the closure of 
the complement of $T$ that is $(n-1)$-convex and includes no concave 
affine manifold of type II, then the decomposition agree with 
the decomposition into components of submanifolds in {\rm (2).}
\end{cor} 

If $A_1 = \emp$, then we define $M^{\sf s}=M$ and if  $A_2 = \emp$, 
then we define $N^{\sf s} = N$. In  Section  \ref{sec:bitwof},  we  
will  give  a 2-dimensional example with a nontrivial splitting
(see Example \ref{exmp:twofacenec}).   

We note that $M$, $M^{\sf s}$, $N$, $N^{\sf s}$ have totally geodesic 
or empty boundary, as we will see in the proof. The final decomposed 
pieces of $N^{\sf s}$ are not so. Concave affine manifolds of type II 
have in general boundary concave seen from its inside and 
the $(n-1)$-convex real projective manifolds have convex boundary 
seen from inside (see Section \ref{sec:conv2}). 

Compare this corollary with what we have proved in \cite{cdcr1} and 
\cite{cdcr2} in the language of this paper, as the term ``decomposition''
is used somewhat differently there. 
\begin{thm}\label{thm:cdcr}
Let $\Sigma$ be a compact real projective surface with totally geodesic
or empty boundary. Suppose $\chi(\Sigma) < 0$. Then $\Sigma$ 
decomposes along the union of disjoint simple closed curves into 
convex real projective surfaces.
\end{thm}

Our Corollary \ref{cor:spl1} is strong enough to imply Theorem
\ref{thm:cdcr}, but we need to work out the classification of 
concave affine $2$-manifolds to do so. As a corollary to 
Corollary \ref{cor:spl1}, we get one further result. 
\begin{cor}\label{cor:cdcr}
Suppose $\chi(\Sigma) = 0$. Then $\Sigma$ decomposes into convex 
annuli, M\"obius bands, and concave affine $2$-manifolds of 
Euler characteristic $O$.
\end{cor}

This paper will be written as self-contained as possible on 
projective geometry and will use no highly developed 
machinery but will use perhaps many aspects 
of discrete group actions and geometric convergence in the 
Hausdorff sense joined in a rather complicated manner. Objects 
in this papers are all very concrete ones. To grasp these ideas, 
one only need to have some graduate student in geometry understanding 
and visualization of higher-dimensional projective and spherical 
geometry.

A {\em holonomy cover}\/ of $M$ is given as the cover of $M$ 
corresponding to the kernel of the developing map. We often need 
not look at the universal cover but the holonomy cover as it carries 
all information and we can define the developing map and holonomy 
homomorphism from it. The so-called Kuiper completion or projective 
completion of the universal or holonomy cover is the completion 
with respect to a metric pulled from $\SI^n$ by a developing map.
This notion was introduced by Kuiper for conformally flat manifolds
(see Kuiper \cite{Kp:50}).

In Section 1, we will give preliminary definitions and define and 
classify convex sets in $\SI^n$. In Section 2, we discuss the Kuiper 
or projective completions $\che M$ or $\che M_h$ of the universal 
cover $\tilde M$ or the holonomy cover $M_h$ respectively and convex 
subsets of them, and discuss how two convex subsets may intersect, 
showing that in the generic case they can be read from their images 
in $\SI^n$. 

The main ideas in this paper is to get good geometric objects 
in the universal cover of $M$. Loosely speaking, we illustrate 
our plan as follows:
\begin{itemize}
\item[(i)] For a compact manifold $M$ which is not $(n-1)$-convex, 
obtain an $n$-crescent in $\che M_h$.
\item[(ii)] Divide the case into two cases where $\che M_h$ includes
hemispheric $n$-crescents and where there are only $n$-crescents
which are bihedral. 
\item[(iii)] We derive a certain equivariance properties of 
hemispheric $n$-crescents or the unions of a collection of 
bihedral $n$-crescents equivalent to each other under 
the equivalence relation generated by the overlapping relation. 
That is, we show that any two of such sets either agree, are disjoint, 
or meet only in the boundary.
\item[(iv)] We show that the boundary where the two collections 
meet covers a closed codimension-one submanifold called 
the two-faced submanifolds. If we split along these, then 
the collection is now truly equivariant. From 
the equivariance, we obtain a submanifold covered by them 
called the concave affine manifold. This completes 
the proof of the Main Theorem.
\item[(v)] Apply the Main Theorem in sequence  
to prove Corollary \ref{cor:spl1}; that is, we split 
along the two-faced manifolds and obtain concave affine manifolds
for hemispheric $n$-crescent case and then bihedral 
$n$-crescent case.
\end{itemize}

In Section \ref{sec:n-1conv}, we prove a central theorem that 
given a real projective manifold which is not $(n-1)$-convex, we 
can find an $n$-crescent in the projective completion. The argument 
is the blowing up or pulling back argument as we saw in \cite{cdcr1}.

In Section \ref{sec:trcres}, we generalize the transversal intersection 
of crescents to that of $n$-crescents  (see \cite{cdcr1}). This shows 
that they intersect in a manageable manner as their sides in the frontier 
extend each other and the remaining sides intersecting transversally.

In   Section \ref{sec:hemispheres}, when  $\che  M_h$ includes 
a hemispheric $n$-crescent, we show how to obtain a two-faced 
$(n-1)$-submanifold. This is accomplished by the fact that 
two hemispheric crescents are either  disjoint, equal, or meet  
only in the  boundary, i.e., at a totally geodesic $(n-1)$-manifold. 
In  Section \ref{sec:bitwof}, we assume that $\che  M_h$ includes 
no hemispheric $n$-crescent but  includes  bihedral $n$-crescents. 
We define  equivalence classes of bihedral $n$-crescents. Two 
bihedral $n$-crescents   are  equivalent if there exists a chain 
of bihedral $n$-crescents overlapping with the next ones in the 
chain. This enables us to define $\Lambda(R)$ the union of 
$n$-crescents equivalent  to a given $n$-crescent $R$. Given 
$\Lambda(R)$ and $\Lambda(S)$  for two $n$-crescents $R$ and $S$, 
they are either disjoint, equal, or meet at a totally  geodesic 
$(n-1)$-submanifold. We obtain a two-faced $(n-1)$-submanifold 
from the totally geodesic $(n-1)$-submanifolds, which covers 
closed totally geodesic $(n-1)$-submanifolds in $M$.

In Section \ref{sec:crescres}, we show what happens to 
$n$-crescents if we take submanifolds or splits manifolds in the 
corresponding completions of the holonomy cover. They are 
all preserved.

In  Section \ref{sec:pfmain}, we  prove the Main Theorem: If there 
is no two-faced submanifold of type I, then two hemispheric $n$-crescents are 
either disjoint or equal. The union of all hemispheric $n$-crescents 
is invariant under deck transformations and hence covers a submanifold 
in  $M$, a concave affine manifold of type I.  
If there is no two-faced submanifold of type II, then $\Lambda(R)$ and 
$\Lambda(S)$ for two $n$-crescents $R$ and $S$ are either disjoint or 
equal. Again since the deck transformation group acts on the union of 
$\Lambda(R)$ for all $n$-crescents $R$, the union covers a manifold 
in $M$, which is a concave affine manifold of type II. 

In  Section \ref{sec:corspl1}, we prove Corollary \ref{cor:spl1}. 
That is, we decompose real projective manifolds. We show that when we have 
a two-faced submanifold, we can cut $M$ along these.  The result does   
not have a two-faced submanifold and hence can be decomposed   
into $(n-1)$-convex ones and properly  concave  affine manifolds 
as in Section \ref{sec:pfmain}. In Section \ref{sec:Liegp}, we will
show some consequence or modification of our result for affine Lie
groups. In Appendix A, we show that a real projective manifold is 
convex if and only if it is a quotient of a convex domain in $\SI^n$.
In Appendix B, we study some questions on shrinking sequences of convex 
balls in $\SI^n$ that are needed in Section \ref{sec:n-1conv}.

A  real projective structure on a  Lie group is {\sl left-invariant}\/
if left-multiplications preserve the real projective structure. 
The following theorem is also applicable to real projective structures 
on homogeneous manifolds invariant with respect to a proper group actions
(see Theorem \ref{thm:homg}).

\begin{thm}\label{thm:leftinv}
Let $G$ be a Lie group with  left-invariant  real projective
structure. Then either $G$ is   $(n-1)$-convex or $\tilde G$ is
projectively diffeomorphic to the universal cover of the complement of
a closed convex set in $\bR^n$ with induced real projective structure. 
\end{thm} 

The $(n-1)$-convexity of affine structures are defined similarly.
This theorem easily translates to one on affine Lie groups. 
\begin{cor}\label{cor:leftinvco}
Suppose that  $G$ has a left-invariant  affine structure. 
Then  either $G$  is  $(n-1)$-convex  or  $\tilde G$  is
affinely diffeomorphic  to the universal cover  of  a complement  
of a closed convex set in $\bR^n$ with induced affine structure. 
\end{cor}

We benefited greatly from conversations with Thierry Barbot, 
Yves Benoist, Yves Carri\`ere, William Goldman, Craig Hodgson, 
Michael Kapovich, Steven Kerckhoff, Hyuk Kim, Francois Labourie, 
John Millson, and William Thurston. 

\vfill
\break

\section{Convex subsets of the real projective sphere}
\label{sec:convSI}

In this foundational section, we will discuss somewhat slowly 
the real projective geometry of $\rp^n$ and the sphere $\SI^n$, and 
discuss convex subsets of $\SI^n$. We will give classification of 
convex subsets and give topological properties of them. We end with 
the geometric convergence of convex subsets. (We assume that the reader 
is familiar with convex sets in affine spaces, which are explain in 
Berger \cite{Berger:87} and Eggleston \cite{Egg:77} in detailed 
and complete manner.)

The real projective space $\rp^n$ is the quotient space of 
$\bR^{n+1} - \{O\}$ by the equivalence relation $\sim$ given by 
$x \sim y$ iff $x = sy$ for two nonzero vectors $x$ and $y$ and 
a nonzero real number $s$. The group $\GL(n+1, \bR)$ acts on 
$\bR^{n+1} - \{O\}$ linearly and hence on $\rp^n$, but not effectively.
However, the group $\PGL(n+1, \bR)$ acts on $\rp^n$ effectively. 
The action is analytic, and hence any element acting trivially in 
an open set has to be the identity transformation. (We will assume 
that $n \geq 2$ for convenience.)

Real projective geometry is a study of the invariant properties of
the real projective space $\rp^n$ under the action of $\PGL(n+1, \bR)$.
Given an element of $\PGL(n+1, \bR)$ we identify it with the
corresponding projective automorphism of $\rp^n$. 

Here by a real projective manifold, we mean an $n$-manifold with 
a maximal atlas of charts to $\rp^n$ where the transition functions 
are projective. This lifts all local properties of real projective 
geometry to the manifold. A real projective map is an immersion 
from a real projective $n$-manifold to another one which is projective
under local charts. More precisely, a function $f: M \ra N$ for 
two real projective $n$-manifolds $M$ and $N$ is real projective if 
it is continuous and for each pair of charts $\phi:U \ra \rp^n$ for $M$ 
and $\psi: V \ra \rp^n$ for $N$ such that $U$ and $f^{-1}(V)$ overlap, 
the function 
\[\psi\circ f \circ \phi^{-1}: \phi(U \cap f^{-1}(V)) \ra 
\psi(f(U) \cap V) \] 
is a restriction of an element of $\PGL(n+1, \bR)$ (see Ratcliff
\cite{RT:94}).

It will be very convenient to work on the simply connected, spheres 
$\SI^n$ the double cover of $\rp^n$ as $\SI^n$ is orientable and 
it is easier to study convex sets. We may identify the standard 
unit sphere $\SI^n$ in $\bR^{n+1}$ with the quotient space of
$\bR^{n+1} - \{O\}$ by the equivalence relation $\sim$ given by 
$x \sim y$ if $x = sy$ for nonzero vectors $x$ and $y$ and $s > 0$. 
As above $\GL(n+1, \bR)$ acts on $\SI^n$. The subgroup 
$\SLt_\pm(n+1, \bR)$ of linear maps of determinant $\pm 1$ acts on 
$\SI^n$ effectively. We see easily that $\SLt_\pm(n+1, \bR)$ is a 
double cover of $\PGL(n+1, \bR)$. We denote by $\Aut(\SI^n)$ the 
isomorphic group of automorphisms of $\SI^n$ which is induced by 
an element of $\SLt_\pm(n+1, \bR)$.

Since $\rp^n$ has an obvious chart to itself, namely the identity map,
it has a maximal atlas containing this chart. Hence, $\rp^n$ has 
a real projective structure. Since $\SI^n$ is a double cover of 
$\rp^n$, and the covering map $p$ is a local diffeomorphism, it follows 
that $\SI^n$ has a real projective structure. We see easily that 
each element of $\Aut(\SI^n)$ are real projective maps. Conversely,
each real projective automorphism of $\SI^n$ is an element of 
$\Aut(\SI^n)$ as the actions are locally identical with those 
of elements of $\Aut(\SI^n)$. There is a following 
convenient commutative diagram:
\begin{equation}\label{eqn:autolift}
\begin{array}{ccc}
\SI^n & \stackrel{g}{\longra} & \SI^n  \\
\downarrow p&			& \downarrow p \\
\rp^n & \stackrel{g'}{\longra} & \rp^n 
\end{array}
\end{equation}
where given a real projective automorphism $g$, a real projective map 
$g'$ always exists and given $g'$, we may obtain $g$ unique up to the 
antipodal map $A_{\SI^n}$ which sends $x$ to $-x$ for each unit vector 
$x$ in $\SI^n$. (Note that $\SI^n$ with this canonical real projective
structure is said to be a {\bf {\em real projective sphere}}.)   

The standard sphere has a standard Riemannian metric $\mu$ of 
curvature $1$. We denote by $\bdd$ the distance metric on $\SI^n$ 
induced from $\mu$. The geodesics of this metric are arcs on a great 
circles parameterized by $\bdd$-length. This metric is 
projectively flat, and hence geodesics of the metric agrees with 
projective geodesics up to choices of parametrization.

A {\em convex line}\/ is an embedded geodesic in $\SI^n$ of $\bdd$-length
less than or equal to $\pi$. A {\em convex set}\/ is a subset of 
$\SI^n$ such that any two points of $A$ is connected by a convex 
segment in $A$.
A {\em simply convex}\/ subset of $\SI^n$ is a convex subset such 
that every pair of point is connected by a convex segment of 
$\bdd$-length $< \pi - \eps$ for a positive number $\eps$. 
(Note that all these are projectively invariant properties.)
A singleton, i.e., the set consisting of a point, is convex 
and simply convex.

A {\em great $0$-dimensional sphere}\/ is the set of points antipodal 
to each other. This is not convex. A {\em great $i$-dimensional sphere}\/ 
in $\SI^n$ for $i \geq 1$ is convex but not simply convex. 
A {\em great $i$-dimensional hemisphere}, $i \geq 1$, is the closure 
of a component of a great $i$-sphere $\SI^i$ removed with 
a great $(i-1)$-sphere $\SI^{i-1}$ in $\SI^i$. It is a convex 
but not simply convex subset. A {\em $0$-dimensional hemisphere}\/
is simply a singleton.

Given a codimension one subspace $\rp^{n-1}$ of $\rp^n$, the complement 
of $\rp^n$ can be identified with an affine space $\bR^n$ so that 
geodesic structures agree, i.e., the projective geodesics are affine ones 
and vice versa up to parameterization. Given an affine space $\bR^n$, 
we can compactify it to a real projective space $\rp^n$ by adding points
(see Berger \cite{Berger:87}). Hence the complement $\rp^n - \rp^{n-1}$ 
is called an {\em affine patch}. An open $n$-hemisphere in $\SI^n$ 
maps homeomorphic onto $\rp^n - \rp^{n-1}$ for a subspace $\rp^{n-1}$. 
Hence, the open $n$-hemisphere has a natural affine structure 
of $\bR^n$ whose geodesic structure is same as that of 
the projective structure. An open $n$-hemisphere is sometimes 
called an {\em affine patch}. A bounded set in $\bR^n$ convex in the 
affine sense is convex in $\SI^n$ by our definition when $\bR^n$ 
is identified with the open $n$-hemisphere in this manner.


We give a definition given in \cite{RT:94}: A pair of points $x$ and $y$
is {\em proper}\/ if they are not antipodal. A minor geodesic connecting 
a proper pair $x$ and $y$ is the shorter arc in the great circle passing
through $x$ and $y$ with boundary $x$ and $y$.

The following theorem shows the equivalence of our definition to
the definition given in \cite{RT:94} except for pairs of antipodal points.
\begin{thm}\label{thm:eqconv}
A set $A$ is a convex set or a pair of antipodal points if and only if 
for each proper pair of points $x$, $y$ in $A$, $A$ includes a minor 
geodesic $\ovl{xy}$ in $A$ connecting $x$ and $y$.
\end{thm}
\begin{pf}
If $A$ is convex, then given two proper pair of points the convex segment 
in $A$ connecting them is clearly a minor geodesic. A pair of antipodal
points has no proper pair.

Conversely, let $x$ and $y$ be two points of $A$. If $x$ and $y$ are 
proper then since a minor geodesic is convex, we are done. If $x$ and 
$y$ are antipodal, and $A$ equals $\{x, y\}$, then we are done.
If $x$ and $y$ are antipodal, and there exists a point $z$ in $A$ 
distinct from $x$ and $y$, then $A$ includes the minor segment $\ovl{xz}$
and $\ovl{yz}$ and hence $\ovl{xz} \cup \ovl{yz}$ is a convex segment 
connecting $x$ and $y$; $A$ is convex.
\end{pf}

By above theorem, we see that our convex sets satisfy the properties in
Section 6.2 of \cite{RT:94}. Let $A$ be a nonempty
convex subset of $\SI^n$. The {\em dimension}\/ of $A$ is defined to be the
least integer $m$ such that $A$ is included in a great $m$-sphere in 
$\SI^n$. If $\dim(A) = m$, then $A$ is included in a unique great 
$m$-sphere which we denote by $<A>$. The interior of $A$, denoted 
by $A^o$, is the topological interior of $A$ in $<A>$, and the boundary 
of $A$, denoted by $\partial A$, is the topological boundary of $A$ 
in $<A>$. The closure of $A$ is denoted by $\clo(A)$ and is 
a subset of $<A>$. $\clo(A)$ is convex and so is $A^o$. (These are 
from Ratcliff \cite{RT:94}.) Moreover, the intersection of two convex 
sets is either convex or is a pair of antipodal points by the above 
theorem. Hence, the intersection of two convex sets is convex if it 
contains at least three points, and it contains a pair of nonantipodal 
points, or one of the sets does not contain a pair of antipodal points.

A {\em convex hull}\/ of a set is the minimal convex containing the set.

\begin{lem}\label{lm:intconv}
Let $A$ be a convex set. $A^o$ is not empty unless $A$ is empty.
\end{lem}
\begin{pf}
Let $<A>$ have dimension $k$. Then $A$ has to have at least $k+1$ points 
$p_1, \dots, p_{k+1}$ in general position as unit vectors in $\bR^{n+1}$ 
since otherwise every $(k+1)$-tuple of vectors are dependent and $A$ is 
a subset of a great sphere of lower dimension. The convex hull of 
the points $p_1, \dots, p_{k+1}$ is easily shown to be a spherical 
simplex with vertices $p_1, \dots, p_{k+1}$. The simplex is obviously
a subset of $A$, and the interior of the simplex is included in $A^o$.
\end{pf}

We give the following classification of convex sets in the following 
two theorems.
\begin{thm}\label{thm:class1} 
Let $A$ be a convex subset of $\SI^n$. Then $A$ is one of 
the following sets\/{\rm :}
\begin{itemize}
\item a great sphere $\SI^i$, $1 \leq i \leq n$,
\item an $i$-dimensional hemisphere $H^i$, $0 \leq i \leq n$,
\item a proper convex subset of an $i$-hemisphere $H^i$.
\end{itemize}
\end{thm} 
\begin{pf} 
We will prove by induction on dimension $m$ of $<A>$. The theorem is obvious 
for $m =0, 1$. Suppose that the theorem holds for $m =k-1$, $k \geq 2$.

Suppose now that the dimension of $A$ equals $m$ for $m = k$.
Let us choose a hypersphere $\SI^{m-1}$ in $<A>$ intersecting 
with $A^o$. Then $A_1 = A \cap \SI^{m-1}$ is as one of the above 
(1), (2), (3). The dimension of $A_1$ is at least one, i.e., 
$m-1 \geq 1$. Suppose $A_1 = \SI^{m-1}$. As $A^o$ has two points $x, y$ 
respectively in components of $<A> -\SI^{m-1}$, taking the union of 
segments from $x$ to points of $\SI^{m-1}$, and segments from $y$ 
to points of $\SI^{m-1}$, we obtain that $A = <A>$.

If $A_1$ is as in (2) or (3), then choose an $(m-1)$-hemisphere $H$
containing $A_1$ with boundary a great $(m-2)$-sphere $\partial H$.
Consider the collection $\cal P$ of all $(m-1)$-hemispheres 
including $\partial H$. Then $\cal P$ has a natural real projective 
structure of a great circle, and let $A'$ be the set of 
the $(m-1)$-hemispheres in $\cal P$ whose interior meets $A$. Then 
since a convex segment in $<A> - \partial H$ project
to a convex segment in the circle $\cal P$, it follows that $A'$ 
has the property that any two proper pair of points of $A'$ is 
connected by a minor geodesic, and by Theorem \ref{thm:eqconv} $A$
is either a pair of antipodal points or a convex subset. 

Let $-H$ denote the closure of the complement of $<A> - H$. 
Then the interior of $-H$ do not meet $A$ as it does not meet $A_1$.
Hence $A'$ is a subset of ${\cal P} - \{-H\}$.

If $A'$ is a pair of antipodal points, then $A'$ must be $\{H, -H\}$,
and this is a contradiction. Since $A'$ is a proper convex subset of 
$\cal P$, $A'$ must be a convex subset of a $1$-hemisphere $I$ in 
$\cal P$. This means that only the interior of $(m-1)$-hemispheres 
in $I$ meets $A$, and there exists an $m$-hemisphere in $<A>$ 
including $A$. Thus $A$ either equals this $m$-hemisphere or 
a proper convex subset of it.    
\end{pf}

\begin{thm}\label{thm:class2} 
Let $A$ be a proper convex subset of an $i$-dimensional hemisphere 
$H^i$ for $i \geq 1$. Then exactly one of the following holds\/{\rm :} 
\begin{itemize} 
\item $\partial A$ contains a unique maximal great $j$-sphere $\SI^j$ for 
some $0 \leq j \leq i-1$, which must be in $\partial H^i$ and its closure is 
the union of $(j+1)$-hemispheres with common boundary $\SI^j$, or 
\item $A$ is a simply convex subset of $H^i$, in which case 
$A$ can be realized as a bounded convex subset of perhaps another 
open $i$-hemisphere $K^i$ identified with an affine space $\bR^i$. 
\end{itemize}
\end{thm}
\begin{pf}
We assume without loss of generality that $A$ is closed by taking 
a closure of $A$ if necessary. If $A$ includes a pair of antipodal
points, then $A$ do not satisfy the Kobayashi's criterion \cite{cc:93}. 
Then by Section 1.4 of \cite{cc:93}, we have the first item. 

If $A$ includes no pair of antipodal points, then let $m$ be the 
dimension of $<A>$ and we do the induction over $m$. If $m =0, 1$,
then the second item is obvious. Suppose we have the second item holding 
for $m = k-1$, where $k \geq 2$. Now let $m = k$, and choose a great sphere 
$\SI^{m-1}$ meeting $A^o$, and let $A_1 = A \cap \SI^{m-1}$. Then 
$A_1$ is another simply convex set. Hence, $A_1$ is a bounded 
convex subset of an open $(m-1)$-hemisphere $K$ identified as an affine 
space $\bR^{m-1}$. Hence $A_1$ does not meet $\partial K$. As in the proof 
of Theorem \ref{thm:class1}, we let $\cal P$ be the set of all 
$(m-1)$-hemispheres with boundary in $\partial K$, which has a natural 
real projective structure of a great circle. As in the proof, we see 
that the subset $A'$ of $\cal P$ consisting of hemispheres whose 
interior meets $A$ is a convex subset of a $1$-hemisphere in $\cal P$. 
The boundary of $A'$ consists of two hemispheres $H_1$ and $H_2$. Since 
$A'$ is connected, $H_1$ and $H_2$ bound a convex subset $L$ in $<A>$ 
and $H_1$ and $H_2$ meet in a $\mu$-angle less than or equal to $\pi$. 

If the angle between $H_1$ and $H_2$ equals $\pi$, then $H_1 \cup H_2$ 
is a great $(m-1)$-sphere, and $H_1^o$ and $H_2^o$ includes two points 
$p, q$ of $A$ respectively which are not antipodal. Since $A$ is convex, 
$\ovl{pq}$ is a subset of $A$; since $p$ and $q$ is not antipodal, 
$\ovl{pq}$ meets $\partial K$ by geometry, a contradiction. 

Since the angle between $H_1$ and $H_2$ is less then $\pi$, it is now
obvious that there exists an $m$-hemisphere $H$ containing $A$ and meeting
$L$ only at $\partial K$. Hence $A$ is a convex subset of $H^o$. 
Since $A$ is compact, $A$ is a bounded convex subset of $H^o$.
\end{pf}

An {\sl $m$-bihedron} in $\SI^n$ is the closure of a component of 
a great sphere $\SI^m$ removed with two great spheres 
of dimension $m-1$ in $\SI^m$ ($m \geq 1$). A $1$-bihedron is 
a simply convex segment.

\begin{lem}\label{lem:classf} A compact convex subset $K$ of $\SI^n$ 
including an $(n-1)$-hemisphere is either the sphere $\SI^n$ itself, 
a great $(n-1)$-sphere, an $n$-hemisphere, or an $n$-bihedron.
\end{lem}
\begin{pf} 
Let $H$ be the $(n-1)$-hemisphere in $K$ and $s$ the great circle 
perpendicular to $H$ at the center of $H$. Then since $K$ is convex, 
$s \cap K$ is a convex subset of $s$ or a pair of antipodal point
as in the proof of Theorem \ref{thm:class1}. If $s \cap K = s$, 
then every segment from a point of $s$ to a point of $H$ belongs 
to $K$ by convexity. Thus, $K = \SI^n$. Depending on whether 
$s \cap K$ is a pair of antipodal points, a great segment or a simply 
convex segment (see \cite{cdcr1}), $K$ is a great $(n-1)$-sphere, 
an $n$-hemisphere or an $n$-bihedron. 
\end{pf} 

We will now discuss the topological properties of convex sets.
\begin{thm}\label{thm:topconv}
Let $A$ be a convex $m$-dimensional subset of $\SI^n$ other than 
a great sphere. Then $A^o$ is homeomorphic to an open $m$-ball,
$\clo(A)$ the compact $m$-ball, and $\partial A$ to the sphere 
of dimension $m-1$.
\end{thm}
\begin{pf}
If $A$ is zero or one dimensional, then the theorem is obvious. 
Assume that $m \geq 2$ from now on. $A$ is a subset of a closed 
$m$-hemisphere $H$. Choose a point $p$ of $H$ in $A^o$, which must 
be a point of $A^o \cap H^o$. Then for each point $q$ of $A$, 
$[p, q)$, $[p, q) = \ovl{pq} - \{q\}$ is a subset of $A^o$ 
(see Theorem 6.2.2 of \cite{RT:94}) where $\ovl{pq}$ denotes 
the unique convex segment connecting $p$ and $q$. Hence, for each ray 
$r$ from $p$ to a point of $\partial H$, $r \cap A$ is a ray with 
endpoints $r$ and a point $q_r$ in $r \cap \partial A$ so that 
$[p, q_r) \subset A^o$. We define a function $f: \partial H \ra \bR$ by
letting $f(x)$ equal the $\bdd$-length of $\ovl{pr_q}$ where $r$ is 
the ray from $p$ to $x$. Then $f$ is obviously bounded. 

We claim that $f$ is a continuous function. Suppose not. Then 
there exists a sequence $y_i$ in $\partial H$ converging to a point $y$
of $\partial H$ so that $|f(y_i) - f(y)| > \delta$ for a small 
positive number $\delta$. Then we see that the corresponding 
$q_{r_i}$ for $r_i$ a ray connecting $p$ and $y_i$ converges 
the limit $y^*$ on the ray $r$ distinct from $y$. Since $\partial A$ 
is closed, we see that $y^*$ is a point of $\partial A$. Whether $y$ 
lies within $[p, q_r)$ or in $r - [p, q_r]$, we get contradiction by 
Theorem 6.2.2 of \cite{RT:94}.

Now, we can follow Section 11.3.1 of Berger \cite{Berger:87} to see that 
$A^o$ is homeomorphic to $H^o$, $A$ to $H$, and $\partial A$ to the 
sphere of dimension $m-1$.
\end{pf}

Let $A$ be an arbitrary subset of $\SI^n$. A {\em supporting}\/ hypersphere 
$L$ for $A$ is a great $(n-1)$-sphere containing $x$ in $A$ such that 
the two closed hemispheres determined by $L$ includes $A$ and $x$ 
respectively. We say that $L$ is the {\em supporting}\/ hypersphere 
for $A$ at $x$. 

\begin{thm}\label{th:supp}
Let $A$ be a convex subset of $\SI^n$, other than $\SI^n$ itself.
Then for each point of $\partial A$, there exists a supporting hypersphere 
for $A$ at $x$. 
\end{thm}
\begin{pf}
If the dimension $i$ of $A$ is $0$, this is trivial. Assume $i \geq 1$.
If $A$ is a great $i$-sphere or an $i$-hemisphere $i \geq 1$, it is obvious. 
If not, then $A$ is contained an $i$-hemisphere, say $H$. Then $A^o$ is 
a convex subset of the affine space $H^o$. Hence, there exists 
a supporting hyperplane $K$ for $A^o$ at $x$ by Proposition 11.5.2 of 
\cite{Berger:87}. The hyperplane $K$ equals $L \cap H^o$ for a great 
$(i-1)$-sphere $L$ in $<A>$. Thus any great $(n-1)$ sphere $P$ meeting 
$<A>$ at $\partial H$ is the supporting hypersphere for $A$ at $x$.
\end{pf}

We define a Hausdorff distance between all compact subsets of $\SI^n$. 
We say that two compact subsets $X, Y$ have distance less than $\eps$, 
if $X$ is in a $\bdd$-$\eps$-neighborhood of $Y$ and $Y$ is in one of $X$.
This defines a metric on the space of all compact subsets of $\SI^n$. 

Suppose that a sequence of compact sets $K_i$ converges to $K_\infty$. 
If $x \in K_\infty$, then by definition for any positive number $\eps$, 
there exists an $N$ so that for $i > N$, there exists a point $x_i \in K_i$ 
so that $\bdd(x, x_i) < \eps$. Also, given a point $x$ of $\SI^n$, so that 
a sequence $x_i \in K_i$ converges to $x$, then $x$ lies in $K_\infty$.
If otherwise, $x$ is at least $\delta$ away from $K_\infty$ for 
$\delta >0$, and so the $\delta/2$-$\bdd$-neighborhood of $K_\infty$ 
is disjoint from an open neighborhood $J$ of $x$. But since 
$x_i \in J$ for $i$ sufficiently large, this contradicts $K_i \ra 
K_\infty$.

\begin{thm}\label{thm:gconv}
Given a sequence of compact convex subsets $K_i$ of $\SI^n$, we can 
always choose a subsequence converging to a subset $K_\infty$. 
$K_\infty$ is compact and convex. Also if $K_i$ are $i$-balls, then 
$K_\infty$ is a convex ball of dimension less than or equal to $i$. 
If $\dim K_\infty = n$, then we have $\bigcup_{i = 1}^\infty K_i^o 
\supset K_i^o$. In this case $\partial K_i \ra \partial K_i$.
\end{thm}
\begin{pf}
The first statment follows from the well-known compactness of the spaces
of compact subsets of compact Hausdorff spaces under Hausdorff metrics.

For each point $x$ of $K_\infty$, there exists a sequence $x_i \in K_i$ 
converging to $x$. Choose arbitrary two points $x$ and $y$ of $K_\infty$, and 
sequences $x_i \in K_i$ and $y_i \in K_i$ converging to $x$ and $y$ 
respectively. Then there exists a segment $\ovl{x_iy_i}$ of $\bdd$-length
$\leq \pi$ in $K_i$ connecting $x_i$ and $y_i$. Since the sequence 
of $\ovl{x_iy_i}$ is a sequence of compact subsets of $\SI^n$, we may 
assume that a subsequence converges to a compact subset $L$ of $\SI^n$. 
By the above paragraph $L \subset K_\infty$. By elementary geometry, 
it is easy to see that $L$ is a segment of $\bdd$-length $\leq \pi$. 
Thus $K_\infty$ is convex. 

If $K_i$ are $i$-balls, then $K_i \subset H_i$ for $i$-hemispheres $H_i$. 
We choose a subsequence $i_j$ of $i$ so that $H_{i_j}$ converges to
an $i$-hemisphere $H$. If follows that $K_\infty$ is a subset of $H$ 
by the paragraph above our lemma since $K_{i_j}$ converges to $K_\infty$.
Thus, $K_\infty$ is an compact convex subset of $H_\infty$, which 
shows that $K_\infty$ is a convex ball of dimension $\leq i$.

The third and final statements follow as in Section 2 of Appendix of 
\cite{cdcr1}. The dimension does not play a role.
\end{pf}

\vfill
\break

\section{Convex subsets in the Kuiper completion of 
the universal and holonomy covers.} 
\label{sec:conv2}
In this second foundational section, we begin by lifting 
the development pair to the real projective sphere $\SI^n$. 
To make our discussion much more simpler, 
we will define a completion, called a Kuiper 
completion or projective completion, by inducing the Riemmanian metric 
of the sphere to the universal cover $\tilde M$ or the holonomy cover 
$M_h$ of $M$ and then completing them in the Cauchy sense. Then 
we define the ideal set to be the completion removed
with $\tilde M$ or $M_h$, i.e., points infinitely far away from 
points of $\tilde M$ or $M_h$.  

We will define convex sets in these completions, which are always 
isomorphic to ones in $\SI^n$. Then we will introduce $n$-crescents, 
which are convex $n$-balls in the completions where a side or 
an $(n-1)$-hemisphere in the boundary lies in the ideal sets. 
We show how two convex subsets of the completion may intersect; their 
intersection properties are described by their images in $\SI^n$ 
under the developing map. Finally, we describe the dipping intersection, 
the type of intersection which will be useful in this paper, and on 
which our theory of $n$-crescents depend heavily as we will see in Section
\ref{sec:trcres}.

We will assume that our manifolds in this paper have dimension $\geq 2$ 
unless stated otherwise. Let $M$ be a real projective $n$-manifold.
Then $M$ has a development pair $(\dev, h)$ of an immersion 
$\dev:\tilde M \ra \rp^n$, called a developing map, 
and a  holonomy   homomorphism   $h:  \pi_1(M)  \ra  \PGL(n+1,   \bR)$
satisfying $\dev\circ \gamma =  h(\gamma)\circ \dev$ for every $\gamma
\in \pi_1(M)$. Such a pair is determined up to an action of an element 
$\vth$ of $\PGL(n+1, \bR)$ as follows:
\begin{equation}\label{eqn:devconj}
(\dev, h(\cdot)) \mapsto 
(\vth\circ \dev, \vth \circ h(\cdot) \circ \vth^{-1}).
\end{equation}
Developing maps are obtained by analytically extending coordinate charts 
in the atlas. Holonomy homomorphisms are obtained from the chosen 
developing map. See Ratcliff \cite{RT:94} for more details.
The development pair characterizes the real projective structure, and
hence another way to give a real projective structure to a manifold 
is to find a pair $(f, k)$ where $f$ is an immersion $\tilde M \ra \rp^n$ 
which is equivariant with respect to the homomorphism $k$ from the 
group of deck transformations to $\PGL(n+1, \bR)$.

We  assume that the manifold-boundary  $\delta  M$ of a real projective
manifold $M$ is {\em totally geodesic}\/ unless stated otherwise. This 
means that for each point of $\delta M$, there exist an open 
neighborhood $U$ and a lift $\phi: U \ra \SI^n$ of a chart $U \ra\rp^n$  
so  that $\phi(U)$ is a nonempty intersection of a closed 
$n$-hemisphere with a simply convex open  set. (By  an  $n$-hemisphere, 
we mean  a closed hemisphere unless we mention otherwise.) $\delta M$ 
is said to be {\em convex} if there exists an open neighborhood $U$ 
and a chart $\phi$ for each point of $\delta M$ so that $\phi(U)$ 
is a convex domain in $\SI^n$. $\delta M$ is said to be {\em concave} 
if there exists a chart $(U, \phi)$ for each point of $\delta M$ so 
that $\phi(U)$ is the complement of a convex open set in an open 
simply convex subset of $\SI^n$. Note that if $M$ has totally geodesic 
boundary, then so do all of its covers. The same is true for convexity 
and concavity of boundary.

\begin{lem}\label{lem:bdtangent} 
Suppose that a connected totally geodesic $(n-1)$-submanifold $S$ of 
$M$ of codimension $\geq 1$ intersects $\delta M$ in its interior 
point. Then $S \subset \delta M$. 
\end{lem}
\begin{pf} 
The intersection point must be a tangential intersection point. Since 
$\delta \tilde M$ is a closed subset of $\tilde M$, the set of 
intersection of $S$ and $\delta \tilde M$ is an open and closed subset 
of $S$. Hence it must be $S$. 
\end{pf}

\begin{rem}
If $\delta M$ is assume to be convex, the conclusion holds also. 
This was done in \cite{cdcr1} in dimension $2$. The proof for 
the convex boundary case is the same as the dimension $2$.
\end{rem}

\begin{rem}\label{rem:rpmaps2}
Given any two real projective maps $f_1, f_2: N \ra \rp^n$ on a real 
projective manifold $N$, they 
differ by an element of $\PGL(n+1, \bR)$, i.e., $f_2 = \zeta \circ f_1$ 
for a projective automorphism $\zeta$ as they are charts restricted to 
an open set, and they must satisfy the equation there, and by analyticity 
everywhere. Given two real projective automorphisms $f_1, f_2: N 
\ra \SI^n$, we have that $p\circ f_1 = \zeta \circ p \circ f_2$ for 
$\zeta$ in $\PGL(n+1, \bR)$. By equation \ref{eqn:autolift}, there 
exists an element $\zeta'$ of $\Aut(\SI^n)$ so that $p \circ \zeta' 
= \zeta \circ p$ where $\zeta'$ and $A_{\SI^n}\circ \zeta'$ are the only 
automorphisms satisfying the equation. This means that $p \circ f_1 = 
p \circ \zeta'\circ f_2$, and hence it follows easily that $f_1 = \zeta' 
\circ f_2$ or $f_1 = A_{\SI^n}\circ \zeta' \circ f_2$ by analyticity of 
developing maps. Hence, any two real projective maps $f_1, f_2: N \ra \SI^n$
differ by an element of $\Aut(\SI^n)$. 
\end{rem}

We agree to lift  our developing map $\dev$ to the  standard sphere 
$\SI^n$, the  double  cover of $\rp^n$, where we denote the lift by $\dev'$.  
Then for any deck transformation $\vth$ of $\tilde M$, we have 
$\dev' \circ \vth = h'(\vth) \circ \dev'$ by the above remark. Hence 
$\vth \mapsto h'(\vth)$ is a homomorphism, and we see easily that $h'$ is 
a lift of $h$ for the covering homomorphism $\Aut(\SI^n) \ra \PGL(n+1, \bR)$.

The pair $(\dev', h')$ will from now on be denoted by $(\dev, h)$,  
and  they   satisfy   $\dev\circ \gamma = h(\gamma)\circ \dev$ for 
every $\gamma \in \pi_1(M)$, and moreover, given a real projective 
structure, $(\dev, h)$ is determined up to an action of $\vth$ of 
$\Aut(\SI^n)$ as in equation \ref{eqn:devconj} by the above remark.

The sphere $\SI^n$ has a standard metric  $\mu$ so that its projective
structure is   projectively equivalent   to  it; i.e.,   the geodesics
agree. Let  us denote  by  $\bdd$ the   distance metric  induced  from
$\mu$. From the immersion $\dev$, we  induce a Riemannian metric $\mu$
of $\tilde M$,  and let $\bdd$ denote  the induced  distance metric on
$\tilde M$.  The Cauchy completion of $(\tilde M, \bdd)$ is denoted by
$(\che M,  \bdd)$, which we say is the {\em Kuiper completion}\/ or
{\em projective completion}\/ of $\tilde M$. We define the frontier 
$\ideal{M}=\che M-\tilde M$. 

Note that $\dev$ extends  to a  distance decreasing map, which we 
denote by $\dev$ again. Since for each $\vth \in \Aut(\SI^n)$, 
$\vth$ is quasi-isometric with respect to $\bdd$, and each deck 
transformations $\vpi$ of $\tilde M$ locally mirror the metrical 
property of $h(\vpi)$, it follows that the deck transformations are 
quasi-isometric (see \cite{cdcr1}). Thus, each deck transformation of 
$\tilde M$ extends to a self-homeomorphism of $\che M$. The extended map 
will be still called a {\em deck transformation}\/ and will be denoted by
the same symbol $\vpi$ if so was the original deck transformation denoted. 
Finally, the equation $\dev \circ \vth = h(\vth)\circ \dev$ still holds 
for each deck transformation $\vth$.

The kernel $K$ of $h:\pi_1(M) \ra \Aut(\SI^n)$ is well-defined since 
$h$ is well-defined up to conjugation. Since $\dev \circ \vth = \dev$ 
for $\vth \in K$, we see that $\dev$ induces a well-defined immersion
$\dev': \tilde M /K \ra \SI^n$. We say that $\tilde M/K$ 
the {\em holonomy}\/ cover of $M$, and denote it by $M_h$. 
We identify $K$ with $\pi_1(M_h)$.
Since any real projective map $f: M_h \ra \SI^n$ equals $\vth \circ \dev'$ 
for $\vth$ in $\Aut(\SI^n)$ by Remark \ref{rem:rpmaps2}, it follows that 
$\dev \circ \vpi$ equals $h'(\vpi)\circ \dev$ for each deck transformation 
$\vpi \in \pi_1(M)/\pi_1(M_h)$. Thus, $\vpi \mapsto h'(\vpi)$ is 
a homomorphism $h': \pi_1(M)/\pi_1(M_h) \ra \Aut(\SI^n)$, which is 
easily seen to equal $h' = h\circ \Pi$ for the quotient homomorphism 
$\Pi:\pi_1(M) \ra \pi_1(M)/\pi_1(M_h)$.

Moreover, by Remark \ref{rem:rpmaps2}, $(\dev', h')$ is determined up 
to an action of $\vth$ in $\Aut(\SI^n)$ as in equation \ref{eqn:devconj}. 
Conversely, such a pair $(f, k)$ where $f:M_h \ra \SI^n$ equivariant with 
respect to the homomorphism $k: \pi_1(M)/\pi_1(M_h) \ra \Aut(\SI^n)$ 
determines a real projective structure on $M$. From now on, we will denote 
$(\dev', h')$ by $(\dev, h)$, and call them a {\em development pair}.

Given $\dev$, we may pull-back $\mu$, and complete the distance metric 
$\bdd$ to obtain $\che M_h$, the completion of $M_h$, which is again 
called a {\em Kuiper completion}. We define the frontier $\hideal{M}$ to be 
$\che M_h-M_h$. As before the developing map $\dev$ extends to a distance 
decreasing map, again denoted by $\dev$, and each deck transformation 
extends to a self-homeomorphism $\che M_h \ra \che M_h$, which we call 
a {\em deck transformation} still. Finally, the equation $\dev\circ\vth 
= h(\vth) \circ \dev$ still holds for each deck transformation $\vth$.

\begin{figure}[h]
\centerline{\epsfxsize=3in \epsfbox{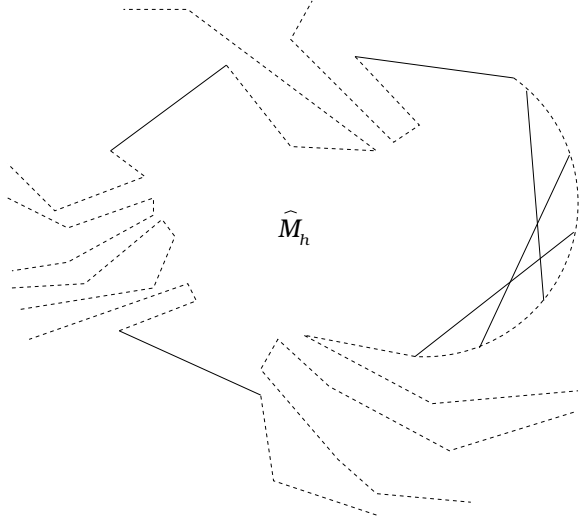}}
\caption{\label{fig:mfl} A figure of $\che M_h$. 
The thick dark lines indicate $\partial M_h$ and the dotted lines 
the ideal boundary $\hideal{M}$, and $2$-crescents in them in 
the right. They can have as many ``pods'' and what looks like 
``overlapping''. Such pictures happen if we graft annuli
into convex surfaces.}
\end{figure}
\typeout{<<pfm.eps>>}

A subset $A$ of $\che M$ is a {\sl convex segment} if $\dev| A$ is 
an imbedding onto a convex segment in $\SI^n$. $M$ is {\sl convex} 
if given two points of the universal cover $\tilde M$, there exists 
a convex segment in $\tilde M$ connecting these two points. A subset 
$A$ of $\che M$ is {\sl convex} if given points $x$ and $y$ of $A$, 
$A$ includes a convex segment containing $x$ and $y$. We say that $A$ 
is a {\sl tame} subset if it is a convex subset of $\tilde M$ or 
a convex subset of a compact convex subset of $\che M$. If $A$ is tame, 
then $\dev| A$ is an imbedding onto $\dev(A)$ and $\dev|\clo(A)$ for 
the closure $\clo(A)$ of $A$ onto a compact convex set $\clo(\dev(A))$.
The {\em interior} $A^o$ of $A$ is defined to be the set corresponding
to $\clo(\dev(A))^o$ and the boundary $\partial A$ the subset of $\clo(A)$ 
corresponding to $\partial \clo(\dev(A))$. Note that $\partial A$ may 
not equal the manifold boundary $\delta A$ if $A$ has a (topological) 
manifold structure. But if $A$ is a compact convex set, then $\dev(A)$ is 
a manifold by Theorem \ref{thm:topconv}, i.e., a sphere or a ball, 
and $\partial A$ has to equal $\delta A$. In this case, we shall use 
$\delta A$ over $\partial A$.

\begin{defn}\label{defn:i-ball}
An {\em $i$-ball}\/ $A$ in $\che M$ is a compact subset of $\che M$ 
such that $\dev|A$ is a homeomorphism to an $i$-ball (not necessarily 
convex) in a great $i$-sphere and its manifold interior $A^o$ is 
a subset of $\tilde M$. A {\em convex $i$-ball}\/ is an $i$-ball 
that is convex.
\end{defn} 
Note that a tame set in $\che M$ which is homeomorphic to an $i$-ball
is not necessarily an $i$-ball in this sense; that is, its interior
may not be a subset of $\tilde M$. We will say it is 
a {\em tame topological $i$-ball} but not $i$-ball or convex $i$-ball.

We define the terms convex segments, convex subset, tame subset, $i$-ball 
and convex $i$-ball in $\che M_h$ in the same manner as for $\che M$.

{\em We will from now on will be working on $\che M_h$ only\/{\rm ;} 
however, all of the materials in this section will work for $\che M$ 
as well, and much of the materials in the remaining section will work 
also\/{\rm ;} however, we will not say explicitly as the readers can 
easily figure out these details.}

\begin{defn}\label{defn:crsc}
An $n$-ball $A$ of $\che M_h$ is said to be an {\em $n$-bihedron}\/ 
if $\dev|A$ is a homeomorphism onto an $n$-bihedron. The $n$-bihedron  
is  bounded   by two $(n-1)$-dimensional hemispheres; the corresponding 
subsets of $A$ are  said to be the {\em sides}\/ of $A$ 
(see \cite{iconv:94}\/). 
An $n$-ball  $A$ of $\che M_h$ is said to be an {\sl $n$-hemisphere} if 
$\dev| A$ is a homeomorphism onto an $n$-hemisphere in $\SI^n$. 

A bihedron is said to be an {\sl $n$-crescent}\/ if one of its side is a
subset of $\hideal{M}$ and the other side is  not. An $n$-hemisphere is
said to  be  an  {\sl  $n$-crescent}\/  if a    subset in the   boundary
corresponding to an  $(n-1)$-hemisphere  under $\dev$ is a  subset  of
$\hideal{M}$ and the boundary itself is not a subset of $\hideal{M}$. 
\end{defn}

In this paper, we will often omit `$n$-' from $n$-bihedron. To distinguish, 
a {\em bihedral $n$-crescent}\/ is an $n$-crescent that is a bihedron, 
and a {\em hemispheric $n$-crescent}\/ is an $n$-crescent that is otherwise. 

Contrast to Definition \ref{defn:i-ball}, we define an {\em $m$-bihedron}
for $1 \leq m < n$, to be only a tame topological $m$-ball whose image 
under $\dev$ is an $m$-bihedron in a great $m$-sphere in $\SI^n$, and 
an {\em $m$-hemisphere}, $0 \leq m \leq n-1$, to be one whose image 
under $\dev$ is an $m$-hemisphere. So, we do not necessarily have 
$A^o \subset M_h$ when $A$ is one of these.

\begin{exmp}\label{exmp:ncrsc}
Let us give two trivial examples of real projective $n$-manifolds to 
demonstrate $n$-crescents (see \cite{cdcr1} for $2$-dimensional examples).

Let $\bR^n$ be an affine patch of $\SI^n$ with standard affine coordinates 
$x_1, x_2, \dots, x_n$ and $O$ the origin. Consider $\bR^n - \{O\}$ quotient 
out by the group $<g>$ where $g: x \ra 2x$ for $x \in \bR^n -\{0\}$.
Then the quotient is a real projective manifold diffeomorphic to 
$\SI^{n-1} \times \SI^1$. Denote the manifold by $N$, and we see that 
$N_h$ can be identified with $\bR^n - \{O\}$. Thus, $\che N_h$ equals the 
closure of $\bR^n$ in $\SI^n$; that is, $\che N_h$ equals 
an $n$-hemisphere $H$, 
and $\hideal{N}$ is the union of $\{O\}$ and the boundary great sphere 
$\SI^{n-1}$ of $H$. Moreover, the closure of the set given by $x_1 + x_2 + 
\cdots + x_n > 0$ in $H$ is an $n$-bihedron and one of its face is 
included in $\SI^{n-1}$. Hence, it is an $n$-crescent. 

Let $H_1$ be the open half-space given by $x_1 > 0$, and $l$ the line 
$x_2 = \cdots = x_n = 0$. Let $g_1$ be the real projective transformation 
given by $(x_1, x_2, \dots, x_n) \mapsto (2x_1, x_2, \dots, x_n)$ 
and $g_2$ that given by 
\[(x_1, x_2, \dots, x_n) \mapsto (x_1, 2x_2, \dots, 2x_n).\] 
Then the quotient manifold $L$ of $H_1 - l$ by the commutative group 
generated by $g_1$ and $g_2$ is diffeomorphic to 
$\SI^{n-2} \times \SI^1 \times \SI^1$, and we may identify 
its holonomy cover $L_h$ with $H_1 - l$ and $\che L_h$ with 
the closure $\clo(H_1)$ of $H_1$ in $\SI^n$. 
Clearly, $\clo(H_1)$ is an $n$-bihedron bounded by an $(n-1)$-hemisphere 
that is the closure of the hyperplane given by $x_1 = 0$ and an 
$(n-1)$-hemisphere in the boundary of the affine patch $\bR^n$. 
Therefore, $\hideal{L}$ is the union of $H_1 \cap l$ and two 
$(n-1)$-hemispheres that form the boundary of $\clo(H_1)$. 
$\clo(H_1)$ is not an $n$-crescent since $\clo(H_1)^o \cap \hideal{L} 
\supset l \cap H_1^o \ne \emp$. In fact, $\clo(H_1)$ includes 
no $n$-crescents.
\end{exmp}

Let $R$ be an $n$-crescent. If $R$ is an $n$-bihedron, then we define 
$\alpha_R$ to be the interior of the side of $R$ in $\hideal{M}$ and $\nu_R$ 
the other side. If $R$ is an $n$-hemisphere, then we define $\alpha_R$ to 
be the union of the interiors of all $(n-1)$-hemispheres in 
$\partial R \cap \hideal{M}$ and define $\nu_R$ the complement of $\alpha_R$ 
in $\partial R$. Clearly, $\nu_R$ is a tame topological $(n-1)$-ball.

Let us now discuss about how two convex sets may meet. Let $F_1$ 
and $F_2$ be two convex $i$-, $j$-balls in $\che M_h$ respectively. 
We say that $F_1$ and $F_2$ overlap if $F_1^o \cap F_2 \ne \emp$, which 
is equivalent to $F_1 \cap F_2^o\ne \emp$ or $F_1^o \cap F_2^o \ne \emp$. 

\begin{thm}\label{thm:ovlconv} 
If $F_1$ and $F_2$ overlap, then $\dev|F_1 \cup F_2$ is an imbedding 
onto $\dev(F_1) \cup \dev(F_2)$ and $\dev| F_1 \cap F_2$ onto 
$\dev(F_1) \cap \dev(F_2)$. Moreover, if $F_1$ and $F_2$ 
are $n$-balls, then $F_1 \cup F_2$ is an $n$-ball, and 
$F_1 \cap F_2$ is a convex $n$-ball.
\end{thm}
\begin{pf} 
The proof is a direct generalization of that of Theorem 1.7
of \cite{cdcr1}. 
\end{pf}

The above theorem follows from
\begin{prop}\label{prop:ovlballs}
Let $A$ be a $k$-ball in $\che M_h$ and $B$ an $l$-ball. Suppose 
that $A^o \cap B^o \ne \emp$, $\dev(A)\cap \dev(B)$ is a compact 
manifold in $\SI^n$ with interior equal to $\dev(A^o) \cap \dev(B^o)$ 
and $\dev(A^o) \cap \dev(B^o)$ is arcwise-connected. Then 
$\dev|A \cup B$ is a homeomorphism onto $\dev(A) \cup \dev(B)$.
\end{prop}
\begin{pf} 
This follows as in its affine version Lemma 6 in \cite{uaf:97}.
\end{pf}

In the following, we describe a useful geometric situation modelled on 
``dipping a bread into a bowl of milk''. Let $D$ be a convex $n$-ball in 
$\che M_h$ such that $\partial D$ includes a tame subset $\alpha$ 
homeomorphic to an $(n-1)$-ball. We say that a convex $n$-ball $F$ 
is {\sl dipped into} $(D, \alpha)$ if the following statements hold: 
\begin{itemize}
\item $D$ and $F$ overlap.
\item $F \cap \alpha$ is a convex $(n-1)$-ball $\beta$ 
with $\delta \beta \subset \delta F$ and $\beta^o\subset F^o$. 
\item $F - \beta$ has two convex components $O_1$ and $O_2$ such that 
$\clo(O_1) = O_1 \cup \beta = F - O_2$ and
$\clo(O_2) = O_2 \cup \beta = F - O_1$. 
\item $F \cap D$ is equal to $\clo(O_1)$ or $\clo(O_2)$.
\end{itemize}
(The second item sometimes is crucial in this paper.)
We say that $F$ is {\sl dipped into $(D, \alpha)$ nicely} 
if the following statements hold:
\begin{itemize}
\item $F$ is dipped into $(D, \alpha)$.
\item $F\cap D^o$ is identical with $O_1$ and $O_2$.
\item $\delta(F\cap D) = \beta \cup \xi$ for a topological $(n-1)$-ball 
$\xi$, not necessarily convex or tame, in the topological boundary $\Bd F$ 
of $F$ in $\che M_h$ where $\beta \cap \xi = \delta \beta$. 
\end{itemize}
As a consequence, we have $\delta \beta \subset \Bd F$. (As above
this is a crucial point.)
(The nice dipping occurs when the bread does not touch the bowl.)


\begin{figure}[h]
\centerline{\epsfxsize=3.7in \epsfbox{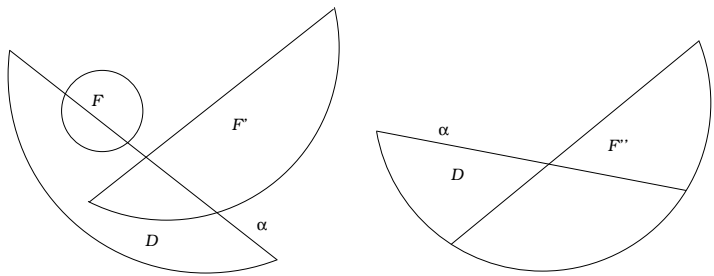}}
\caption{\label{fig:dipping} 
Various examples of dipping intersections. Loosely speaking $\alpha$ 
plays the role of the milk surface, $F, F',$ and $F''$ the breads, 
and $D^o$ the milk. The left one indicates nice dippings, and 
the right one not a nice one.}
\end{figure}
\typeout{<<pdip.eps>>}

The direct generalization of Corollary 1.9 of \cite{cdcr1} gives us:
\begin{cor}\label{cor:dipping} Suppose that $F$ and $D$ overlap, 
and $F^o \cap (\delta D - \alpha^o) = \emp$. 
Assume the following two equivalent conditions\/{\rm :} 
\[ \bullet \quad F^o \cap \alpha \ne \emp,  
\qquad \bullet \quad F \not\subset D. \] 
Then $F$ is dipped into $(D, \alpha)$. 
If $F \cap (\delta D - \alpha^o) = \emp$ furthermore{\rm ,} 
then $F$ is dipped into $(D, \alpha)$ nicely. \qed
\end{cor}

\begin{exmp}\label{exmp:dipexmp}
In Example \ref{exmp:ncrsc}, choose a compact convex ball $B$ in 
$\bR^n -\{O\} = N_h$ intersecting $R$ in its interior but not 
included in $R$. Then $B$ dips into $(R, P)$ nicely where $P$ 
is the closure of the plane given by $x_1 + \cdots + x_n = 1$. 
Also let $S$ be the closure of the half plane given by $x_1 > 0$. 
Then $S$ dips into $(R, P)$ but not nicely.

Consider the closure of the set in $\che N_h$ given by $0 < x_1 < 1$ 
and that of the set $0 < x_2 < 1$. Then these two sets do not dip 
into each other for any choice of $(n-1)$-balls in their respective 
boundaries to play the role of $\alpha$.
\end{exmp} 

Since $\dev$ restricted to a small open sets are charts, and the 
boundary of $M_h$ is convex, each point $x$ of $M_h$
has a compact ball-neighborhood $B(x)$ so that $\dev| B(x)$ is 
an imbedding onto a compact convex ball in $\SI^n$ (see Section 
1.11 of \cite{cdcr1}). $\dev(B(x))$ can be assumed to be 
a $\bdd$-ball with center $\dev(x)$ and radius $\eps> 0$ intersected 
with an $n$-hemisphere $H$ so that $\delta M_h \cap B(x)$ corresponds 
to $\partial H \cap \dev(B(x))$. Of course, $\delta M_h \cap B(x)$ or
$\partial H \cap \dev(B(x))$ may be empty. We say that such $B(x)$ 
is an {\em $\eps$-tiny ball} of $x$ and $\eps$ the {\sl $\bdd$-radius} 
of $B(x)$. Thus, for an $\eps$-tiny ball $B(x)$, $\delta M_h \cap B(x)$ 
is a compact convex $(n-1)$-ball, and the topological boundary 
$\Bd B(x)$ equals the closure of $\delta B(x)$ removed with this 
$(n-1)$-ball.

\begin{lem}\label{lem:tinydip}
If $B(x)$ and an $n$-crescent $R$ overlap, then either $B(x)$ 
is a subset of $R$ or $B(x)$ is dipped into $(R, \nu_R)$ nicely. 
\end{lem}
\begin{pf} 
Since $\clo(\alpha_R) \subset \hideal{M}$ and $B(x) \subset M_h$, 
Corollary \ref{cor:dipping} implies the conclusion.
\end{pf} 

\vfill
\break

\section{$(n-1)$-convexity}
\label{sec:n-1conv}
In this section, we introduce $m$-convexity. Then we state
Theorem \ref{thm:n-1conv} central to this section, which relates
the failure of $(n-1)$-convexity with an existence of $n$-crescents, or
half-spaces. The proof of theorem is similar to what is in 
Section 5 in \cite{cdcr1}. We first choose a sequence of points 
$q_i$ converging to a point $x$ in $F_1 \cap \hideal{M}$. Then we 
pull back $q_i$ to points $p_i$ in the closure of a fundamental 
domain by a deck transformation $\vth_i^{-1}$. Then analogously 
to \cite{cdcr1}, we show that $T_i = \vth_i^{-1}(T)$ ``converges to'' 
a nondegenerate convex $n$-ball. Showing that $\dev(T_i)$ converges 
to an $n$-bihedron or an $n$-hemisphere is much more complicated 
than in \cite{cdcr1}. The idea of the proof is to show that the 
sequence of the images under $\vth_i$ of the $\eps$-$(n-1)$-$\bdd$-balls 
in $\vth^{-1}(F_1)$ with center $p_i$ have to degenerate to a point 
when $x$ is chosen specially. So when pulled back by $\vth_i^{-1}$, 
the balls become standard ones again, and $F_1$ must blow up to be 
an $(n-1)$-hemisphere under $\vth_i^{-1}$.

An $m$-simplex $T$ in $\che  M$ is a tame  subset of $\che M_h$ such that
$\dev|T$ is an imbedding onto an affine $m$-simplex in an affine patch
in $\SI^n$. 

\begin{defn} 
We say that $M$ is {\sl  $m$-convex\/}, $0 < m  < n$, if the following
holds.  If $T \subset \che M_h$ be an $(m +1)$-simplex with sides $F_1,
F_2, \dots, F_{m+2}$ such that $T^o \cup F_2 \cup \cdots \cup F_{m+2}$ 
does not meet $\hideal{M}$, then $T$ is a subset of $M_h$.
\end{defn}

\begin{thm}\label{thm:equivconv}
Let $T$ be an affine $(m+1)$-simplex in an affine space with faces 
$F_1, F_2, \dots,$ \break $F_{m+2}$. 
The following are equivalent\/{\rm :}
\begin{enumerate}
\item[(a)] $M$ is $m$-convex.
\item[(b)] Any {\rm (}nonsingular\/{\rm )} real projective map $f$ from 
$T^o \cup F_1 \cup F_2 \cup \cdots \cup F_{m+2}$ to $M$ extends to one 
from $T$. 
\item[(c)] a cover of $M$ is $m$-convex.
\end{enumerate}
\end{thm}
\begin{pf}
The proof of the equivalence of (a) and (b) is the same as the affine 
version Lemma 1 in \cite{uaf:97}. The equivalence of (b) and (c) 
follows from the fact that a real projective map to $M$ always 
lifts to its cover.
\end{pf}

\begin{thm}\label{thm:equiv2conv}
$M$ is not convex if and only if there exists an $(m+1)$-simplex with 
a side $F_1$ such that $T \cap \hideal{M} =  F_1^o \cap \hideal{M} 
\ne \emp$.
\end{thm}
\begin{pf} 
The proof is same as Lemma 3 in \cite{uaf:97}.
\end{pf}

\begin{figure}[h]
\centerline{\epsfxsize=3.5in \epsfbox{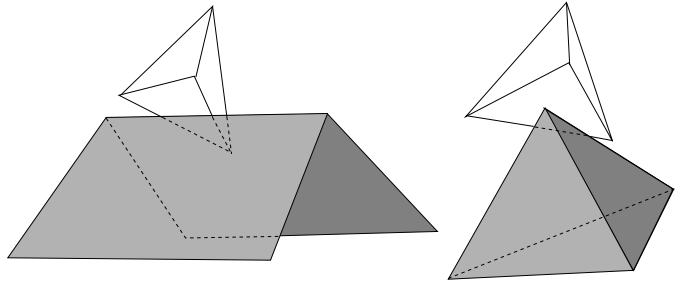}}
\caption{\label{fig:mconv} The tetrahedron in the left fails to 
detect non-$2$-convexity but the right one is detecting 
non-$2$-convexity.}
\end{figure}
\typeout{<<pfm-conv.eps>>}

\begin{rem}\label{rem:ijconv}
It is easy to see that $i$-convexity implies $j$-convexity 
whenever $i\leq j < n$. (See Remark 2 in \cite{uaf:97}. The proof 
is the same.)
\end{rem}


\begin{thm}\label{thm:1conv}
The following are equivalent\/{\rm :}
$\bullet$ $M$ is $1$-convex. $\bullet$ $M$ is convex.
$\bullet$ $M$ is real projectively isomorphic to a quotient of 
a convex domain in $\SI^n$.
\end{thm}

The proof is similar to Lemma 8 in \cite{uaf:97}. Since there are 
minor differences between affine and real projective manifolds, we 
will prove this theorem in Appendix A. 
 
Let us give examples of real projective $n$-manifolds one of which 
is not $(n-1)$-convex and the other $(n-1)$-convex.

As in Figure \ref{fig:mconv}, $\bR^3$ removed with
a complete affine line or a closed wedge, i.e. a set defined 
by the intersection of two half-spaces with non-parallel boundary 
planes is obviously $2$-convex. But $\bR^3$ removed with 
a discrete set of points or a convex cone defined as the intersection 
of three half-spaces with boundary planes in general position is not 
$2$-convex.

We recall Example \ref{exmp:ncrsc}.
Let $\bR^n$ be an affine patch of $\SI^n$ with standard affine 
coordinates $x_1, x_2, \dots, x_n$ and $O$ the origin. Consider 
$\bR^n - \{O\}$ quotient out by the group $<g>$ where 
$g: x \ra 2x$ for $x \in \bR^n -\{0\}$. Then the quotient is 
a real projective manifold diffeomorphic to $\SI^{n-1} \times \SI^1$. 
Denote the manifold by $M$, and we see that $M_h$ can be identified 
with $\bR^n - \{0\}$. Thus, $\che M_h$ equals the closure of 
$\bR^n$ in $\SI^n$; that is, $\che M_h$ equals an $n$-hemisphere $H$, 
and $\hideal{M}$ is the union of $\{O\}$ and the boundary great sphere 
$\SI^{n-1}$ of $H$. Consider an $n$-simplex $T$ in $\bR^n$ given by 
$x_i \leq 1$ for every $i$ and $x_1 + x_2 + \cdots + x_n \geq 0$. 
Then the face of $T$ corresponding to $x_1 + x_2 + \cdots + x_n = 0$ 
contains the ideal point $O$ in its interior. Therefore, $M$ is 
not $(n-1)$-convex. Moreover, the closure of the set given by 
$x_1 + x_2 + \cdots + x_n > 0$ in $\che M_h = H$ is an $n$-bihedron 
and one of its face is included in $\SI^{n-1}$. Hence, it is 
an $n$-crescent.

Let $H_1$ be the open half-space given by $x_1 > 0$, and $l$ 
the line $x_2 = \cdots = x_n = 0$. Let $g_1$ be the real projective 
transformation given by $(x_1, x_2, \dots, x_n) \mapsto 
(2x_1, x_2, \dots, x_n)$ and $g_2$ that given by 
$(x_1, x_2, \dots, x_n) \mapsto (x_1, 2x_2, \dots, 2x_n)$. 
Then the quotient manifold $M$ of $H_1 - l$ by the commutative group 
generated by $g_1$ and $g_2$ is diffeomorphic to $\SI^{n-2} \times 
\SI^1 \times \SI^1$, and we may identify $M_h$ with $H_1 - l$ 
and $\che M_h$ with the closure $\clo(H_1)$ of $H_1$ in $\SI^n$. 
Clearly, $\clo(H_1)$ is an $n$-bihedron bounded by 
an $(n-1)$-hemisphere that is the closure of the hyperplane given 
by $x_1 = 0$ and an $(n-1)$-hemisphere in the boundary of 
the affine patch $\bR^n$. Therefore, $\hideal{M}$ is the union of 
$H_1 \cap l$ and two $(n-1)$-hemispheres that form the boundary of 
$\clo(H_1)$. The intersection of an $n$-simplex $T$ in $\SI^n$ with 
the boundary $(n-1)$-hemispheres or $l$ is not a subset of the interior 
of a face of $T$. It follows from this that $M$ is $(n-1)$-convex.  

The main purpose of this section is to prove the following 
principal theorem:
\begin{thm}\label{thm:n-1conv} 
Suppose that a compact real projective manifold $M$ with empty 
or totally geodesic boundary is  not  $(n-1)$-convex.  
Then  the completion $\che M_h$ of the holonomy cover $M_h$ 
includes an $n$-crescent. 
\end{thm} 
We may actually replace the word ``total geodesic'' with ``convex''
and the proof is same step by step. However, we need this result 
at only one point of the paper so we do not state it. We can also 
show that the completion $\che M$ of the universal cover $\tilde M$ 
also includes an $n$-crescent. 
The proof is identical with $\che M$ replacing $\che M_h$. 
Another way to do this is of course as follows: once we obtain 
an $n$-crescent in $\che M_h$ we may lift it to one in $\che M$ 
but we omit showing how this can be done.

\begin{rem}\label{rem:notbihem}
As $M$ is $(n-1)$-convex, we may assume that $\tilde M$ or $M_h$ 
is not projectively diffeomorphic to an open $n$-bihedron or 
an open $n$-hemisphere. This follows since if otherwise, $\tilde M$ 
is convex and hence $(n-1)$-convex. (We will need this weaker 
statement later.) 
\end{rem}

A point $x$ of a convex subset $A$ of $\SI^n$ is said to be 
{\em exposed} if there exists a supporting great $(n-1)$-sphere $H$ 
at $x$ such that $H \cap A = \{x\}$ (see Section \ref{sec:convSI} 
and Berger \cite[p. 361]{Berger:87}).

To prove Theorem \ref{thm:n-1conv}, we follow Section 5 of Choi 
\cite{cdcr1}:  Since $M$ is not $(n-1)$-convex, $\che M_h$ includes 
an $n$-simplex $T$ with a face $F_1$ such that $T \cap \hideal{M} 
= F_1^o \cap \hideal{M} \ne \emp$ by Theorem \ref{thm:equiv2conv}, 
where $\dev|T: T \ra \dev(T)$ is an imbedding onto the $n$-simplex 
$\dev(T)$. Let $K$ be the convex hull of $\dev(F_1 \cap \hideal{M})$ 
in $\dev(F_1)^o$, which is simply convex as $\dev(F_1)$ is simply 
convex. 

As $K$ is simply convex, we see that $K$ can be considered as a bounded 
convex subset of an affine patch, i.e., an open $n$-hemisphere. We see 
easily that $K$ has an exposed point in the affine sense in the open 
hemisphere, which is easily seen to be an exposed point in our sense
as a hyperplane in the affine patch is the intersection of a hypersphere 
with the affine patch. 

Let $x'$ be an exposed point of $K$. Then 
$x' \in \dev(F_1 \cap \hideal{M})$, 
and there exists a line $s'$ in the complement of $K$ in $\dev(F_1)^o$ 
ending at $x'$. Let $x$ and $s$ be the inverse images of $x'$ and $s'$ 
in $F_1^o$ respectively. 

Let $F_i$ for $i=2, \dots, n+1$ denote the faces of $T$ other than $F_1$.
Let $v_i$ for each $i$, $i= 1, \dots, n+1$, denote the vertex of $T$ 
opposite $F_i$. Let us choose a monotone sequence of points $q_i$ on $s$ 
converging to $x$ with respect to $\bdd$. 

Choose a fundamental domain $F$ in $M_h$ such that for every point $t$ 
of $F$, there exists a $2\eps$-tiny ball of $t$ in $M_h$ for a positive
constant $\eps$ independent of $t$. We assume $\eps \leq \pi/8$ for 
convenience. Let us denote by $F_{2\eps}$ the closure of the 
$2\eps$-$\bdd$-neighborhood of $F$, and $F_\eps$ that of the 
$\eps$-$\bdd$-neighborhood of $F$.

For each natural number $i$, we choose a deck transformation $\vth_i$ 
and a point $p_i$ of $F$ so that $\vth_i(p_i) = q_i$. We let 
$v_{j, i}, F_{j, i},$ and $T_i$, $i = 1, 2, \dots$, $j = 1, \dots, n+1$, 
denote the images under $\vth_i^{-1}$ of $v_j, F_j,$ and $T$ 
respectively. Let $n_i$ denote the outer-normal vector to $F_{1, i}$ 
at $p_i$ with respect to the spherical Riemannian metric $\mu$ of $M_h$.

We choose subsequences so that each sequence consisting of 
\[\dev(v_{j, i}), \dev(F_{j, i}), \dev(T_i), n_i, \mbox{ and } p_i\] 
converge geometrically with respect to $\bdd$ for each $j$, 
$j= 1, \dots, n+1$ respectively. Since $p_i \in F$ for each $i$, 
the limit $p$ of the sequence of $p_i$ belongs to $\clo(F)$. We choose 
an $\eps$-tiny ball $B(p)$ of $p$. We \awlg that $p_i$ belongs to 
the interior $\inte B(p)$ of $B(p)$. Since the action of the deck 
transformation group is properly discontinuous and 
$F_{j, i} = \vth_i^{-1}(F_j)$ for a compact set $F_j$, there exists 
a natural number $N$ such that 
\begin{equation}\label{eqn:disjoint} 
F_{2\eps} \cap F_{j, i} = \emp \hbox{ for each } j, i, j> 1, i > N;
\end{equation} 
so $B(p) \cap F_{j, i} = \emp$ for $j > 1$. (This corresponds to 
Lemma 5.4 in \cite{cdcr1}.) Hence, $B(p) \subset T_i$ or $B(p)$ 
dips into $(T_i, F_{1, i})$ for each $i$, $i > N$, by Corollary 
\ref{cor:dipping}.

\begin{figure}[h]
\centerline{\epsfxsize=3.5in \epsfbox{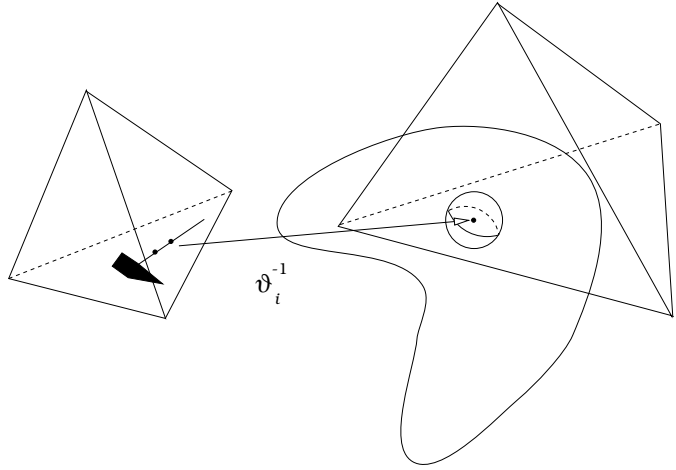}}
\caption{\label{fig:pull1} The pull-back process }
\end{figure}
\typeout{<<pfpull.eps>>}

Lemma 3 of the appendix of \cite{cdcr1} holds for manifolds of 
higher dimensions as the dimension of the sphere $\SI^2$ do not matter 
in the proof. Thus, there exists an integer $N_1$, $N_1 > N$, such 
that $T_i$ includes a common open ball for $i > N_1$. Let $T_\infty$ 
be the limit of $\dev(T_i)$. Since $\dev(T_i)$ includes a common ball 
for $i > N_1$, it follows that $T_\infty$ is a closed convex $n$-ball 
in $\SI^n$ (see Section 2 of the appendix of \cite{cdcr1}). 

Let $F_{j, \infty}$ denote the limit of $\dev(F_{j, i})$. Then 
$\bigcup_{j=1}^{n+1} F_{j, \infty}$ is the boundary $\partial T_\infty$
by Theorem \ref{thm:gconv}. 

Theorem 4 of the appendix of \cite{cdcr1} also holds for higher 
dimensional real projective manifolds as there are no dimensional 
assumptions that matter. Thus, $\che M_h$ includes a convex 
$n$-ball $T^u$ and convex sets $F^u_j$ such that $\dev$ restricted 
to them are imbeddings onto $T_\infty$ and $F_{j, \infty}$ 
respectively. We have $F^u_j \subset \hideal{M}$ for $j \geq 2$ 
from the same theorem since $F_{j, i}$ is ideal. 

We will prove below that $F_{1, \infty}$ is an $(n-1)$-hemisphere.
It follows that $T_\infty$ is a compact convex $n$-ball in $\SI^n$ 
including the $(n-1)$-hemisphere $F_{1, \infty}$ in its boundary 
$\delta T_\infty$. By Lemma \ref{lem:classf}, $T_\infty$ is an 
$n$-bihedron or an $n$-hemisphere. As $\bigcup_{j \geq 2} F^u_j$ 
is a subset of $\hideal{M}$, if $F^u_1 \subset \hideal{M}$, 
then $\che M_h = T^u$ and $M_h$ equals the interior of $T^u$ by 
the following lemma \ref{lem:include}. Since $M_h$ is not 
projectively diffeomorphic to an open $n$-bihedron or an open 
$n$-hemisphere by premise, $F^u_1$ is not a subset of $\hideal{M}$. 
Since $T^u$ is bounded by $F^u_1$ and 
$F^u_2 \cup \cdots \cup F^u_{n+1} \subset \hideal{M}$, 
it follows that $T$ is an $n$-crescent. This completes the proof of 
Theorem \ref{thm:n-1conv}.
\qed

\begin{lem}\label{lem:include} 
Suppose that $\che M_h$ includes an $n$-ball $B$ with 
$\delta B \subset \hideal{M}$. Then $M_h$ equals $B^o$. 
\end{lem}
\begin{pf}
Since $\Bd B \cap M_h \subset \delta B$ and $\delta B \subset \hideal{M}$, 
it follows that $\Bd B \cap M_h$ is empty. Hence, $M_h \subset B$.
\end{pf}

{\em We will now show that $F^u_1$ is an $(n-1)$-dimensional hemisphere}. 
This corresponds to Lemma 5.5 of \cite{cdcr1} showing that one of the side
is a segment of $\bdd$-length $\pi$. (Note that this process may 
require us to choose further subsequences of $T_i$. However, since 
$\dev(F_{1, i})$ is assumed to converge to $F_{1, \infty}$, we see that
we need to only show that a {\em subsequence} of $\dev(F_{1, i})$ converges 
to an $(n-1)$-hemisphere.)

The sequence $\dev(q_i) = h(\vth_i)\dev(p_i)$ converges to $x'$. 
Since $p_i \in F$, $M_h$ includes an $\eps$-tiny ball $B(p_i)$ 
and a $2\eps$-tiny ball $B'(p_i)$ of $p_i$. Let 
$W(p_i) = F_{1, i} \cap B(p_i)$ and $W'(p_i) = F_{1, i} \cap B'(p_i)$. 
We assume that $i > N_1$ from now on. 

We now show that $W(p_i)$ and $W(p'_i)$ are ``whole'' $(n-1)$-balls 
of $\bdd$-radius $\eps$ and $2\eps$, i.e., they map to such balls in 
$\SI^n$ under $\dev$ respectively, or they are not ``cut off'' by 
the boundary $\delta F_{1, i}$:

If $p_i \in \delta M_h$, then the component $L$ of $F_{1, i} \cap M_h$ 
containing $p_i$ is a subset of $\delta M_h$ by Lemma 
\ref{lem:bdtangent}. This component is a submanifold of $\delta M_h$ 
with boundary $\delta F_{1, i}$. Since $\delta F_{1, i}$ is a subset of 
$\bigcup_{j \geq 2} F_{j, i}$, and $B(p_i)$ is disjoint from it by 
equation \ref{eqn:disjoint}, $\delta M_h \cap B(p_i)$ is a subset of $L^o$. 
Thus, $W(p_i)$ equals the convex $(n-1)$-ball $\delta M_h \cap B(p_i)$
with boundary in $\Bd B(p_i)$ and is a $\bdd$-ball in $F_{1, i}^o$ of 
dimension $(n-1)$ of $\bdd$-radius $\eps$ and center $p_i$. It certainly
maps to an $(n-1)$-ball of $\bdd$-radius $\eps$ with center $\dev(p_i)$.

If $p_i \in M_h^o$, then since $F_{1, i}$ passes through $p_i$, 
and $F_{j, i} \cap B(p_i) = \emp$ for $j \geq 2$, it follows that 
$B(p_i)$ dips into $(T_i, F_{1, i})$ nicely by Corollary 
\ref{cor:dipping}. Thus $W(p_i)$ is an $(n-1)$-ball with boundary 
in $\Bd B(p_i)$, and an $\eps$-$\bdd$-ball in $F_{1, i}^o$ of 
dimension $(n-1)$ with center $p_i$. 

Similar reasoning shows that $W'(p_i)$ is a $2\eps$-$\bdd$-ball in 
$F_{1, i}^o$ of dimension $(n-1)$ with center $p_i$ for each $i$.

Since $\vth_i(W(p_i)) \subset F_1$, and $\dev(F_1)$ is a compact set,    
we \awlg by choosing subsequences of $\vth_i$ that the sequence of 
the subsets $\dev(\vth_i(W(p_i)))$ of $\dev(F_1)$, equal to 
$h(\vth_i)(\dev(W(p_i)))$, converges to a set $W_\infty$ containing $x'$ 
in $\dev(F_1)$. Since $\dev|T^u$ is an imbedding onto $T_\infty$, there 
exists a compact tame subset $W^u$ in $F_1$ such that $\dev$ restricted 
to $W^u$ is an imbedding onto $W_\infty$. $\vth_i(W(p_i))$ is a subjugated 
sequence of the sequence whose elements equal $T$ always. Since $W(p_i)$ 
is a subset of a compact set $F_\eps$, it follows that $\vth_i(W(p_i))$ 
is ideal (see Lemma 5.4 of \cite{cdcr1}), and $W^u \subset \hideal{M}$ 
by Lemma 4 of appendix of \cite{cdcr1}. We obtain $W^u \subset F_1 \cap 
\hideal{M}$.

\begin{figure}[h]
\centerline{\epsfxsize=4in \epsfbox{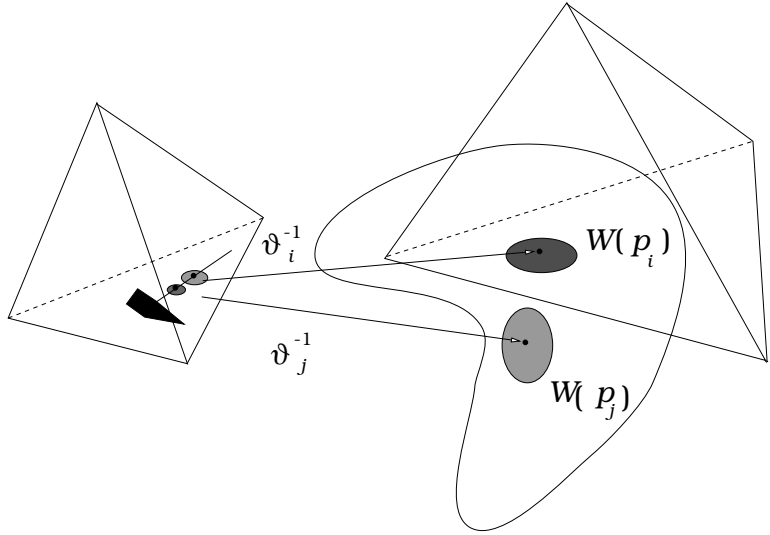}}
\caption{\label{fig:pull2} The pull-back process with $W(p_i)$s.}
\end{figure}
\typeout{<<pfpull2.eps>>}

For the proof, the fact that $x'$ is exposed will play a role:
\begin{prop} $W_\infty$ consists of the single point $x'$.
\end{prop}

Suppose not. Then as $\dev(\vth_i(W(p_i))$ does not converge to a point, 
there has to be a sequence $\{\dev(\vth_i(z_i))\}$, $z_i\in W(p_i)$, 
converging to a point $z'$ distinct from the limit $x'$ of $\{\dev(q_i)\}$.
Since we have $\dev(q_i) = \dev(\vth_i(p_i))$, we choose $s_i$ to be 
the $\bdd$-diameter of $W(p_i)$ containing $z_i$ and $p_i$, 
as a center. We obtained a sequence of segments $s_i \in W(p_i)$ 
passing through $p_i$ of $\bdd$-length $2\eps$ so that the sequence of 
segments $\dev(\vth_i(s_i))$ in $\dev(F_1)$ converges to a nontrivial 
segment $s$ containing $x'$, and $z'$, satisfying $s \subset W_\infty 
\subset \dev(\hideal{M} \cap F_1)$. 

Since $s$ is a nontrivial segment, the $\bdd$-length of 
$h(\vth_i)(\dev(s_i))$ is bounded below by a positive constant 
$\delta$ independent of $i$. Since $h(\vth_i)(\dev(s_i))$ is 
a subset of the $(n-1)$-simplex $\dev(F_1)$, which is a simply 
convex compact set, the $\bdd$-length of $h(\vth_i)(\dev(s_i))$ 
is bounded above by $\pi-\delta'$ for some small positive constant 
$\delta'$. Let $s'_i$ be the maximal segment in $W'(p_i)$ including 
$s_i$. Then the $\bdd$-length of $h(\vth_i)(\dev(s'_i))$ also 
belongs to the interval $[\delta, \pi-\delta']$.

\begin{lem}\label{lem:nshrink}
Let $\SI^1$ be a great circle and $o, s, p, q$ distinct points on 
a segment $I$ in $\SI^1$ of $\bdd$-length $< \pi$ with endpoints $o, s$
and $p$ between $o$ and $q$. Let $f_i$ be a sequence of projective maps 
$I \ra \SI^n$ so that $\bdd(f_i(o), f_i(s))$ and $\bdd(f_i(p), f_i(q))$
lie in the interval $[\eta, \pi-\eta]$ for some positive constant 
$\eta$ independent of $i$. Then all of the $\bdd$-distances between 
$f_i(o), f_i(s), f_i(p)$, and $f_i(q)$ are bounded below by a positive 
constant independent of $i$.     
\end{lem}
\begin{pf}
Recall the well-known formula for cross-ratios (see \cite{Berger:87} 
and \cite{mar:96}):
\[
[f_i(o), f_i(s); f_i(q), f_i(p)] = 
\frac{\sin(\bdd(f_i(o), f_i(q)))}{\sin(\bdd(f_i(s), f_i(q)))}
\frac{\sin(\bdd(f_i(s), f_i(p)))}{\sin(\bdd (f_i(o), f_i(p)))}. 
\]     
Suppose that $\bdd(f_i(o), f_i(p)) \ra 0$. Then since 
\begin{eqnarray*}
\bdd(f_i(o), f_i(q))&=&\bdd(f_i(o), f_i(p))+\bdd(f_i(p), f_i(q))\geq\eta, \\ 
\bdd(f_i(o), f_i(q)) &\leq & \bdd(f_i(o), f_i(s)) \leq \pi-\eta,
\end{eqnarray*}
it follows that $\sin(\bdd(f_i(o), f_i(q)))$ is bounded below and above 
by $\sin(\eta)$ and $1$ respectively. Similarly, so is 
$\sin(\bdd(f_i(s), f_i(p)))$. Therefore, the right side of the equation 
goes to $+\infty$, while the left side remains constant since $f_i$ 
is projective. This is a contradiction, and $\bdd(f_i(o), f_i(p))$ is 
bounded below by a positive constant.
 
Similarly, we can show that $\bdd(f_i(s), f_i(q))$ is bounded below by 
a positive constant. The lemma follows from these two statements.
\end{pf}  

Let $\SI^1$ be the unit circle in the plane $\bR^2$. Let $\theta$, 
$-\pi/2 < \theta < \pi/2$, denote the point of $\SI^1$ corresponding to 
the unit vector having an oriented angle of $\theta$ with $(1, 0)$ 
in $\bR^2$. Since $s_i$ and $s'_i$ are the diameters of balls of 
$\bdd$-radius $\eps$ and $2\eps$ with center $p_i$ respectively,
for the segment $[-2\eps, 2\eps]$ consisting of points 
$\theta$ satisfying $-2\eps \leq \theta \leq 2\eps$ in $\SI^1$, 
we parameterize $s'_i$ by a projective map $f_i:[-2\eps, 2\eps] \ra s'_i$, 
isometric with respect to $\bdd$, so that the endpoints of $s'_i$ 
correspond to $-2\eps$ and $2\eps$, the endpoints of $s_i$ to $-\eps$ 
and $\eps$, and $p_i$ to $0$.

Lemma \ref{lem:nshrink} applied to $k_i = h(\vth_i)\circ \dev 
\circ f_i$ shows that $\bdd(k_i(2\eps), k_i(\eps))$ and 
$\bdd(k_i(-2\eps), k_i(-\eps))$ are bounded below by a positive constant 
since $\bdd(k_i(\eps), k_i(-\eps))$ and $\bdd(k_i(2\eps), k_i(-2\eps))$ are
bounded below by a positive number $\delta$ and above by $\pi - \delta$.
Since $k_i(2\eps)$ and $k_i(-2\eps)$ are endpoints of $h(\vth_i)(\dev(s'_i))$ 
and $k_i(\eps)$ and $k_i(-\eps)$ those of $h(\vth_i)(\dev(s_i))$, 
a subsequence of $h(\vth_i)(\dev(s'_i))$ converges to a segment $s'$ 
in $\dev(F_1)$ including $s$ in its interior. Hence, $s'$ contains $x'$ 
in its interior. 

Since $s'_i$ is a subset of $F_{2\eps}$, a compact subset of $M_h$, 
it follows that the corresponding subsequence of $\vth_i(s'_i)$ is ideal 
in $F_1$. Hence $s' \subset \dev(F_1 \cap \hideal{M}) \subset K$
by Theorem 4 of Appendix of \cite{cdcr1}. Since $x'$ is not an endpoint 
of $s'$ but an interior point, this contradicts our earlier choice of 
$x'$ as an exposed point of $K$. \qed 

Since $W_\infty$ consists of a point, it follows that the sequence of the 
$\bdd$-diameter of \break $h(\vth_i)(\dev(W(p_i)))$ converges to zero, and 
the sequence converges to the singleton $\{x'\}$. 

Let us introduce a $\bdd$-isometry $g_i$, which is a real projective 
automorphism of $\SI^n$, for each $i$ so that each $g_i(\dev(W(p_i)))$ 
is a subset of the great sphere $\SI^{n-1}$ including $\dev(F_1)$, and 
hence $h(\vth_i) \circ g_i^{-1}$ acts on $\SI^{n-1}$. We \awlg that the 
$\bdd$-isometries $g_i$ converges to an isometry $g$ of $\SI^n$. Thus, 
$h(\vth_i)\circ g_i^{-1}(g_i(\dev(W(p_i))))$ converges to $x'$, and
$g_i\circ h(\vth_i)^{-1}(\dev(F_1))$ converges to $g(F_{1, \infty})$
as we assumed in the beginning of the pull-back process. By Proposition 
\ref{prop:blowup}, we see that $g(F_{1, \infty})$ is an $(n-1)$-hemisphere, 
and we are done. 

The proof of the following proposition is left to Appendix 
\ref{app:B} as the proof may distract us too much.
\begin{prop}\label{prop:blowup} 
Suppose we have a sequence of $\eps$-$\bdd$-balls $B_i$ in a real 
projective sphere $\SI^m$ for some $m \geq 1$ and a sequence of 
projective maps $\vpi_i$. Assume the following\/{\rm :}
\begin{itemize}
\item The sequence of $\bdd$-diameters of $\vpi_i(B_i)$ 
goes to zero. 
\item $\vpi_i(B_i)$ converges to a point, say $p$. 
\item For a compact $m$-ball neighborhood $L$ of $p$, 
$\vpi_i^{-1}(L)$ converges to a compact set $L_\infty$.
\end{itemize}
Then $L_\infty$ is an $m$-hemisphere.
\end{prop}

\vfill
\break

\section{The transversal intersection of $n$-crescents}
\label{sec:trcres}

From now on, we will assume that $M$ is compact and with totally 
geodesic or empty boundary. We will discuss about the transversal 
intersection of $n$-crescents, generalizing that of crescents in 
two-dimensions \cite{cdcr1}. 

First, we will show that if two hemispheric $n$-crescents overlap, 
then they are equal. For transversal intersection of two bihedral 
$n$-crescents, we will follow Section 2.6 of \cite{cdcr1}. 

Our principal assumption is that $M_h$ is not projectively diffeomorphic 
to an open $n$-hemisphere or $n$-bihedron, which will be sufficient 
for the results of this section to hold. This is equivalent to 
assuming that $\tilde M$ is not projectively diffeomorphic to these.
{\em This will be our assumption in Sections 5 to 8}. In 
applying the results of Sections 5 to 8 in Sections 9, 10, we 
need this assumption.

For the following theorem, we may even relax this condition even 
further:
\begin{thm}\label{thm:hcresint}
Suppose that $M_h$ is not projectively diffeomorphic to an 
open hemisphere. Suppose that $R_1$ and $R_2$ are two overlapping 
$n$-crescents that are hemispheres. Then $R_1 = R_2,$ and hence 
$\nu_{R_1} = \nu_{R_2}$ and $\alpha_{R_1} = \alpha_{R_2}$.
\end{thm}
\begin{pf}
We use Lemma \ref{lem:twocres} as in \cite{cdcr1}: By Theorem 
\ref{thm:ovlconv}, $\dev|R_1 \cup R_2$ is an imbedding onto the union 
of two $n$-hemispheres $\dev(R_1)$ and $\dev(R_2)$ in $\SI^n$. If $R_1$ 
is not equal to $R_2$, then $\dev(R_1)$ differs from $\dev(R_2)$, 
$\dev(R_1)$ and $\dev(R_2)$ meet each other in a convex $n$-bihedron, 
$\dev(R_1) \cup \dev(R_2)$ is homeomorphic to an $n$-ball, and the boundary 
$\delta (\dev(R_1) \cup \dev(R_2))$ is the union of two $(n-1)$-hemispheres 
meeting each other in a great $(n-2)$-sphere $\SI^{n-2}$.

Since $\alpha_{R_1}$ and $\alpha_{R_2}$ are disjoint from any of $R_1^o$ 
and $R_2^o$ respectively, the images of $\alpha_{R_1}$ and $\alpha_{R_2}$ 
do not intersect any of $\dev(R_1^o)$ and $\dev(R_2^o)$ respectively 
by Theorem \ref{thm:ovlconv}. Therefore, $\dev(\alpha_{R_1})$ and 
$\dev(\alpha_{R_2})$ are subsets of $\delta (\dev(R_1) \cup \dev(R_2))$. 
Since they are open $(n-1)$-hemispheres, the complement of 
$\dev(\alpha_{R_1}) \cup \dev(\alpha_{R_2})$ in 
$\delta (\dev(R_1) \cup \dev(R_2))$ equals $\SI^{n-1}$, and 
$\dev(\alpha_{R_1}) \cup \dev(\alpha_{R_2})$ is dense in 
$\delta(\dev(R_1) \cup \dev(R_2))$. Since $\dev| R_1 \cup R_2$ is an 
imbedding, it follows that $R_1 \cup R_2$ is an $n$-ball, and the closure 
of $\alpha_{R_1} \cup \alpha_{R_2}$ equals $\delta (R_1 \cup R_2)$.
Hence, $\delta (R_1 \cup R_2) \subset \hideal{M}$. By Lemma \ref{lem:include}
it follows that $M_h = R_1^o \cup R_2^o$, and $M_h$ is boundaryless. 
By the following lemma, this is a contradiction. Hence, $R_1 = R_2$.  
\end{pf}

\begin{lem}\label{lem:twocres}
Let $N$ be a closed real projective $n$-manifold. Suppose that 
$\dev: \che N_h \ra \SI^n$ is an imbedding onto the union of 
$n$-hemispheres $H_1$ and $H_2$ meeting each other in an $n$-bihedron 
or an $n$-hemisphere. Then $H_1 = H_2$, and $N_h$ is projectively 
diffeomorphic to an open $n$-hemisphere. 
\end{lem}
\begin{pf} 
Let $(\dev, h)$ denote the development pair of $N$, and $\Gamma$ 
the deck transformation group. As $\dev|N_h$ is a diffeomorphism 
onto $H_1^o \cup H_2^o$, a simply connected set, we have 
$N_h =\tilde N$.

Suppose that $H_1 \ne H_2$. Then $H_1\cup H_2$ is bounded by 
two $(n-1)$-hemispheres $D_1$ and $D_2$ meeting each other on 
a great sphere $\SI^{n-2}$, their common boundary.   
Since the interior angle of intersection of $D_1$ and $D_2$ is greater 
than $\pi$, $\delta H_i - D_i$ is an open hemisphere included in 
$\dev(\tilde N)$ for $i= 1, 2$. Defining $O_i = \delta H_i - D_i$ for 
$i=1, 2$, we see that $O_1 \cup O_2$ is $h(\Gamma)$-invariant since 
$\delta(H_1 \cup H_2)$ is $h(\Gamma)$-invariant. This means that the inverse 
image $\dev^{-1}(O_1 \cup O_2)$ is $\Gamma$-invariant. 

Let $O'_i = \dev^{-1}(O_i)$. Then elements of $\Gamma$ either act on each 
of $O'_1$ and $O'_2$ or interchange them. Thus, $\Gamma$ includes a 
subgroup $\Gamma'$ of index one or two acting on each of $O'_1$ and $O'_2$. 
Since $N_h$ is a simply connected open ball, and so is $O'_1$, it follows 
that the $n$-manifold $\tilde N /\Gamma'$ and an $(n-1)$-manifold 
$O'_1/\Gamma'$ are homotopy equivalent. Since $\tilde N/\Gamma'$ is a finite 
cover of a closed manifold $N$, $\tilde N/\Gamma'$ is a closed manifold. 
Since the dimensions of $\tilde N/\Gamma'$ and $O'_1/\Gamma$ are not the 
same, this is shown to be absurd by computing ${\bf Z}_2$-homologies.   
Hence we obtain that $H_1 = H_2$, and since $\dev(\tilde N)$ equals the 
interior of $H_1$, $\tilde N$ is diffeomorphic to an open $n$-hemisphere.
\end{pf}

\begin{figure}[h]
\centerline{\epsfxsize=3.75in \epsfbox{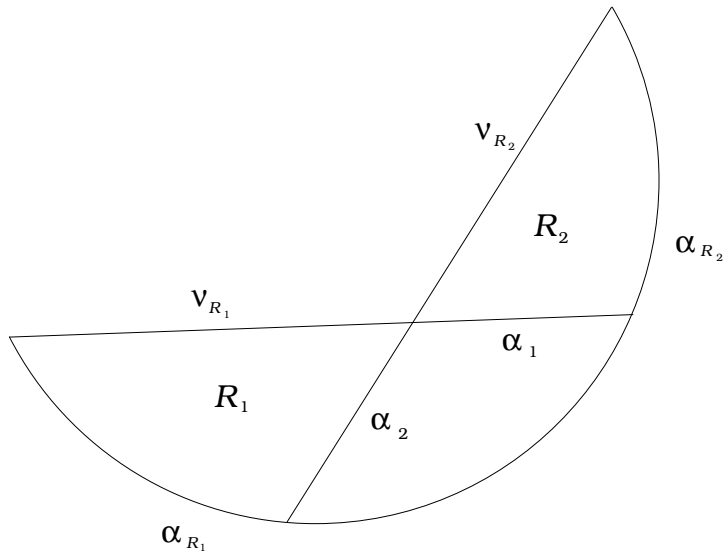}}
\caption{\label{fig:tr} 
Transversal intersections in dimension two.}
\end{figure}
\typeout{<<pftr.eps>>}

\begin{figure}[h]
\centerline{\epsfxsize=4.5in \epsfbox{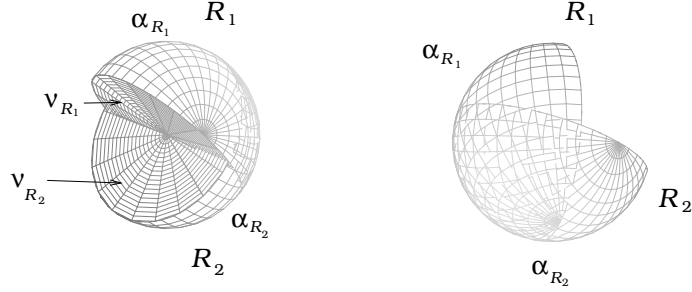}}
\caption{\label{fig:tr2} 
A three-dimensional transversal intersection seen 
in two view points}
\end{figure}

Suppose that $R_1$ is an $n$-crescent that is an $n$-bihedron. Let $R_2$ 
be another bihedral $n$-crescents with sets $\alpha_{R_2}$ and $\nu_{R_2}$.
We say that $R_1$ and $R_2$ intersect {\sl transversally}\/ 
if $R_1$ and $R_2$ overlap and the following conditions hold ($i=1, j= 2$; 
or $i= 2, j= 1$): 
\begin{enumerate}
\item $\nu_{R_1} \cap \nu_{R_2}$ is an $(n-2)$-dimensional hemisphere.
\item For the intersection $\nu_{R_1}\cap \nu_{R_2}$ denoted by $H$,
$H$ is an $(n-2)$-hemisphere,  $H^o$ is a subset of the interior
$\nu_{R_i}^o$, and $\dev(\nu_{R_i})$ and $\dev(\nu_{R_j})$ intersect 
transversally at points of $\dev(H)$. 
\item $\nu_{R_i} \cap R_j$ is a tame $(n-1)$-bihedron with boundary 
the union of $H$ and an $(n-2)$-hemisphere $H'$ in the closure of 
$\alpha_{R_j}$ with its interior $H^{\prime o}$ in $\alpha_{R_j}$.  
\item $\nu_{R_i}\cap R_j$ is the closure of a component of 
$\nu_{R_i} - H$ in $\che M_h$. 
\item $R_i \cap R_j$ is the closure of a component of $R_j - \nu_{R_i}$.
\item Both $\alpha_{R_i} \cap \alpha_{R_j}$ and 
$\alpha_{R_i}\cup \alpha_{R_j}$ are homeomorphic to open 
$(n-1)$-dimensional balls, which are locally totally geodesic 
under $\dev$.
\end{enumerate} 
Note that since $\alpha_{R_i}$ is tame, $\alpha_{R_i} \cap \alpha_{R_j}$ is
tame. (See Figures \ref{fig:tr} and \ref{fig:tr2}.)

By Corollary \ref{cor:transdev}, the above condition mirrors the property 
of intersection of $\dev(R_1)$ and $\dev(R_2)$ where $\dev(\alpha_{R_1})$ 
and $\dev(\alpha_{R_2})$ are included in a common great sphere $\SI^{n-1}$ 
of dimension $(n-1)$, $\dev(R_1)$ and $\dev(R_2)$ included in a common 
$n$-hemisphere bounded by $\SI^{n-1}$ and $\dev(\nu_{R_1})^o$ and 
$\dev(\nu_{R_2})^o$ meets transversally (see Theorem \ref{thm:ovlconv}).

\begin{exmp}
In the example \ref{exmp:dipexmp}, $R$ is an $n$-crescent with the 
closure of the plane $P$ given by the equation $x_1 + \cdots + x_n = 0$ 
equal to $\nu_R$. $\alpha_R$ equals the interior of the intersection of 
$R$ with $\delta H$. $\nu_S$ is the closure of the plane given by 
$x_1 = 0$ and $\alpha_S$ the interior of the intersection of $S$ 
with $\delta H$. Clearly, $R$ and $S$ intersect transversally.
\end{exmp}

\typeout{<<pftr3.eps>>}

Using the reasoning similar to Section 2.6. of \cite{cdcr1}, we obtain: 
\begin{thm}\label{thm:transversal} 
Suppose that $R_1$ and $R_2$ are overlapping. Then either $R_1$ and 
$R_2$ intersect transversally or $R_1 \subset R_2$ and $R_2 \subset R_1$.
\end{thm}

\begin{rem}
In case $R_1 \subset R_2$, we see easily that $\alpha_{R_1} = \alpha_{R_2}$ 
since the sides of $R_1$ in $\hideal{M}$ must be in that of $R_2$. 
Hence, we also see that $\nu^o_{R_1} \subset R_2^o$ as the topological
boundary of $R_1$ in $R_2$ must lie in $\nu_{R_1}$.
\end{rem}

The proof is entirely similar to that in \cite{cdcr1}. A good thing to
have in mind is the configurations of the images of two $n$-crescents
in $\SI^n$ meeting in many ways. We will see that only the configuration 
as indicated above will happen.

Assume that we have $i =1 \hbox{ and } j = 2$ or have $i =2 \hbox{ and } 
j =1$, and $R_1 \not\subset R_2$ and $R_2 \not\subset R_1$. Since 
$\clo(\alpha_{R_i}) \subset \hideal{M}$, Corollary \ref{cor:dipping} and 
Theorem \ref{thm:ovlconv} imply that $R_j$ dips into $(R_i, \nu_{R_i})$ and 
$\dev| R_i \cup R_j$ is an imbedding onto $\dev(R_1) \cup \dev(R_j)$.
Hence, the following statements hold: 
\begin{itemize} 
\item $\nu_{R_i} \cap R_j$ is a convex $(n-1)$-ball 
$\alpha_i$ such that 
\[\delta \alpha_i \subset \delta R_j, \alpha_i^o \subset R_j^o.\]
\item $R_i \cap R_j$ is the convex $n$-ball that is the closure of 
a component of $R_j - \alpha_i$.
\end{itemize}
Since $\alpha_i^o$ is disjoint from $\nu_{R_j}$, $\alpha_i^o$ is 
a subset of a component $C$ of $\nu_{R_i} - \nu_{R_j}$.

\begin{lem}\label{lem:tech}
If $\nu_{R_i}$ and $\nu_{R_j}$ meet\/{\rm ,} then they do so 
transversally\/{\rm ;} i.e, their images under $\dev$ meet transversally. 
If $\nu_{R_i}$ and $\alpha_{R_j}$ meet\/{\rm ,} then they do so
transversally. 
\end{lem}
\begin{pf} 
Suppose that $\nu_{R_i}$ and $\nu_{R_j}$ meet and they are tangential. 
Then $\dev(\nu_{R_i})$ and $\dev(\nu_{R_j})$ both lie on a common great 
$(n-1)$-sphere in $\SI^n$ by the geometry of $\SI^n$ and hence 
$\nu_{R_i} \cap \nu_{R_j} = \nu_{R_i} \cap R_j$ by Theorem \ref{thm:ovlconv}. 
Since $\nu_{R_i} \cap R_j$ includes an open $(n-1)$-ball $\alpha_i^o$, 
this contradicts $\alpha_i^o \subset R_j^o$. 

Suppose that $\nu_{R_i}$ and $\alpha_{R_j}$ meet and they are tangential. 
Then $\nu_{R_i} \cap \clo(\alpha_{R_j}) = \nu_{R_i} \cap R_j$ as before. 
This leads to contradiction similarly.
\end{pf}

We now determine a preliminary property of $\alpha_i$. Since 
$\alpha_i$ is a convex $(n-1)$-ball in $\nu_{R_i}$ with topological 
boundary in $\delta R_i \cup \delta R_j$, we obtain
\begin{eqnarray*}
\delta \alpha_i 
&\subset & \delta \nu_{R_i} \cup (\delta R_j \cap \nu_{R_i}^o) \\
&\subset & \delta \nu_{R_i} \cup (\nu_{R_j} \cap \nu_{R_i}^o)
\cup (\alpha_{R_j} \cap \nu_{R_i}^o).
\end{eqnarray*} 
Hence, we have
\begin{eqnarray*}
\delta \alpha_i = (\delta \alpha_i \cap \delta \nu_{R_i}) 
\cup (\delta \alpha_i \cap \nu_{R_j} \cap \nu_{R_i}^o)
\cup (\delta \alpha_i \cap \alpha_{R_j} \cap \nu_{R_i}^o).
\end{eqnarray*}
If $\alpha_{R_j}$ meets $\nu_{R_i}^o$, then since $\alpha_{R_j}$ is 
transversal to $\nu_{R_i}^o$ by Lemma \ref{lem:tech}, $\alpha_{R_j}$ must 
intersect $R_i^o$ by Theorem \ref{thm:ovlconv}. Since $\alpha_{R_j} \subset 
\hideal{M}$, this is  a contradiction. Thus, $\alpha_{R_j} \cap \nu_{R_i}^o 
= \emp$. We conclude 
\begin{equation}\label{eqn:boundai}
\delta \alpha_i = (\delta \alpha_i \cap \delta \nu_{R_i})
\cup (\delta \alpha_i \cap \nu_{R_j} \cap \nu_{R_i}^o).
\end{equation}

Let us denote by $H$ the set $\nu_{R_i} \cap \nu_{R_j}$. Consider 
for the moment the case where $\nu_{R_j} \cap \nu_{R_i}^o \ne \emp$. 
Then $\nu_{R_j} \cap \nu_{R_i}^o$ is a tame topological $(n-2)$-ball 
from the transversality in Lemma \ref{lem:tech} and Theorem 
\ref{thm:ovlconv}. If $H$ has boundary points, i.e, points of $\delta H$, 
in $\nu_{R_i}^o$, then $\nu_{R_i}^o - H$ would have only one component. 
Since the boundary of $\alpha_i$ in $\nu_{R_i}^o$ is 
included in $H$, the component must equal $\alpha_i^o$. Since $\alpha_i$ 
is a convex $(n-1)$-ball, $\alpha_i^o$ is convex; this is a contradiction. 
It follows that $H$ is an $(n-2)$-ball with boundary in $\delta \nu_{R_i}$ 
and the interior $H^o$ in $\nu_{R_i}^o$, i.e., $H$ separates $\nu_{R_i}$ into
two convex components, and the closures of each of them are $(n-1)$-bihedrons.
Since $\nu_{R_i}$ is an $(n-1)$-hemisphere, $H$ is an $(n-2)$-hemisphere. 
Moreover, $\alpha_i$ is the closure of a component of $\nu_{R_i} - H$. 
(Use Theorem \ref{thm:ovlconv} to see the possible configurations 
of the images of these sets in $\SI^n$.)

We need to consider only the following two cases by interchanging $i$ and 
$j$ if necessary:
\begin{itemize}
\item[(i)] $\nu_{R_j}^o \cap \nu_{R_i}^o \ne \emp$.
\item[(ii)] $\nu_{R_j} \cap \nu_{R_i}^o = \emp$
or $\nu_{R_i} \cap \nu_{R_j}^o = \emp$. 
\end{itemize}

(i) Since $\alpha_i$ is the closure of a component of $\nu_{R_i} - H$, 
$\alpha_i$ is an $(n-1)$-bihedron bounded by an $(n-2)$-hemisphere $H$ 
and another $(n-2)$-hemisphere $H'$ in $\delta \nu_{R_i}$. Since $H'$ 
is a subset of the closure of $\alpha_{R_i}$, it belongs to $\hideal{M}$ 
and hence disjoint from $R_j^o$ while $R_j^o \subset M_h$.

Since $H'$ is a subset of $R_j$, we have $H' \subset \delta R_j$.
Since $\nu_{R_i}$ is transversal to $\nu_{R_j}$, $H'$ is not a subset 
of $\nu_{R_j}$; thus, $H^{\prime o}$ is a subset of $\alpha_{R_j}$, 
and $H'$ that of $\clo(\alpha_{R_j})$ by Theorem \ref{thm:ovlconv}.
This completes the proof of the transversality properties (1)--(4) in 
case (i).

(5) By dipping intersection properties, $R_i \cap R_j$ is the closure of 
a component of $R_j -\alpha_i$ and hence that of $R_j - \nu_{R_i}$.

(6) Since $H^{\prime o}$ is a subset of $\alpha_{R_j}$, 
$\alpha_{R_j} - H'$ has two components $\beta_1$ and $\beta_2$, 
homeomorphic to open $(n-1)$-balls. By (5), we may assume without loss of 
generality that $\beta_1 \subset R_i$, and $\beta_2$ is disjoint from 
$R_i$. Since $\beta_1 \subset \hideal{M}$, we have $\beta_1 \subset 
\delta R_i$. As $\beta_1$ is a component of $\alpha_{R_j}$ removed 
with $H'$, we see that $\beta_1$ is an open $(n-1)$-bihedron 
bounded by $H'$ in $\delta \nu_{R_i}$ and an $(n-2)$-hemisphere $H''$ 
in $\delta \nu_{R_j}$.

Since the closure of $\beta_1$ belongs to $R_i$, we obtain 
that $H'' \subset R_i$ and $H''$ is a subset of $\alpha_j$, where 
$\alpha_j = \nu_{R_j} \cap R_i$. As $H''$ is a subset of 
$\delta \nu_{R_j}$, and $\alpha_j$ is the closure of a component 
of $\nu_{R_j}$ removed with $H$, we obtain $H'' \subset \delta \alpha_j$.

Since (i) holds for $i$ and $j$ exchanged, we obtain, by a paragraph above 
the condition (i), $H^o$ belongs to $\nu_{R_i}^o \cap \nu_{R_j}^o$. 
By (1)--(4) with values of $i$ and $j$ exchanged, $\alpha_j$ is an 
$(n-1)$-bihedron bounded by $H$ and an $(n-2)$-hemisphere $H'''$ with 
interior in $\alpha_{R_i}$ and is the closure of a component of 
$\nu_{R_j} - H$. Since $H''$ is an $(n-1)$-hemisphere in 
$\delta \alpha_j$, and so is $H'''$, it follows that $H'' = H'''$.

Since $\beta_1$ has the boundary the union of $H'$ in $\delta \nu_{R_i}$ 
and $H''$, $H'' = H'''$, with interior in $\alpha_{R_i}$, and $\beta_1$ 
is a convex subset of $R_i$, looking at the bihedron $\dev(R_i)$ and 
the geometry of $\SI^n$ show that $\beta_1 \subset \alpha_{R_i}$. Thus, 
we obtain $\beta_1 \subset \alpha_{R_i} \cap \alpha_{R_j}$. We see that 
$\dev(\alpha_{R_i})$ and $\dev(\alpha_{R_j})$ are subsets of a common 
great $(n-1)$-sphere; it follows easily by Theorem \ref{thm:ovlconv}
that $\beta_1 = \alpha_{R_i} \cap \alpha_{R_j}$. Hence, 
$\alpha_{R_i} \cap \alpha_{R_j}$ and $\alpha_{R_i} \cup \alpha_{R_j}$ 
are homeomorphic to open $(n-1)$-balls, and under $\dev$ they map
to totally geodesic $(n-1)$-balls in $\SI^n$. 

(ii) Assume $\nu_{R_j}^o \cap \nu_{R_i} = \emp$ without loss of 
generality. Then $H$ is a subset of $\delta \nu_{R_j}$. Since 
$R_i$ dips into $(R_j, \nu_{R_j})$, we have that $\alpha_j \ne \emp$. 
Since the boundary of $\alpha_j$ in $\nu_{R_j}^o$ is included in $H$ 
(see equation \ref{eqn:boundai}), we have $\alpha_j = \nu_{R_j}$. 
It follows that $\nu_{R_j}^o$ is a subset of $R_i^o$, and 
$R_i \cap R_j$ is the closure of a component of $R_i -\nu_{R_j}$. 
Since $\nu_{R_j}$ is an $(n-1)$-hemisphere, and $R_i$ is an 
$n$-bihedron, the uniqueness of $(n-1)$-spheres in an $n$-bihedron 
(Theorem \ref{thm:class2}) shows that $\delta \nu_{R_i} = \delta \nu_{R_j}$.
Thus, the closures of components of $R_i - \nu_{R_j}$ are $n$-bihedrons 
with respective boundaries $\alpha_{R_i} \cup \nu_{R_j}$ and 
$\nu_{R_i} \cup \nu_{R_j}$. By Corollary \ref{cor:dipping}, 
$R_i \cap R_j$ is the closure of either the first component or the 
second one. 

In the first case, $R_i^o \cap R_j^o$ is an open subset of $R_j^o$ 
since $R_i^o$ is open in $M_h^o$. The closure of $R_i^o$ in $M_h$ 
equals $R_i^o \cup (\nu_{R_i} \cap M_h) = R_i \cap M_h$. Since 
$\nu_{R_i}^o$, which contains $\nu_{R_i} \cap M_h$, do not meet 
$R_j$ in the first case being in the other component of 
$R_i - \nu_{R_j}$, we see that the intersection of the closure of 
$R_i^o$ in $M_h$ with $R_j^o$ is same as $R_i^o \cap R_j^o$. Thus, 
$R_i^o \cap R_j^o$ is open and closed subset of $R_j^o$. Hence 
$R_j^o \subset R_i^o$ and $R_i \subset R_j$. This contradicts our 
hypothesis.

In the second case, $\dev| R_i \cup R_j$ is a homeomorphism to 
$\dev(R_i) \cup \dev(R_j)$. As $\alpha_{R_i}$ and $\alpha_{R_j}$ are 
subsets of $\hideal{M}$, their images under $\dev$ may not meet that of 
$R_i^o \cup R_j^o$. Hence, $\dev(R_i) \cup \dev(R_j)$ is an $n$-ball 
bounded by $\dev(\clo(\alpha_{R_i}))$ and $\dev(\clo(\alpha_{R_j}))$.
We obtain that $R_i \cup R_j$ is the $n$-ball bounded by 
$(n-1)$-dimensional hemispheres $\clo(\alpha_{R_i})$ and 
$\clo(\alpha_{R_j})$. 

Since $\clo(\alpha_{R_i})$ and $\clo(\alpha_{R_j})$ are subsets of 
$\hideal{M}$, Lemma \ref{lem:include} shows that $\che M_h=R_i\cup R_j$ 
and $M_h = R_i^o \cup R_j^o$; thus, $M_h = \tilde M$ and $M$ is 
a closed manifold. The image $\dev(R_1)\cup \dev(R_2)$ is bounded by two 
$(n-1)$-hemispheres meeting each other on a great sphere $\SI^{n-2}$, 
their common boundary. Since $M_h$ is not projectively diffeomorphic 
to an open $n$-hemisphere or an open $n$-bihedron, the interior angle of 
intersection of the two boundary $(n-1)$-hemisphere should be greater 
than $\pi$. However, Lemma \ref{lem:twocres} contradicts this. \qed

\begin{rem} 
Using the same proof as above, we may drop the condition on the Euler 
characteristic from Theorem 2.6 of \cite{cdcr1} if we assume that 
$\tilde M$ is not diffeomorphic to an open $2$-hemisphere or an open lune. 
This is weaker than requiring that the Euler characteristic of $M$ is 
less than zero. So, our theorem is an improved version of Theorem 2.6 of 
\cite{cdcr1}.
\end{rem}

\begin{cor}\label{cor:transdev}
Let $R_1$ and $R_2$ be bihedral $n$-crescents and they overlap. 
Then $\dev(\alpha_{R_1})$ and $\dev(\alpha_{R_2})$ are included 
in a common great $(n-1)$-sphere $\SI^{n-1},$ and $\dev(R_1)$ 
and $\dev(R_2)$ are subsets of a common great $n$-hemisphere 
bounded by ${\bf S}^{n-1}$. Moreover{\rm ,} 
$\dev(R_i - \clo(\alpha_{R_i}))$ 
is a subset of the interior of this $n$-hemisphere for $i= 1, 2$.
\qed
\end{cor}
%
%

\vfill
\break

\section{$n$-crescents that are $n$-hemispheres and the two-faced
submanifolds}
\label{sec:hemispheres}
In this section, we introduce the two-faced submanifolds arising 
from hemispheric $n$-crescents. We showed above that if two 
hemispheric $n$-crescents overlap, then they are equal. We show that 
if two hemispheric $n$-crescents meet but do not overlap, then 
they meet at the union of common components of their $\nu$-boundaries, 
which we call copied components. The union of all copied components 
becomes a properly imbedded submanifold in $M_h$ and covers a properly 
imbedded submanifold in $M$. This is the two-faced submanifold.

\begin{lem}\label{lem:bdandcres}
Let $R$ be an $n$-crescent. A component of $\delta M_h$ is either 
disjoint from $R$ or is a component of $\nu_R \cap M_h$. Moreover, 
a tiny ball $B(x)$ of a point $x$ of $\delta M_h$ is a subset of $R$ 
if $x$ belongs to $\nu_R \cap M_h${\rm ,} and{\rm ,} consequently{\rm ,} 
$x$ belongs to the topological interior $\inte R$.
\end{lem}
\begin{pf} 
If $x \in \delta M_h$, then a component $F$ of the open $(n-1)$-manifold 
$\nu_R \cap M_h$ intersects $\delta M_h$ tangentially, and by Lemma 
\ref{lem:bdtangent}, it follows that $F$ is a subset of $\delta M_h$.  
Since $F$ is a closed subset of $\nu_R \cap M_h$, $F$ is a closed subset 
of $\delta M_h$. Since $F$ is an open manifold, $F$ is open in $\delta M_h$. 
Thus, $F$ is a component of $\delta M_h$.
 
Since $x \in \inte B(x)$, $B(x)$ and $R$ overlap. As $\clo(\alpha_R)$ 
is a subset of $\hideal{M}$, we have $\Bd R \cap B(x) \subset \nu_R$ and 
$\nu_R \cap B(x) = F \cap B(x)$ for a component $F$ of $\nu_R \cap M_h$ 
containing $x$. Since $F$ is a component of $\delta M_h$, we obtain 
$F \cap B(x) \subset \delta B(x)$; since we have
$\Bd R \cap B(x) \subset \delta B(x)$, it follows that $B(x)$ is 
a subset of $R$. 
\end{pf}

Suppose that $\che M_h$ includes an $n$-crescent $R$ that is 
an $n$-hemisphere. Then $M_h \cap R$ is a submanifold of $M_h$ 
with boundary $\delta R \cap M_h$. Since $R$ is an $n$-crescent, 
$\delta R \cap M_h$ equals $\nu_R \cap M_h$. Let $B_R$ denote 
$\nu_R \cap M_h$. 

Let $S$ be another hemispheric $n$-crescent, and $B_S$ the set 
$\nu_S \cap M_h$. By Theorem \ref{thm:hcresint}, we see that either 
$S \cap R^o =\emp$ or $S = R$. Suppose that $S \cap R \ne \emp$ 
and $S$ does not equal $R$. Then $B_S \cap B_R \ne \emp$. Let $x$ be 
a point of $B_S \cap B_R$ and $B(x)$ the tiny ball of $x$. 
Since $\inte B(x) \cap R \ne \emp$, it follows that $B(x)$ dips into 
$(R, \nu_R)$ or $B(x)$ is a subset of $R$ by Lemma \ref{lem:tinydip}.   
Similarly, $B(x)$ dips into $(S, \nu_{S})$ or $B(x)$ is a subset of $S$.
If $B(x)$ is a subset of $R$, then $S$ intersects the interior of $R$. 
Theorem \ref{thm:hcresint} shows $S = R$, a contradiction. Therefore, 
$B(x)$ dips into $(R, \nu_R)$ and similarly into $(S, \nu_S)$. 
If $\nu_R$ and $\nu_S$ intersect transversally, then $R$ and $S$ 
overlap. This means a contradiction $S = R$. Therefore, $B_S$ and 
$B_R$ intersect tangentially at $x$.

If $x \in \delta M_h$, Lemma \ref{lem:bdandcres} shows that $B(x)$ is 
a subset of $R$. This contradicts a result of the above paragraph. 
Thus, $x \in M_h^o$. Hence, we conclude that $B_R \cap B_S \cap M_h 
\subset M_h^o$.    

Since $B_R$ and $B_S$ are closed subsets of $M_h$, and $B_R$ and $B_S$ are 
totally geodesic and intersect tangentially at $x$, it follows that 
$B_R \cap B_S$ is an open and closed subset of $B_R$ and $B_S$ respectively.  
Thus, for components $A$ of $B_R$ and $B$ of $B_S$, either we have 
$A = B$ or $A$ and $B$ are disjoint. Therefore, we have proved:

\begin{thm}\label{thm:inthcrs}
Given two hemispheric $n$-crescents $R$ and $S$, we have either 
$R$ and $S$ disjoint, or $R$ equals $S$, or $R \cap S$ equals 
the union of common components of $\nu_R \cap M_h$ and 
$\nu_S \cap M_h$.
\end{thm}
Readers may easily find examples where $\nu_R \cap M_h$ 
and $\nu_S \cap M_h$ are not equal in the above situations
especially if $M$ is an open manifold. 

\begin{defn}\label{defn:copied1}
Given a hemispheric $n$-crescent $T$, we say that a component of 
$\nu_T \cap M_h$ is {\em copied}\/ if it equals a component of 
$\nu_U \cap M_h$ for some hemispheric $n$-crescent $U$ not equal to
to $T$. 
\end{defn}

Let $c_R$ be the union of all copied components of $\nu_R \cap M_h$
for a hemispheric $n$-crescent $R$. Let $A$ denote 
$\bigcup_{R \in \cal H} c_R$ where $\cal H$ is the set of all 
hemispheric $n$-crescents in $M_h$. $A$ is said to be 
the {\em pre-two-faced submanifold arising from hemispheric 
$n$-crescents}.

\begin{prop}\label{prop:spmfld}
Suppose that $A$ is not empty. Then $A$ is a properly imbedded 
totally geodesic $(n-1)$-submanifold of $M_h^o$ and $p| A$ is 
a covering map onto a closed totally geodesic imbedded 
$(n-1)$-manifold in $M^o$.   
\end{prop}

First, given two $n$-crescents $R$ and $S$, $c_R$ and $c_S$ meet either in 
the union of common components or in an empty set: Let $a$ and $b$ be 
respective components of $c_R$ and $c_S$ meeting each other. Then $a$ 
is a component of $\nu_R \cap M_h$ and $b$ that of $\nu_S \cap M_h$. 
Since $R \cap S \ne \emp$, either $R$ and $S$ overlap or $a=b$ by the 
above argument. If $R$ and $S$ overlap, $R = S$ and hence $a$ and $b$ 
must be the identical component of $\nu_R \cap M_h$ and hence $a = b$. 
Hence, $A$ is a union of mutually disjoint closed path-components that 
are components of $c_R$ for some $n$-crescent $R$. In other words, $A$ is 
a union of path-components which are components of $c_R$ for some $R$.

Second, given a tiny ball $B(x)$ of a point $x$ of $M_h$, we claim that 
no more than one path-component of $A$ may intersect $\inte B(x)$: Let 
$a$ be a component of $c_R$ intersecting $\inte B(x)$. Since copied 
components are subsets of $M_h^o$, $a$ intersects $B(x)^o$ and hence 
$B(x)$ is not a subset of $R$. By Lemma \ref{lem:tinydip}, 
$\nu_R \cap B(x)$ is a compact convex $(n-1)$-ball with boundary in 
$\Bd B(x)$. Since it is connected, $a \cap B(x) = \nu_R \cap B(x)$, and 
$B(x) \cap R$ is the closure of a component $C_1$ of $B(x) - (a \cap B(x))$ 
by Corollary \ref{cor:dipping}. Since $a$ is copied, $a$ is a component 
of $\nu_S \cap M_h$ for an $n$-crescent $S$ not equal to $R$, 
and $B(x) \cap S$ is the closure of a component $C_2$ of 
$B(x)-(a \cap B(x))$. Since $R$ and $S$ do not overlap, it follows that 
$C_1$ and $C_2$ are the two disjoint components of $B(x)-(a \cap B(x))$.

Suppose that $b$ is a component of $c_T$ for an $n$-crescent $T$ 
and $b$ intersects $\inte B(x)$ also. If the $(n-1)$-ball $b \cap B(x)$  
intersects $C_1$ or $C_2$, then $T$ overlaps $R$ or $S$ respectively 
and hence $T = R$ or $T= S$ respectively by Theorem \ref{thm:hcresint}; 
therefore, we have $a = b$. This is absurd. Hence $b \cap B(x) \subset 
a\cap B(x)$ and $T$ overlaps with either $R$ or $S$. Since these are 
hemispheric crescents, we have either $T = R$ or $T = S$ respectively; 
therefore, $a = b$. 

Since given a tiny ball $B(x)$ for a point $x$ in $M_h$, no more than 
one distinct path-component of $A$ may intersect $\inte B(x)$, each 
path-component of $A$ is an open subset of $A$. This shows that 
$A$ is a totally geodesic $(n-1)$-submanifold of $M_h^o$, closed
and properly imbedded in $M_h^o$.
 
Let $p: M_h \ra M$ be the covering map. Since $A$ is the deck 
transformation group invariant, we have $A = p^{-1}(p(A))$ and $p|A$ 
covers $p(A)$. The above results show that $p(A)$ is a closed totally 
geodesic manifold in $M^o$.  

\begin{defn}\label{defn:twoface1}
$p(A)$ for the union $A$ of all copied components of hemispheric 
$n$-crescents in $\che{M_h}$ is said to be the {\em two-faced 
$(n-1)$-manifold}\/ of $M$ {\em arising from hemispheric $n$-crescents}\/
(or {\em type I}\/).
\end{defn} 

Each component of $p(A)$ is covered by a component of $A$, i.e.,
a copied component of $\nu_R \cap M_h$ for some crescent $R$. Since 
$\alpha_R$ is the union of the open $(n-1)$-hemispheres in $\delta R$, 
$\nu_R \cap M_h$ lies in an open $(n-1)$-hemisphere, i.e., 
an affine patch in the great $(n-1)$-sphere $\delta R$. Hence, each 
component of $p(A)$ is covered by an  open domain in $\bR^n$.  

We end with the following observation: 
\begin{prop}\label{prop:hnonintf}
Suppose that $A = \bigcup_{R \in \cal H} c_R$. Then $A$ is disjoint from 
$S^o$ for each hemispheric $n$-crescent $S$ in $\che M_h$. 
\end{prop}
\begin{pf}
If not, then a point $x$ of $c_R$ meets $S^o$ for some hemispheric 
$n$-crescent $S$. But if so, then $R$ and $S$ overlap, and $R = S$.
This is a contradiction.
\end{pf}

\vfill
\break

\section{$n$-crescents that are $n$-bihedrons and the two-faced 
submanifolds}\label{sec:bitwof}

In this section, we will define an equivariant set $\Lambda(R)$ 
for a bihedral $n$-crescent $R$. We discuss its properties which 
are exactly same as those of its two-dimensional version in 
\cite{cdcr1}. Then we discuss the two-faced submanifolds that 
arises from $\Lambda(R)$'s: We show that $\Lambda(R)$ and 
$\Lambda(S)$ for two $n$-crescents are either equal or disjoint 
or meet at their common boundary components in $M_h$. The union 
of all such boundary components for $\Lambda(R)$ for every bihedral 
$n$-crescents $R$ is shown to be a properly imbedded submanifold 
in $M_h$ and cover a compact submanifold in $M$. The proof of this 
fact is similar to those in the previous section.

We will suppose in this section that $\che M_h$ includes no hemispheric 
crescent; we assume that all $n$-crescents in $\che M_h$ are bihedrons. 
Two bihedral $n$-crescents in $\che M_h$ are equivalent if they overlap. 
This generates an equivalence relation on the collection of all bihedral 
$n$-crescents in $\che M_h$; that is, $R \sim S$ if and only if there 
exists a sequence of bihedral $n$-crescents $R_i$, $i=1, \dots, n$, such 
that $R_1 = R, R_n = S$ and $R_{i-1} \cap R_i^o \ne \emp$ for $i= 2, 
\dots, n$.

We define 
\[ \Lambda(R) := \bigcup_{S \sim R} S, \qquad
\delta_\infty \Lambda(R) := \bigcup_{S \sim R} \alpha_S, \qquad 
\Lambda_1(R) := \bigcup_{S \sim R} (S - \nu_R). \]

\begin{figure}[h]
\centerline{\epsfxsize=3.5in \epsfbox{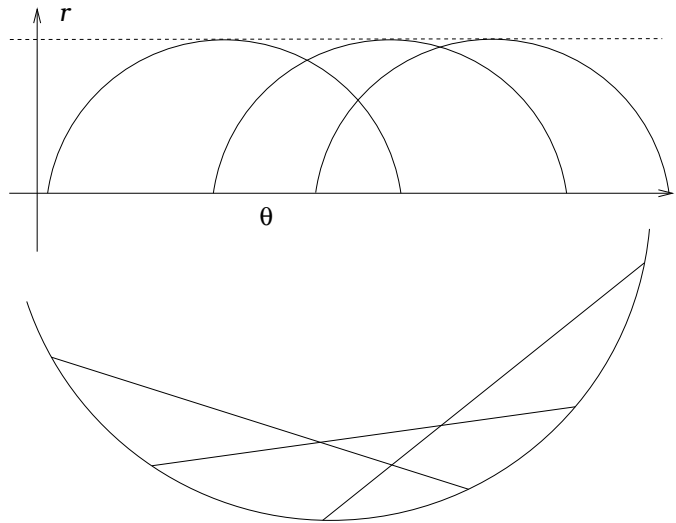}}
\caption{\label{fig:lambdaR} Figures of $\Lambda(R)$.}
\end{figure}

\begin{exmp}\label{exmp:lambda}
Consider the universal cover $L$ of $H^o - \{O\}$ where $H$ is 
a $2$-hemisphere in $\SI^2$. Then it has an induced real projective
structure with developing map equal to the covering map $c$.
There is a nice parameterization $(r, \theta)$ of $L$ where $r$ denotes
the $\bdd$-distance of $c(x)$ from $O$ and $\theta(x)$ the oriented 
total angle from the lift of the positive $x$-axis for $x \in L$, 
i.e., one obtained by integrating the $1$-form $d\theta$. Here 
$r$ is from $[0, \pi/2]$ and $\theta$ in $(-\infty, \infty)$.
$L$ is hence a holonomy cover of itself as it is simply connected.
$\che L$ may be identified with the universal cover of $H - \{O\}$
with a point $O'$ added to make it a complete space where 
$O'$ maps to $O$ under the extended developing map $c$.
(We use the universal covering space since the holonomy 
cover gives us trivial examples. Of course, the holonomy cover 
of a universal cover is itself.)

A crescent in $\che L$ is the closure of a lift of an affine half space 
in $H^o - \{O\}$. (Recall that $H^o$ has an affine structure.) 
A special type of a crescent is the closure of the set given by 
$\theta_0 \leq \theta \leq \theta_0 +\pi$. Given a crescent $R$ 
in $\che L$, we see that $\Lambda(R)$ equals $\che L$.

We may also define another real projective manifold $N$ by an equation 
$f(\theta) < r < \pi/2$ for a function $f$ with values in $(0, \pi/2)$. 
Then $\che N$ equals the closure of $N$ in $\che L$. 
Given a crescent $R$ in $\che N$, we see that $\Lambda(R)$ may not
equal to $\che N$ in especially in case $f$ is not a convex function 
(as seen in polar coordinates). See figure 
\ref{fig:lambdaR}.

For a higher dimensional example, let $H$ be a $3$-hemisphere in 
$\SI^3$, and $l$ a segment of $\bdd$-length $\pi$ passing through 
the origin. Let $L$ be the universal cover of $H^o -l$. Then $L$
becomes a real projective manifold with developing map the covering 
map $c:L \ra H^o - l$. The holonomy cover of $L$ is $L$ itself. The 
completion $\che L$ of $L$ equals the completion of the universal cover 
of $H - l$ with $l$ attached to make it a complete space. A $3$-crescent 
is the closure of a lift of an open half space in $H - l$. Given 
a $3$-crescent $R$, $\Lambda(R)$ equals $\che L$.

We introduce coordinates on $H^o$ so that $l$ is now the $z$-axis.
Note that $L$ is parameterized by $(r, \theta, \phi)$ where $r(x)$ 
equals the $\bdd$-distance from $O$ to $c(x)$ and $\phi$ the angle 
that $\ovl{Oc(x)}$ makes with the ray in $L$ from the origin in a 
given direction, and $\theta(x)$ the angle from the lift of the 
half-$xz$-plane given by $x > 0$. We may also define other real 
projective manifolds by equation $f(\theta, \phi) < r < \pi/2$ for 
$f: \bR \times (0, \pi) \ra (0, \pi/2)$. The readers may work out how the 
completions might look and what $\Lambda(R)$ may look when $R$ is a 
$3$-crescent. We remark that for certain $f$ which converges to $\pi/2$ 
as $\phi \ra 0$ or $\pi$, we may have no $3$-crescents in the completion
of the real projective manifold given by $f$. 

Even higher-dimensional examples are given in a similar spirit by 
removing sets from such covers. We will see that what we gave 
are really typical examples of $\Lambda(R)$.
\end{exmp}

\typeout{<<pflam.eps>>}

Let us state the properties that hold for these sets: The proofs are 
exactly as in \cite{cdcr1}. 
\begin{eqnarray*}
\inte \Lambda (R) \cap M_h &=& \inte (\Lambda(R)\cap M_h)\\
\Bd \Lambda (R) \cap M_h &=& \Bd (\Lambda(R)\cap M_h) \cap M_h
\end{eqnarray*}
(see Lemma 6.4 \cite{cdcr1}). 
For $\vth$ a deck transformation, from definitions we easily obtain
\begin{eqnarray}\label{eqn:lambda2}
\vth(\Lambda(R)) &=& \Lambda(\vth(R)), \nonumber \\ 
\vth(\delta_\infty \Lambda(R)) &=& \delta_\infty \Lambda(\vth(R)) \nonumber \\
\vth(\Lambda_1(R)) &=& \Lambda_1(\vth(R)) \\ 
\inte \vth(\Lambda (R)) \cap M_h &=& \vth(\inte \Lambda(R)) \cap M_h
= \vth(\inte \Lambda(R) \cap M_h) \nonumber \\ 
\Bd \vth(\Lambda(R)) \cap M_h &=& \vth(\Bd \Lambda(R)) \cap M_h = 
\vth(\Bd \Lambda(R) \cap M_h). \nonumber
\end{eqnarray}

The sets $\Lambda(R)$ and $\Lambda_1(R)$ are path connected. 
$\delta_\infty \Lambda(R)$ is an open $(n-1)$-manifold. 
Since Theorem \ref{thm:transversal} shows that for two overlapping 
$n$-crescents $R_1$ and $R_2$, $\alpha_{R_1}$ and $\alpha_{R_2}$ extend 
each other into a larger $(n-1)$-manifold, there exists a unique great 
sphere $\SI^{n-1}$ including $\dev(\delta_\infty \Lambda(R))$ and 
by Corollary \ref{cor:transdev}, a unique component $A_R$ of 
$\SI^n - \SI^{n-1}$ such that $\dev(\Lambda(R)) \subset \clo(A_R)$ 
and $\dev(\Lambda(R) - \clo(\delta_\infty \Lambda(R))) \subset A_R$. 
For a deck transformation $\vth$ acting on $\Lambda(R)$, $A_R$ is 
$h(\vth)$-invariant. $\Lambda_1(R)$ admits a real projective 
structure as a manifold with totally geodesic boundary 
$\delta_\infty \Lambda(R)$. 

\begin{prop}\label{prop:lamcl}
$\Lambda(R) \cap M_h$ is a closed subset of $M_h$. 
\end{prop}
\begin{pf}
This follows as in Section 6.1 of \cite{cdcr1}. We show that each 
point of $\Bd \Lambda(R) \cap M_h$ belongs to $\Lambda(R)$ by using 
a sequence of points converging to it and a sequence of $n$-crescents 
containing it using sequences of crescents (see Lemma \ref{lem:seqcres}). 
\end{pf}

\begin{lem}\label{lem:concavbd}
$\Bd \Lambda(R) \cap M_h$ is a properly imbedded topological submanifold
of $M_h^o$, and $\Lambda(R) \cap M_h$ is a topological submanifold of 
$M_h$ with concave boundary $\Bd \Lambda(R) \cap M_h$.
\end{lem} 
\begin{pf} 
Let $p$ be a point of $\Bd \Lambda(R) \cap M_h$. Since $\Lambda(R)$ is 
closed, $p$ is a point of a crescent $R'$ equivalent to $R$.

Let $B(p)$ be an open tiny ball of $p$. Since by Lemma \ref{lem:bdandcres}, 
$\Bd \Lambda(R) \cap M_h$ is a subset of $M_h^o$, $B(p)^o$ is an open 
neighborhood of $p$. Since $B(p)^o \cap \Lambda(R)$ is a closed subset 
of $B(p)^o$, $O = B(p)^o - \Lambda(R)$ is an open subset. 

We claim that $O$ is a convex subset of $B(p)^o$. Let $x, y \in O$. 
Then let $s$ be the segment in $B(p)$ of $\bdd$-length $\leq \pi$ 
connecting $x$ and $y$. If $s^o \cap \Lambda(R) \ne \emp$, then a point 
$z$ of $s^o$ belongs to an $n$-crescent $S$, $S \sim R$. If $z$ belongs 
to $S^o$, since $s$ must leave $S$, $s$ meets $\nu_S$ and is transversal 
to $\nu_S$ at the intersection point. Since a maximal line in the bihedron 
$S$ transversal to $\nu_S$ have an endpoint in $\alpha_S$, at least one 
endpoint of $s$ belongs to $S^o$, which is a contradiction. If $z$ belongs 
to $\nu_S$ and $s$ is transversal to $\nu_S$ at $z$, the same argument 
gives us a contradiction. If $z$ belongs to $\nu_S$ and $s$ is tangential 
to $\nu_S$ at $z$, then $s$ is included in the component of $\nu_S \cap M_h$ 
containing $z$ since $s \subset M_h$ is connected. Since $x$ and $y$ 
belong to $O$, this is a contradiction. Hence $s \subset O$, and $O$ 
is convex. 

Since $O$ is convex and open, $\Bd O$ in $M_h$ is homeomorphic to 
an $(n-1)$-sphere by Lemma \ref{thm:topconv}. The boundary 
$\Bd_{B(p)^o} O$ of $O$ relative to $B(p)^o$ equals $\Bd O \cap B(p)^o$. 
Hence, $\Bd_{B(p)^o} O$ is an imbedded open $(n-1)$-submanifold 
of $B(p)^o$. Since $\Bd \Lambda(R) \cap B(p)^o = \Bd_{B(p)^o} O$, 
$\Bd \Lambda(R) \cap M_h$ is an imbedded $(n-1)$-submanifold.
\end{pf}

\begin{figure}[h]
\centerline{\epsfxsize=5in \epsfbox{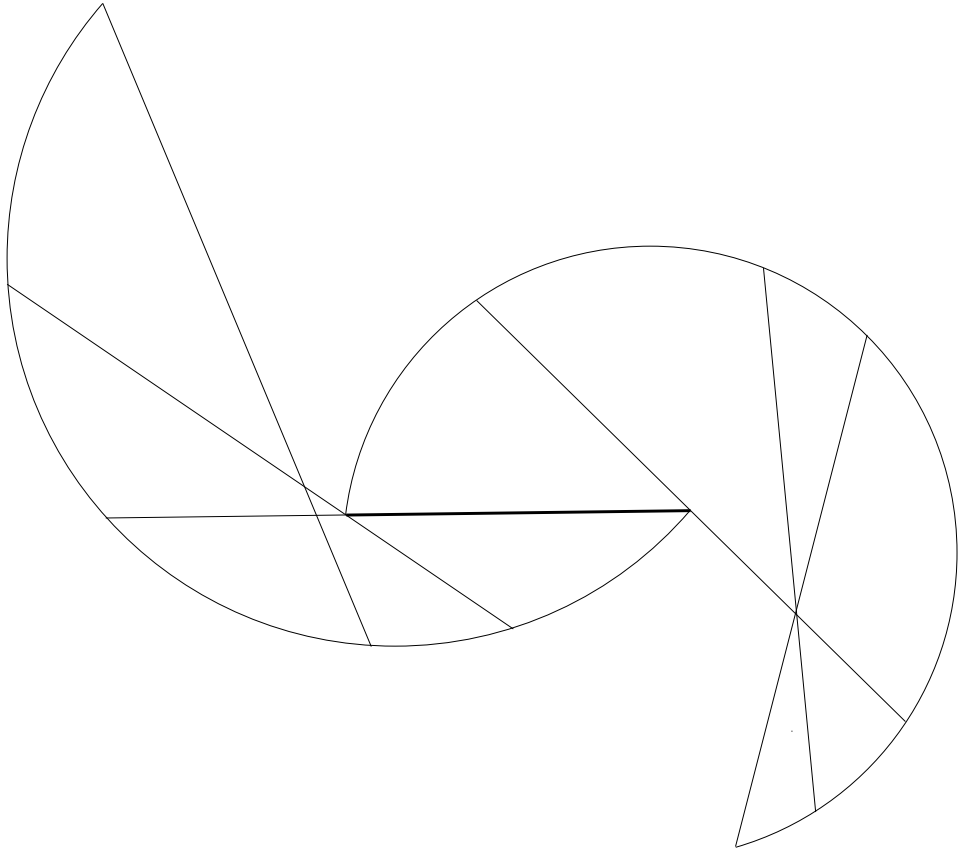}}
\caption{\label{fig:twof} 
A pre-two-faced submanifold.}
\end{figure}
\typeout{<<pftwof.eps>>}

Using the same argument as in Section 6.2 of \cite{cdcr1} 
(see Lemma 6.4 of \cite{cdcr1}), we obtain the following lemma: 
\begin{lem}\label{lem:prestab}
If $\inte\Lambda(R) \cap M_h \cap \Lambda(S) \ne \emp$ for an 
$n$-crescent $S${\rm ,} then $\Lambda(R) = \Lambda(S)$. Moreover{\rm ,} 
if for a crescent $S${\rm ,} $\Lambda(R) \cap M_h$ and $\Lambda(S)\cap M_h$ 
meet and they are distinct, then $\Lambda(R) \cap \Lambda(S) \cap M_h$ 
is a subset of $\Bd \Lambda(R) \cap M_h$ and $\Bd \Lambda(S) \cap M_h$. 
\end{lem}
%
%

Assume now that $\Lambda(R)$ and $\Lambda(S)$ are distinct and meet 
each other. Consequently $R$ and $S$ are not equivalent. Let $x$ be 
a common point of $\Bd \Lambda(R)$ and $\Bd \Lambda(S)$, and $B(x)$ 
a tiny-ball neighborhood of $x$. By Lemma \ref{lem:concavbd}, $x \in M_h^o$ 
and so $x \in B(x)^o$. Let $T$ be a crescent equivalent to $R$ containing 
$x$, and $T'$ that equivalent to $S$ containing $x$. Then $T \cap B(x)$ is 
the closure of a component $A$ of $B(x) - P$ for a totally geodesic 
$(n-1)$-ball $P$ in $B(p)$ with boundary in $\Bd B(x)$ by Lemma 
\ref{lem:tinydip}. Moreover, $\nu_T \cap B(x) = P$ and $T^o \cap B(x) = A$,
and $A$ is a subset of $\inte \Lambda(R)$. Let $B$ denote $B(x)$ removed 
with $A$ and $P$. Similarly, $T' \cap B(x)$ is the closure of a component 
$A'$ of $B(x) -P'$ for a totally geodesic $(n-1)$-ball $P'$, 
$P' = \nu_{T'} \cap B(x)$ with boundary in $\Bd B(x)$, and $A'$ is a subset 
of $T^{\prime o}$ in $\inte \Lambda(S)$. Since $T^o \subset \inte \Lambda(R)$, 
$T' \subset \Lambda(S)$, $T' \cap B(x)$ and $T^o \cap B(x)$ are disjoint,
and $P, P' \ni x$, it follows that $P = P'$ and $B = A'$; that is, 
$P$ and $P'$ are tangential. (We have that $P = P'= \nu_S \cap B(x)= 
\nu_{S'} \cap B(x). \quad )$

Since $B$ is a subset of $\inte \Lambda(S)$, $B$ contains no point 
of $\Lambda(R)$ by Lemma \ref{lem:prestab}, and similarly $A$ contains no 
point of $\Lambda(S)$. Thus, $\Lambda(R) \cap B(x)$ is a subset of the 
closure of $A$, and $\Lambda(S) \cap B(x)$ is that of $B$.  Since 
$A \subset \inte \Lambda(R)$ and $B \subset \inte \Lambda(S)$, it follows 
that 
\begin{eqnarray}\label{eqn:tinyint}
A = \inte \Lambda(R) \cap B(x), & & B = \inte \Lambda(S) \cap B(x), 
\nonumber \\
P \cup A = \Lambda(R) \cap B(x), & &  P \cup B = \Lambda(S) \cap B(x), 
\\
P = \Bd \Lambda(R) \cap B(x) &=& \Bd \Lambda(S) \cap B(x).  
\nonumber
\end{eqnarray} 
Hence, we have $P = \Bd \Lambda(R) \cap \Bd \Lambda(S) \cap B(x)$ 
and $P$ is a totally geodesic $(n-1)$-ball with boundary in $\Bd B(x)$
and our point $x$ belongs to $P^o$, to begin with. Since this holds for 
an arbitrary choice of a common point $x$ of $\Bd \Lambda(R)$ and 
$\Bd \Lambda(S)$, a tiny ball $B(x)$ of $x$, it follows that 
$\Bd \Lambda(R) \cap \Bd \Lambda(S) \cap M_h$ is an imbedded totally 
geodesic open $(n-1)$-submanifold in $M_h^o$. It is properly imbedded 
since $B(x) \cap \Bd \Lambda(R) \cap \Bd \Lambda(S)$ is compact
for every choice of $B(x)$. 

The above paragraph also shows that 
$\Bd \Lambda(R) \cap \Bd \Lambda(S) \cap M_h$ 
is an open and closed subset of $\Bd \Lambda(R) \cap M_h$. 
Therefore, for components $B$ of 
$\Bd \Lambda(R) \cap M_h$ and $B'$ of $\Bd \Lambda(S) \cap M_h$ where 
$R \not\sim S$, either we have $B = B'$ or $B$ and $B'$ are disjoint.
If $B = B'$, the above paragraph shows that $B$ is a properly imbedded 
totally geodesic open $(n-1)$-submanifold in $M_h$.

We say that a component of $\Bd \Lambda(R) \cap M_h$ is {\sl copied}\/ 
if it equals a component of $\Bd \Lambda(S) \cap M_h$ for some 
$n$-crescent $S$ not equivalent to $R$. Let $c_R$ be the union of all 
copied components of $\Bd \Lambda(R) \cap M_h$. 

\begin{lem}\label{lem:twoco}
Each component of $c_R$ is a properly imbedded totally geodesic 
$(n-1)$-manifold{\rm ,} and equals a component of $\nu_T \cap M_h$ for 
fixed $T${\rm ,} $T \sim R$ and that of $\nu_{T'} \cap M_h$ for fixed 
$T'${\rm ,} $T' \sim S$, where $S$ is not equivalent to $R$. 
\end{lem}
\begin{pf}
From above arguments, we see that given $x$ in a component $C$ of $c_R$, 
and a tiny ball $B(x)$ of $x$, there exists a totally geodesic 
$(n-1)$-ball $P$ with $\delta P \subset \Bd B(x)$ so that a component 
of $B(x) - P$ is included in $T$, $T \sim R$ and the other component 
in $T'$ for $T'$ equivalent to $S$ but not equivalent to $R$.

Since $P$ is connected, $P \subset C$. Let $y$ be another point 
of $C$ connected to $x$ by a path $\gamma$ in $C$, a subset of $M_h$. 
Then we can cover $\gamma$ by a finitely many tiny balls. 
By induction on the number of tiny balls, we see that $y$ 
belongs to $\nu_T \cap M_h$ and $\nu_{T'} \cap M_h$ for fixed 
$T$ and $T'$.  

From this induction, we obtain for each point $y$ of $C$ and 
a tiny ball $B(y)$ of $y$, an existence of an $(n-1)$-ball $P$ satisfying 
\[\Bd P \subset \Bd B(y), P \subset \nu_T \cap \nu_{T'}, \hbox{ and } 
P \subset \Bd \Lambda(R) \cap \Bd \Lambda(S) \cap M_h.\] 
Since $y$ belongs to the interior of $P$ and $P \subset C$ for any 
choice of $y$, $C$ is open in $\nu_T \cap M_h$. Since $C$ is a closed 
subset of $M_h$, $C$ is a component of $\nu_T \cap M_h$. Similarly, 
$C$ is a component of $\nu_{T'} \cap M_h$. 
\end{pf} 

Let $A$ denote $\bigcup_{R \in \cal B} c_R$ where $\cal B$ denotes 
the set of representatives of the equivalence classes of bihedral 
$n$-crescents in $\che M_h$. $A$ is said to be the {\em pre-two-faced 
submanifold arising from bihedral $n$-crescents}. $A$ is a union 
of path-components that are totally geodesic $(n-1)$-manifolds closed 
in $M_h^o$.

\begin{prop}\label{prop:twofacebi}
Suppose that $A$ is not empty. Then $A$ is a properly imbedded submanifold 
of $M_h$ and $p| A$ is a covering map onto a closed totally geodesic 
imbedded $(n-1)$-manifold in $M^o$.   
\end{prop}
\begin{pf}
We follow the argument in Section \ref{sec:hemispheres} somewhat 
repetitively. Every pair of two components $a$ of $c_R$ and $b$ 
of $c_S$ for $n$-crescents $R$ and $S$ where $R, S \in \cal B$, 
are either disjoint or identical. Hence, $A$ is a union of 
disjoint closed path-components that are some components of 
$c_R$ for $R \in \cal B$. This is proved exactly 
as in Section \ref{sec:hemispheres}.

Second, given a tiny ball $B(x)$ of a point $x$ of $M_h$, no more than 
one path-component of $A$ may intersect $\inte B(x)$. Let $a$ be 
a component of $c_R$ intersecting $\inte B(x)$. By Lemma \ref{lem:twoco}, 
$a$ is a component of $\nu_S \cap M_h$ for $S \sim R$ and that of 
$\nu_T \cap M_h$ for $T$ not equivalent to $S$. Furthermore, 
$\nu_S \cap B(x)$ is a compact convex $(n-1)$-ball with boundary in 
$\Bd B(x)$. Since it is connected, $a \cap B(x) = \nu_S \cap B(x)$, and 
$B(x) \cap S$ is the closure of a component $C_1$ of $B(x) - (a \cap B(x))$.
Similarly, $a \cap B(x) = \nu_T \cap B(x)$, and $B(x) \cap T$ is the 
closure of a component $C_2$ of $B(x) - (a \cap B(x))$. Since $S$ and 
$T$ do not overlap, it follows that $C_1$ and $C_2$ are the two 
disjoint components of $B(x) - (a \cap B(x))$.

Suppose that $b$ is a component of $c_U$ for $U \in \cal B$ intersecting 
$\inte B(x)$ also. By Lemma \ref{lem:twoco}, $b$ is a component of 
$\nu_{T'} \cap M_h$ for $T' \sim U$. If the $(n-1)$-ball $b \cap B(x)$  
intersects $C_1$ or $C_2$, then $U \sim S$ or $U \sim T$ and 
$\Lambda(U) = \Lambda(R)$ or $\Lambda(U) = \Lambda(T)$ by Lemma 
\ref{lem:prestab}. The characterization of $B(x) \cap \Lambda(R)$ 
and $B(x) \cap \Lambda(T)$ in Lemma \ref{lem:twoco} implies that $a = b$.
This is absurd. Hence, $b \cap B(x) \subset a \cap B(x)$, and $T'$ 
overlaps with at least one of $S$ or $T$, and we have $a = b$ as above.  

Since given a tiny ball $B(x)$ no more than one distinct path-component 
of $A$ may intersect $\inte B(x)$, $A$ is a properly imbedded closed 
submanifold of $M_h^o$. The rest of the proof of proposition is 
the same as that of Proposition \ref{prop:spmfld}. 
\end{pf}

Let $p: M_h \ra M$ be the covering map. Since $A$ is the deck 
transformation group invariant, we have $A = p^{-1}(p(A))$ and $p|A$ 
covers $p(A)$. The above results show that $p(A)$ is a closed totally 
geodesic manifold in $M^o$.  

\begin{defn}\label{defn:twoface2}
$p(A)$ for the union $A$ of all copied components of $\Lambda(R)$ for 
bihedral $n$-crescents $R$ in $\che{M_h}$ is said to be the {\em two-faced 
$(n-1)$-manifold}\/ of $M$ {\em arising from bihedral $n$-crescents}
(or {\em type II}\/).
\end{defn} 

Each component of $p(A)$ is covered by a component of $A$, i.e.,
a component of $\nu_R \cap M_h$ for some bihderal $n$-crescent $R$. 
Hence, each component of $p(A)$ is covered by open domains in $\bR^n$ 
as in Section \ref{sec:hemispheres}.  

We end with the following observation: 
\begin{prop}\label{prop:nonintf}
Suppose $\che M_h$ includes no hemispheric $n$-crescents and 
$A = \bigcup_{R \in \cal B} c_R$. Then $A$ is disjoint from $R^o$ for 
each $n$-crescent $R$.
\end{prop}
\begin{pf}
The proof is same as that of Proposition \ref{prop:hnonintf}.
\end{pf}

\begin{exmp}\label{exmp:twofacenec}
Finally, we give an example in dimension 2. Let $\vth$ be the projective 
automorphism of $\SI^2$ induced by the diagonal matrix with entries $2$, 
$1$, and $1/2$. Then $\vth$ has fixed points $[\pm 1, 0, 0]$, $[0, \pm 1, 0]$, 
and $[0, 0, \pm 1]$ corresponding to eigenvalues $2, 1, 1/2$. 
Given three points $x, y, z$ of $\SI^2$, we let $\ovl{xyz}$ denote the 
segment with endpoints $x$ and $z$ passing through $y$ if there exists 
such a segment. If $x$ and $y$ are not antipodal, then let $\ovl{xy}$ 
denote the unique minor segment with endpoints $x$ and $y$. We look at 
the closed lune $B_1$ bounded by $\ovl{[0, 0, 1][1, 0, 0][0, 0, -1]}$ and 
$\ovl{[0, 0, 1][0, 1, 0][0, 0, -1]}$, which are to be denoted by $l_1$ 
and $l_2$, and the closed lune $B_2$ bounded by
$\ovl{[1, 0, 0][0, -1, 0][-1, 0, 0]}$ and 
$\ovl{[1, 0, 0][0, 0, 1][-1, 0, 0]}$, which are denoted by $l_3$ and $l_4$.
We consider the domain $U$ given by $U = B_1^o \cup B_2^o \cup 
l_1^o \cup l_4^o -\{ [1, 0, 0], [0, 0, 1]\}$. Since there exists a compact 
fundamental domain of the action of $<\vth>$, $U/<\vth>$ is a compact 
annulus $A$ with totally geodesic boundary. $U$ is the holonomy cover of $A$. 
The projective completion of $U$ can be identified with $B_1 \cup B_2$. 
It is easy to see that $B_1$ is a $2$-crescent with $\alpha_{B_1} = l_2^o$ 
and $\nu_{B_1} = l_1$ and $B_2$ one with $\alpha_{B_2} = l_3^o$ and 
$\nu_{B_2} = l_4$. Also, any other crescent is a subset of $B_1$ or $B_2$.   
Hence $\Lambda(B_1) = B_1$ and $\Lambda(B_2) = B_2$ and the pre-two-faced 
submanifold $L$ equals $\ovl{[1, 0,0][0, 0, 1]}^o$. $L$ covers 
a simple closed curve in $A$ give by $\ovl{[1, 0, 0][0, 0, 1]}^o/<\vth>$.
\end{exmp}

\vfill
\break

\section{Preservation of crescents after splitting}
\label{sec:crescres}
In this section, we consider somwhat technical questions:
How does the $n$-crescents in the completions of the holonomy cover 
of a submanifold become in those of the holonomy cover of an ambient 
manifolds? What happens to $n$-crescents in the completion of 
a manifold when we split the manifold along the two-faced manifolds.
The answers will be that there are one-to-one correspondence: 
Propositions \ref{prop:corrsub}, \ref{prop:hoocorr}, and 
\ref{prop:boocorr}. In the process, we will define splitting manifolds
precisely and show how to contruct holonomy covers of split manifolds. 

Let $M_h$ be the holonomy cover of $M$ with development pair $(\dev, h)$ 
and the group of deck transformations $G_M$. Let $N$ be a submanifold 
(not necessarily simply connected) in $M$ of codimension $0$ with an 
induced real projective structure. Then $p^{-1}(N)$ is a codimension 
$0$ submanifold of $M_h$. Choose a component $A$ of $p^{-1}(N)$. Then 
$A$ is a submanifold in $M_h$ and $p|A$ covers $N$ with the 
deck transformation group $G_{A}$ equal to the group of deck 
transformations of $M_h$ preserving $A$. 

We claim that $A$ is a holonomy cover of $N$ with development pair 
$(\dev|A, h')$ where $h'$ is a composition of the inclusion homomorphism 
and $h: G_M \ra \Aut(\SI^n)$. First, for each closed path in $A$ which 
lifts to one in $N$ obviously has a trivial holonomy (see Section 8.4 
in \cite{RT:94}). Given a closed path with a trivial holonomy, it 
lifts to a closed path in $M_h$ with a base point in $A$. Since $A$ 
is a component of $p^{-1}(N)$, it follows that the closed path is in $A$. 
Therefore, $A$ is the holonomy cover of $N$.

Let us discuss about the Kuiper completion of $A$. The distance 
metric on $A$ is induced from the Riemannian metric on $A$ induced from 
$\SI^n$ by $\dev|A$. The completion of $A$ with respect to the metric 
is denoted by $\che A$ and the set of ideal points $A_\infty$; that is, 
$A_\infty = \che A - A$. Note that $\che A$ may not necessarily 
equal the closure of $A$ in $\che M_h$. A good example may consists of 
the complement of the closure of the positive axis in $\bR^2$ as $A$ 
and $M_h$ as $\bR^2$.

Let $i: A \ra M_h$ be an inclusion map. Then $i$ extends to a distance 
nonincreasing map $\cim: \che A \ra \clo(A) \subset \che M_h$. 

\begin{lem}\label{lem:compl}
\begin{itemize}
\item[(i)] $\cim^{-1}(\hideal{M})$ is a subset of $A_\infty$.
\item[(ii)] If $A$ is closed as a subset of $M_h${\rm ,} then 
$\cim(A_\infty) \subset \hideal{M}$. Thus, in this case, 
$\cim(A_\infty) = \hideal{M}$.
\item[(iii)] Let $P$ be a submanifold in $A$ with convex interior $P^o$.
Then the closure $P'$ of $P$ in $\che A$ maps isometrically to the 
closure $P''$ of $P$ in $\che M_h$ under $\cim$. Here $P'$ and $P''$ 
are convex. 
\item[(iv)] $\cim$ maps $P' \cap A_\infty$ homeomorphic onto 
$P''\cap \hideal{M}$. 
\end{itemize}
\end{lem}
\begin{pf} 
(i) If $x$ is a point of $\cim^{-1}(\hideal{M})$, then $x$ 
does not belong to $A$ since otherwise $\cim(x) = i(x) \in M_h$. 

(ii) Suppose not. Then there exists a point $x$ in $M_h$ such that 
$x = \cim(y)$ for $y \in A_\infty$. There exists a sequence of points 
$y_i \in A$ converging to $y$. The sequence of points $i(y_i)$ converges 
to a point $x$ since $i$ is distance decreasing. Since $A$ is closed, 
this means $x \in A$ and $y = x$. This is a contradiction. 

(iii) Since $i| P^o$ is an isometry with respect to $\bdd|A$ and 
$\bdd$ on $M_h$, the third part follows.

(iv) Let $K$ be the inverse image of $P'' \cap \hideal{M}$ under 
$\cim| P'$. By (i), we obtain $K \subset A_\infty$. By (ii), we see 
$\cim(P' \cap A_\infty) \subset P'' \cap \hideal{M}$. 
\end{pf}

Let $R$ be an $n$-crescent in $\che M_h$, and consider a submanifold 
$R' = R \cap M_h$ with convex interior $R^o$. If $R'$ is a subset of 
a submanifold $A$ of $M_h$, then the above lemma shows that the closure 
$R''$ of $R'$ in $\che A$ is isometric to $R$ under $\cim$. By the above 
lemma, we obtain that $R''$ is also a crescent with $\alpha_{R''} 
= \cim^{-1}(\alpha_R)$, and $\nu_{R''} = \cim^{-1}(\nu_R)$. Moreover, 
if $R$ is bihedral (resp. hemispheric), then $R''$ is bihedral (resp. 
hemispheric). 

Conversely, let $R$ be an $n$-crescent in $\che A$. By Lemma \ref{lem:compl},
$\cim| R: R \ra \cim(R)$ is an imbedding, and the closure of $i(R \cap A)$ 
equals $\cim(R)$ and is a convex $n$-ball. By Lemma \ref{lem:compl}, 
$\cim(R)$ is an $n$-crescent, which is bihedral (resp. hemispheric) 
if $R$ is bihedral (resp. hemispheric). (Note that $\alpha_{\cim(R)} 
= \cim(\alpha_R)$ and $\nu_{\cim(R)} = \cim(\nu_R)$.)

Thus, we have proved.
\begin{prop}\label{prop:corrsub}
There exists a one to one correspondence of all bihedral $n$-crescents 
in $\che A$ and those in $\clo(A)$ in $\che M_h$ given by
$R \leftrightarrow R'$ for a bihedral $n$-crescent $R$ in $\che A$ 
and $R'$ one in $\clo(A)$ if and only if $R^o = R^{\prime o}$.
The same statement holds for hemispherical $n$-crescents.
\end{prop}

We give a precise definition of splitting. Let $N$ be a real 
projective $n$-manifold with a properly imbedded $(n-1)$-submanifold 
$A$. We take an open regular neighborhood $N$ of $A$. Then $N$ is an 
$I$-bundle over $A$. Let us enumerate components of $A$ by $A_1, \dots, 
A_n, \dots$ and corresponding components of $N$ by $N_1, \dots, N_n, \dots$ 
which are regular neighborhoods of $A_1, \dots, A_n, \dots $ respectively.
(We do not require the number of components to be finite.) 

For an $i$, $N_i$ is an $I$-bundle over $A_i$. By parameterizing 
each fiber by a real line, $N_i$ becomes a vector bundle 
over $A$ with a flat linear connection. We see that there is a subgroup
$G_i$ of index at most two in $\pi_1(A_i)$ with trivial holonomy. The 
single or double cover $\tilde N_i$ of $N_i$ corresponding to $G_i$ is 
a product $I$-bundle over $\tilde A_i$ the cover of $A_i$ corresponding 
to $G_i$, considered as a submanifold of $\tilde N_i$.

Since $N_i$ is a product or twisted $I$-bundle over $A_i$, $N_i - A_i$ has 
one or two components. If $N_i - A_i$ has two components, then we take 
the closure of each components in $N_i$ and take their disjoint 
union $\hat N_i$ which has a natural inclusion map 
$l_i: N_i-A_i \ra \hat N_i$. If $N_i - A_i$ has one component, then 
$\pi_1(A_i)$ has an index two subgroup. Take the double cover 
$\tilde N_i$ of $N_i$ so that it is now a product $I$-bundle over 
the corresponding double cover $\tilde A_i$ in $\tilde N_i$. 
Then $N_i - A_i$ lifts and imbeds onto a component of $\tilde N_i - 
\tilde A_i$. We denote by $\hat N_i$ the closure of this component in 
$\tilde N_i$. There is a natural lift $l_i: N_i - A_i \ra \hat N_i$, which 
is an imbedding. After we do this for each $i$, $i=1, \dots, n, \dots$, 
we identify $N - A$ and the disjoint union $\coprod_{i = 1}^n \hat N_i$
of all $N_i$ by the maps $l_i$. Then the resulting manifold $M$ is said 
to be obtained from $N$ by {\em splitting along} $A$. Obviously, 
$M$ is compact if $N$ is compact and $M$ has totally geodesic 
boundary if $A$ is totally geodesic. When $A$ is not empty, $M$ is 
said to be the {\em split manifold}\/ obtained from $N$ along $A$ 
(this is just for terminological convenience). Also, we see that 
for each component of $A$, we get either two copies or a double 
cover of the component of $A$ in the boundary of the split manifold $M$. 
They are newly created by splitting. There is a natural quotient map 
$q:M \ra N$ by identifying these new faces to $A$, i.e., 
$q|q^{-1}(A): q^{-1}(A) \ra A$ is a two-to-one covering map.


Now, suppose that $N$ is compact and $N_h$ its holonomy cover with 
development pair $(\dev, h)$. We let $G_N$ the group of deck 
transformations of the covering map $p:N_h \ra N$. Suppose that $A$ has 
finitely many components. Then $p^{-1}(A)$ is a properly imbedded 
$(n-1)$-manifold in $N_h$. If one splits $N_h$ along $p^{-1}(A)$, then 
it is easy to see that the split manifold $M'$ covers the manifold $M$ 
of $N$ split along $A$ with covering map $p'$ obtained from extending $p$. 
However, $M'$ may not be connected; for each component $M_i$ of $M$, we 
choose a component $M'_i$ of $M'$ covering that component. $M_i$ contains 
exactly one component $P_i$ of $N - A$, and $M'_i$ includes exactly one 
component $P'_i$ of $N_h - p^{-1}(A)$ as a dense open subset. Thus, 
$\coprod_{i=1}^n P'_i$ covers $\coprod_{i=1}^n P_i$ and 
$\coprod_{i=1}^n M'_i$ covers $M$. (Note that covers in this section 
are not necessarily connected ones.)

\begin{rem}\label{rem:orient}
The submanifold $p^{-1}(A)$ is orientable since great $(n-1)$-spheres 
in $\SI^n$ are orientable and $\dev$ maps $p^{-1}(A)$ into a great 
$(n-1)$-sphere as an immersion. Thus, there are no twisted 
$I$-bundle neighborhoods of components of $p^{-1}(A)$ as $N_h$ is 
orientable also.
\end{rem}

The developing map $\dev|P_i$ uniquely extends to a map from $M'_i$ 
as an immersion for each $i$; we denote by 
$\dev':\coprod_{i=1}^n M'_i \ra \SI^n$ the map obtained this way. 
The deck transformation group $G_i$ of the covering map 
$p|P'_i: P'_i \ra P_i$ equals the subgroup of $G_N$ consisting of 
deck transformations of $N_h$ preserving $P'_i$. It is easy to see 
that the action of $G_i$ naturally extends to $M'_i$ and becomes the 
group of deck transformations of the covering map $p|M'_i: M'_i \ra M_i$. 
For each $i$, we define $h_i:G_i \ra \Aut(\SI^n)$ by $h \circ l_i$ 
where $l_i:G_i \ra G_N$ is the homomorphism induced from the 
inclusion map. We choose base points of $P_i$ and $P'_i$ respectively, 
and this gives us a map from the fundamental group $\pi_1(P_i)$ to 
$G_i$. We define the holonomy homomorphism $h'_i:\pi_1(P_i) \ra 
\Aut(\SI^n)$ by first mapping to $G_i$, and then followed by $h_i$. 

We claim that $P'_i$ is the holonomy cover of $P_i$ with development
pair $(\dev'|P'_i, h_i)$; that is, for each $P_i$ and $P'_i$ that 
$\ker h'_i$ equals the fundamental group of $P'_i$: Since $N_h$ 
is the holonomy cover of $N$, $h$ is injective; for each element 
of $G_i$ other than the identity maps to a nontrivial element in 
$\Aut(\SI^n)$. Hence $\ker h'_i$ equals the kernel of the map 
$\pi_1(P_i) \ra G_i$ given by path lifting to $P'_i$, and $\ker h'_i$ 
equals $\pi_1(P'_i)$. 

Moreover, since $M'_i$ is obtained from $P'_i$ by attaching boundary, 
it follows that $M'_i$ is the holonomy cover of $M_i$ with development 
pair $(\dev'|M'_i, h_i)$. We will sometimes say that the disjoint 
union $\coprod_{i=1}^n M'_i$ is a {\em holonomy cover} of 
$M = \coprod_{i=1}^n M_i$. 

Suppose that there exists a nonempty pre-two-faced submanifold $A$ 
of $N_h$ arising from hemispheric $n$-crescents. Then we can split 
$N$ by $p(A)$ to obtain $M$ and $N_h$ by $A$ to obtain $M'$, and 
$M'$ covers $M$ under the extension $p'$ of the covering map 
$p: N_h \ra N$. 

We claim that the collection of hemispheric $n$-crescents in $\che N_h$ 
and the completion $\che M'$ of $M'$ are in one to one correspondence. 
Let $q: M' \ra N_h$ denote the natural quotient map identifying the 
newly created boundary components which restricts to the inclusion map 
$N_h - A \ra N_h$. We denote by $A'$ the set $q^{-1}(A)$, which 
are newly created boundary components of $M'$. Let $\che M'$ denote 
the projective completion of $M'$ with the metric $\bdd$ extended from 
$N_h-A$. Then as $q$ is distance decreasing, $q$ extends to a map 
$\che q:\che M' \ra \che N_h$ which is one-to-one and onto on 
$M' - A' \ra N_h-A$.

\begin{lem}\label{lem:techn}
$\che q$ maps $A'$ to $A$, $M'$ to $N_h$, and $M'_\infty$ to 
$\hideal{N}$.
{\rm (} Which implies that $\che q^{-1}(A) = A'$, 
$\che q^{-1}(N_h) = M'$, and 
$\che q^{-1}(\hideal{N}) = M'_\infty$ {\rm .)} 
\end{lem}
\begin{pf}
We obviously have $\che q(A') = q(A') = A$ and 
$\che q(M') = q(M') = N_h$.

If a point $x$ of $\che M'$ is mapped to that of $A$, then 
let $\gamma$ be a path in $N_h - A$ ending at $\che q(x)$ in $A$.
Then we may lift $\gamma$ to a path $\gamma'$ in $M'- A'$.
$\che q(x)$ has a small compact neighborhood $B$ in $N_h$ where 
$\gamma$ eventually lies in, and as the closure $B'$ of 
a component of $B -A$ is compact, there exists a compact 
neighborhood $B''$ in $M'$ mapping homeomorphic to $B'$ under 
$q$ and $\gamma'$ eventually lies in $B''$. This means that $x$ 
lies in $B''$ and hence in $M'$. As $q^{-1}(A) = A'$, $x$ lies 
in $A'$. Thus, $\che q^{-1}(A) = A'$ and points of $M'_\infty$
cannot map to a point of $A$. 

Using path lifting argument, we may show that 
$\che q(M'_\infty)$ is a subset of $A \cup \hideal{N}$
as $\che q| M' - A' \ra N_h -A$ is a local $\bdd$-isometry.
Hence, it follows that $\che q(M'_\infty) \subset \hideal{N}$.
\end{pf}

First, consider the case when $\che N_h$ includes a hemispheric 
$n$-crescent $R$. As by Proposition \ref{prop:hnonintf}, 
$R^o$ is a subset of $N_h - A$, $R^o$ is a subset of $M'$. 
The closure $R'$ of $R^o$ in $\che M'$ is naturally an $n$-hemisphere 
as $\dev|R^o$ is an imbedding onto an open $n$-hemisphere 
in $\SI^n$. As $\che q$ is a $\bdd$-isometry restricted to $R^o$, 
it follows that $\che q|R': R' \ra R$ is an imbedding. 

Lemma \ref{lem:techn} shows that $(\che q|R')^{-1}(\alpha_R)$ is 
a subset of $M'_\infty$. Thus, $R'$ includes an open $(n-1)$-hemisphere 
in $\delta R' \cap M'_\infty$, which shows that $R'$ is 
a hemispheric $n$-crescent. ($\delta R'$ can not belong to 
$M'_\infty$ by Lemma \ref{lem:include}.) 

Now if an $n$-crescent $R$ is given in $\che M'$, then we have 
$R^o \subset M' - A'$, and $\che q(R)$ is obviously an $n$-hemisphere 
as the closure $R'$ of $R^o$ in $\che N_h$ is an $n$-hemisphere and 
equals $\che q(R)$. Since $\che q(\alpha_R)$ is a subset of 
$\hideal{N}$ by Lemma \ref{lem:techn}, $\che q(R)$ is a hemispheric
$n$-crescent.

\begin{prop}\label{prop:hoocorr}
There exists a one-to-one correspondence between all hemispheric 
$n$-crescents in $\che N_h$ and those of $\che M'$ by the 
correspondence $R \leftrightarrow R'$ if and only if 
we have $R^o = R^{\prime o}$.
\end{prop}

\begin{cor}\label{cor:hexists}
If $\che N_h$ includes a hemispheric $n$-crescent, then the projective
completion of the holonomy cover of at least one component of the split 
manifold $M$ along the two-faced submanifold, also includes 
a hemispheric $n$-crescent.
\end{cor}
\begin{pf} 
Let $M_i$ be the components of $M$ and $M'_i$ its holonomy cover as 
obtained earlier in this section; let $P_i$ be the component of 
$N - p(A)$ in $M_i$ and $P'_i$ that of $N_h - A$ in $M'_i$ so that 
$P'_i$ covers $P_i$. We regard two components of $N_h - A$ to be 
equivalent if there exists a deck transformation of $N_h$ mapping one 
to the other. Then $P_i$ is a representative of an equivalence class 
${\cal A}_i$. As the deck transformation group is transitive in 
an equivalence class ${\cal A}_i$, we see that given two elements 
$P^a_i$ and $P^b_i$ in ${\cal A}_i$, the components $M^a_i$ and $M^b_i$ 
of $M$ including them respectively are projectively isomorphic 
as the deck transformation sending $P^a_i$ to $P^b_i$ extends to 
a projective map $M^a_i \ra M^b_i$, and hence to a quasi-isometry 
$\che M^a_i \ra \che M^b_i$. Thus, if no $\che M'_i$ includes 
a hemispheric $n$-crescent, then it follows that $\che M'$ do not 
also. This contradicts Proposition \ref{prop:hoocorr}. 
\end{pf}

Now, we suppose that $\che N_h$ includes no hemispheric $n$-crescents 
$R$ but includes some bihedral $n$-crescents. Let $A$ be the 
pre-two-faced submanifold of $N_h$ arising from bihedral $n$-crescents.
We split $N$ by $p(A)$ to obtain $M$ and $N_h$ by $A$ to obtain $M'$, 
and $M'$ covers $M$ under the extension $p'$ of the covering map 
$p: N_h - A \ra N - p(A)$. 

By same reasonings as above, we obtain 
\begin{prop}\label{prop:boocorr}
There exists a one-to-one correspondence between all bihedral 
$n$-crescents in $\che N_h$ and those of $\che M'$ by correspondence 
$R \leftrightarrow R'$ if and only if $R^o = R^{\prime o}$.
\end{prop}

\begin{cor}\label{cor:bhexist}
If $\che N_h$ includes a bihedral $n$-crescent, then the projective
completion of the holonomy cover of at least one component of the split 
manifold $M$ along the two-faced submanifold{\rm ,} also includes 
a bihedral $n$-crescent.
\end{cor}

\vfill
\break

\section{The proof of the Main theorem}
\label{sec:pfmain}
In this section we prove Theorem \ref{thm:main} using the previous
three-sections, in a more or less straighforward manner. We will
start with hemispheric crescent case and then the bihedral case.

\begin{defn}\label{defn:concaff}
A {\sl concave affine manifold $M$ of type I}\/ is a compact real projective 
manifold such that its completion $\che M_h$ of the holonomy cover $M_h$ 
of $M$ is an $n$-crescent that is an $n$-hemisphere, or a finite disjoint 
union of such real projective manifolds. A {\sl concave affine manifold 
$M$ of type II}\/ is a compact real projective manifold with concave 
or totally geodesic boundary so that its completion $\che M_h$ of 
the holonomy cover $M_h$ equals $\Lambda(R)$ for an $n$-crescent $R$ 
in $\che M_h$ and no $n$-crescent is an $n$-hemisphere, or is a finite 
disjoint union of such real projective manifolds. We allow $M$ to have 
nonsmooth boundary that is concave, i.e., not necessarily totally geodesic. 
\end{defn}

Concave affine manifolds admit natural affine structures:
If $M$ is a concave affine manifold of type I, from the properties 
proved in the above Section \ref{sec:crescres}, $\dev(\che M_h)$ equals 
an $n$-hemisphere $H$. Since the holonomy group acts on $H$, it restricts 
to an affine transformation group in $H^o$. Since a projective 
transformation acting on an affine patch is affine, the interior $M^o$ 
has a compatible affine structure. If $M$ is one of the second type, 
then for each bihedral $n$-crescent $R$, $\dev$ maps $R \cap M_h$ into 
the interior of an $n$-hemisphere $H$ (see Section \ref{sec:bitwof}). 
Hence, it follows that $\dev$ maps $M_h$ into $H^o$. Since given 
a deck transformation $\vth$, we have $\vth(\Lambda(R)) = \Lambda(\vth(R)) 
= \che M_h$, we obtain $\inte \Lambda(\vth(R)) \cap \inte \Lambda(R) 
\cap M_h \ne \emp$ and $R \sim \vth(R)$ by Lemma \ref{lem:prestab}. 
This shows that $\Lambda(R) =\Lambda(\vth(R)) = \vth(\Lambda(R))$ and 
$\delta_\infty \Lambda(R) = \delta_\infty \Lambda(\vth(R)) = 
\vth(\delta_\infty \Lambda(R))$ for each deck transformation $\vth$
by equation \ref{eqn:lambda2}. 
Since $\dev(\delta_\infty \Lambda(R))$ is a subset of a unique great 
sphere $\SI^{n-1}$, it follows that $h(\vth)$ acts on $\SI^{n-1}$ and 
since $\dev(M_h)$ lies in $H^o$, the holonomy group acts on $H^o$. 
Therefore, $M$ has a compatible affine structure. 

If $M$ is a concave affine manifold of type I, then $\delta M$ is 
totally geodesic since $M_h = R \cap M_h$ for an hemispheric $n$-crescent 
$R$ and $\delta M_h = \nu_R \cap M_h$. If $M$ is one of type II, then 
$\delta M$ is concave, as we said in the definition above.

For the purpose of the following lemma, we also define two $n$-crescents 
$S$ and $T$, hemispheric or bihedral, to be {\em equivalent}\/ if there 
exists a chain of $n$-crescents $T_1, T_2, \dots, T_n$ so that $S = T_1$ 
and $T = T_i$ and $T_i$ and $T_{i+1}$ overlaps for each $i=1, \dots, n-1$. 
We will use this definition in this section only.

\begin{lem}\label{lem:seqcres}
Let $x_i$ be a sequence of points of $M_h$ converging to a point $x$ of 
$M_h${\rm ,} and $x_i \in R_i$ for $n$-crescents $R_i$ for each $i$. 
Then for any choice of an integer $N${\rm ,} we have $R_i \sim R_j$ for 
infinitely many $i, j \geq N$. Furthermore{\rm ,} if each $R_i$ is an 
$n$-hemisphere{\rm ,} then $R_i = R_j$ for infinitely many $i, j \geq N$.
Finally $x$ belongs to an $n$-crescent $R$ for $R \sim R_i$ for each 
$i \geq N$.
\end{lem}
\begin{pf} 
Let $B(x)$ be a tiny ball of $x$. Assume $x_i \in \inte B(x)$ for each $i$.
We can choose a smaller $n$-crescent $S_i$ in $R_i$ so that $x_i$ now 
belongs to $\nu_{S_i}$ with $\alpha_{S_i}$ subset of $R_i$ 
as $R_i$ are geometrically ``simple'', i.e., a convex $n$-bihedron or 
an $n$-hemisphere 

Since $B(x)$ cannot be a subset of $S_i$, $S_i \cap B(x)$ is the component of 
the closure of $B(x) - a_i$ for $a_i = \nu_{S_i} \cap B(x)$ an $(n-1)$-ball 
with boundary in $\Bd B(x)$. Let $v_i$ be the outer-normal vector at $x_i$ 
to $\nu_{S_i}$ for each $i$. Choose a subsequence $i_j$, with $i_1 = N$, 
of $i$ so that the sequence $v_{i_j}$ converges to a vector $v$ at $x$. 
A generalization of Lemma 3 of the Appendix of \cite{cdcr1} shows that 
there exists a common open ball $\cal P$ in $R_{i_j}$ for $j \geq K$ for 
some integer $K$. (The proofs are identical.) Thus $S_{i_j} \sim S_{i_k}$ 
for $j, k \geq K$. Since $R_{i_j} \sim S_{i_j}$ as $S_{i_j}$ is a subset of 
$R_{i_j}$, we have that $R_{i_j} \sim R_{i_k}$ for $j, k \geq K$.

We assume that each sequence of $\dev(S_{i_j})$ and  
$\dev(\clo(\alpha_{S_{i_j}}))$ respectively converge to compact sets in 
$\SI^n$ under the Hausdorff topology by choosing a subsequence if necessary.
As $S_{i_j}$ forms a cored sequence, Lemma 4 of the Appendix of \cite{cdcr1}
shows that there exists a convex $n$-ball $S$ containing $\cal P$,  
such that $\dev(S)$ equals the limit of $S_{i_j}$. There exists a set 
$\alpha$ in $S$ which maps to the limit $\alpha_\infty$ of 
$\dev(\clo(\alpha_{S_{i_j}}))$, and as $\clo(\alpha_{S_{i_j}})$ forms 
a subjugated ideal sequence, it follows that $\alpha$ lies in $\hideal{M}$. 
As $\clo(\alpha_{S_{i_j}})$ includes an $(n-1)$-hemisphere, so does
$\alpha_\infty$. Thus $\dev(S)$ is either an $n$-bihderon or 
an $n$-hemisphere. Therefore, $S$ is an $n$-crescent. Obviously, 
$S \sim S_{i_j}$ for $j \geq K$.

Since the sequence consisting of $x_i$ is also a subjugated sequence
of $S_i$ with $\dev(x_i) \ra \dev(x)$, we see that there exists 
a point $y$ in $S$ mapping to $\dev(x)$ under $\dev$. 
As $B(x)$ and $S$ overlap, $\dev| B(x) \cup S$ is an imbedding.
This shows that $x = y$. 
\end{pf}

We begin the proof of the Main Theorem \ref{thm:main}. Actually, what 
we will be proving the following theorem, which together with Theorem 
\ref{thm:n-1conv} implies Theorem \ref{thm:main}.
\begin{thm}\label{thm:main2} 
Suppose that $M$ is a compact real projective $n$-manifold with totally
geodesic or empty boundary, and that $M_h$ is not real projectively 
diffeomorphic to an open $n$-hemisphere or $n$-bihedron. Then 
the following statements hold\/{\rm :}
\begin{itemize}
\item If $\che M$ includes a hemispheric $n$-crescent, then $M$ includes 
a compact concave affine $n$-submanifold $N$ of type I or $M^o$ 
includes the two-faced $(n-1)$-submanifold arising from hemispheric 
$n$-crescent. 
\item If $\che M$ includes a bihedral $n$-crescent, then $M$ 
includes a compact concave affine $n$-submanifold $N$ of type II or 
$M^o$ includes the two-faced $(n-1)$-submanifold arising from 
bihedral $n$-crescent. 
\end{itemize}
\end{thm} 

First, we consider the case when $\che M_h$ includes an $n$-crescent 
$R$ that is an $n$-hemisphere. Suppose that there is no copied 
component of $\nu_T \cap M_h$ for every hemispheric $n$-crescent $T$. 
Recall from Section \ref{sec:hemispheres} that either $R = S$ or 
$R$ and $S$ are disjoint for every pair of hemispheric $n$-crescents 
$R$ and $S$. 


Let $x \in M_h$ and $B(x)$ the tiny ball of $x$. Then only finitely 
many distinct hemispheric $n$-crescents intersect a compact neighborhood
of $x$ in $\inte B(x)$. Otherwise, there exists a sequence of points $x_i$ 
converging to a point $y$ of $\inte B(x)$, where $x_i \in \inte B(x)$ and 
$x_i \in R_i$ for hemispheric $n$-crescents $R_i$ for mutually distinct 
$R_i$, but Lemma \ref{lem:seqcres} contradicts this. 

Consider $R \cap M_h$ for a hemispheric $n$-crescent $R$. Since $R$ is 
a closed subset of $\che M_h$, $R\cap M_h$ is a closed subset of $M_h$.
Let $A$ be the set $\bigcup_{R \in \cal H} R \cap M_h$. Then $A$ is 
a closed subset of $M_h$ by above. 

Since $R \cap M_h$ is a submanifold for each $n$-crescent $R$, $A$ is 
a submanifold of $M_h$, a closed subset. Since the union of all hemispheric
$n$-crescents $A$ is deck transformation group invariant, we have 
$p^{-1}(p(A)) = A$. The above results show that $p| A$ is a covering map 
onto a compact submanifold $N$ in $M$, and $p| R \cap M_h$ is a covering map 
onto a component of $N$ for each hemispheric $n$-crescent $R$. 


Since the components of $A$ are locally finite in $M_h$, it follows 
that $N$ has only finitely many components.
Let $K$ be a component of $N$. In the beginning of Section 
\ref{sec:crescres}, we showed that $R \cap M_h$ is a holonomy cover 
of $K$. Let $\che K$ be the projective completion of $R \cap M_h$. 
The closure of $R \cap M_h$ in $\che K$ is a hemispheric $n$-crescent 
identical with $\che K$ as a set by Proposition \ref{prop:corrsub}. 
Hence, $K$ is a concave affine manifold of type I. 


If there is a copied component of $\nu_T \cap M_h$ for some 
hemispheric $n$-crescent $T$, Proposition \ref{prop:spmfld} implies 
the Main theorem.

Now, we assume that $\che M_h$ includes only $n$-crescents that are 
$n$-bihedrons. Suppose that there is no copied component of 
$\Bd \Lambda(T) \cap M_h$ for every bihedral $n$-crescents $T$, 
$T \in \cal B$. Then either $\Lambda(R) = \Lambda(S)$ or $\Lambda(R)$ 
and $\Lambda(S)$ are disjoint for $n$-crescents $R$ and $S$, $R, S \in 
\cal B$, by Lemma \ref{lem:prestab}.

Using this fact and Lemma \ref{lem:seqcres}, we can show similarly 
to the proof for the hemispheric $n$-crescents that 
$A = \bigcup_{R \in \cal B} \Lambda(R) \cap M_h$ is closed:
We showed that $\Lambda(R) \cap M_h$ is a closed subset of $M_h$ in 
Section \ref{sec:bitwof}. For each point $x$ of $M_h$ and a tiny ball 
$B(x)$ of $x$, there are only finitely many mutually distinct 
$\Lambda(R_i)$ intersecting a compact neighbohood of $x$ in 
$\inte B(x)$ for $n$-crescents $R_i$: Otherwise, we get a sequence 
$x_i$, $x_i \in \inte B(x)$, converging to $y$, $y \in \inte B(x)$, 
so that $x_i \in \Lambda(R_i)$ for $n$-crescents $R_i$ with mutually 
distinct $\Lambda(R_i)$, i.e., $R_i$ is not equivalent to $R_j$ 
whenever $i \ne j$. Then $x_i \in S_i$ for an $n$-crescent $S_i$ 
equivalent to $R_i$. Lemma \ref{lem:seqcres} implies $S_i \sim S_j$ 
for infinitely many $i, j \geq N$. Since $S_i \sim R_i$, this 
contradicts the fact that $\Lambda(R_i)$ are mutually distinct.

The subset $A$ is a submanifold since each $\Lambda(R) \cap M_h$ is 
one for each $R$, $R \in \cal B$. Similarly to the hemisphere case, 
since $p^{-1}(p(A)) = A$, we obtain that $p| A$ is a covering map 
onto a compact submanifold $N$ in $M$, and $p|\Lambda(R) \cap M_h$, 
$R \in \cal B$, is a covering map onto a component of $N$. 
$N$ has finitely many components since the components of $A$ are 
locally finite in $M_h$ by the above paragraph.


Let $K$ be the component of $N$ that is the image of $\Lambda(R) 
\cap M_h$ for $R$, $R \in \cal B$. As in the beginning of Section 
\ref{sec:crescres}, $\Lambda(R) \cap M_h$ is a holonomy cover of $K$. 
Let $\che K_h$ denote the projective completion of 
$\Lambda(R) \cap M_h$. For each crescent $S$, $S \sim R$, 
the closure $S'$ of $S \cap M_h$ in $\che K_h$ is an $n$-crescent 
(see Section \ref{sec:crescres}). It follows that each point $x$ of 
$\Lambda(R) \cap M_h$ is a point of a crescent $S'$ in $\che K$ 
equivalent to the crescent $R'$, the closure of $R \cap M_h$ 
in $\che K$. Since $\che K$ thus is the union of equivalent bihedral 
$n$-crescents, $N$ is a concave affine manifold of type II. 

When there are copied components of $\Bd \Lambda(R) \cap M_h$, 
then Proposition \ref{prop:twofacebi} completes the proof of 
the Main theorem. \qed

\vfill
\break

\section{The proof of the Corollary \ref{cor:spl1}} 
\label{sec:corspl1}

In this section, we will prove Corollary \ref{cor:spl1}. The basic 
tools are already covered in previous three sections. As before, we 
study hemispheric case first, and the proof of the bihedral case is 
entirely similar but will be spelled out. We will also give a proof 
of Corollary \ref{cor:cdcr} here.

We give a cautionary note before we begin: We will assume that 
$M$ is not $(n-1)$-convex, and so $M_h$ is not projectively 
diffeomorphic to an open $n$-bihedron or an open $n$-hemisphere,
so that we can apply various results in Sections 5 -- 8, such as 
the intersection properties of hemispheric and bihedral $n$-crescents. 
We will carry out various decompositions of $M$ in this section. Since 
in each of the following step, the result are real projective manifolds 
with nonempty boundary if a nontrivial decomposition had occurred, it 
follows that their holonomy covers are not projectively diffeomorphic 
to open $n$-bihedrons and open $n$-hemispheres. So our theory in 
Sections 5 -- 8 continues to be applicable.

Also, in this section, covering spaces need not be connected. This 
only complicates the matter of identifying the fundamental groups with 
the deck transformation groups, where we will be a little cautious. Note 
that even for disconnected spaces we can define projective completions 
as long as immersions to $\SI^n$, i.e., developing maps, are defined 
since we can always pull-back the metrics in that case.

We show a diagram of manifolds that we will be obtaining in 
the construction. The ladder in the first row is continued to 
the next one. Consider it as a one continuous ladder.
\begin{equation*}
\begin{array}{*{5}{c}}
M &\Rar_{p(A_1)} &M^{\sf s} & \Rar & N \coprod K \\
\uar \; p & {} & \uar & {}& \uar \\
M_h &\Rar_{A_1} & M_h^{\sf s} &\Rar & 
N_h \coprod \bigcup_{R \in \cal H} R \cap N_h 
\end{array}
\end{equation*}
\begin{equation}
\begin{array}{*{4}{c}}\label{eqn:arrows}
\Rar_{p(A_2)} & N^{\sf s} \coprod K & \Rar & S \coprod T \coprod K\\
{} &\uar &{} &\uar \\
\Rar_{A_2} & N_h^{\sf s} \coprod \bigcup_{R \in \cal H} R \cap N_h
&\Rar & S_h \coprod
\bigcup_{R \in \cal B} \Lambda(R) \cap N_h^{\sf s} \coprod 
\bigcup_{R \in \cal H} R \cap N_h
\end{array}
\end{equation}
where the notation $\Rar_{A}$ means to split along a submanifold $A$
if $A$ is compact and means to split and take appropriate components 
to obtain a holonomy cover if $A$ is noncompact, $\Rar$ means to 
decompose and to take appropriate components, $\coprod$ means a 
disjoint union and other symbols will be explained as we go along. 
Note that when any of $A_1, K, A_2, T$ is empty, then the operation 
of splitting or decomposition does not take place and the next 
manifolds are identical with the previous ones. For convenience, 
we will assume that all of them are not empty in the proof. 

Let $M$ be a compact real projective $n$-manifold with totally 
geodesic or empty boundary. Suppose that $M$ is not $(n-1)$-convex, 
and we will now be decomposing $M$ into various canonical pieces. 

Since $M$ is not $(n-1)$-convex, $\che M_h$ includes an $n$-crescent
(see Theorem \ref{thm:n-1conv}). By Theorem \ref{thm:main2}, $M$ has a 
two-faced $(n-1)$-manifold $S$, or $M$ includes a concave affine manifold.

Suppose that $\che M_h$ has a hemispheric $n$-crescent, and that $A_1$ is 
a pre-two-faced submanifold arising from hemispheric $n$-crescents.
(As before $A_1$ is two-sided.) Let $M^{\sf s}$ denote the result of 
the splitting of $M$ along $p(A_1)$, and $M'$ that of $M_h$ along $A_1$, 
and $A'_1$ the boundary of $M'$ corresponding to $A_1$, ``newly created from 
splitting.'' We know from Section \ref{sec:bitwof} that there exists 
a holonomy cover of $M^{\sf s}$ that is a union of suitable components 
of $M'$. This completes the construction of the first column of 
arrows in equation \ref{eqn:arrows}.

We now show that $M^{\sf s}$ now has no two-faced submanifold. Let $\che M'$ 
denote the projective completion of $M'$. Suppose that two hemispheric 
$n$-crescents $R$ and $S$ in $\che M'$ meet at a common component $C$ of 
$\nu_R \cap M'$ and $\nu_S \cap M'$ and that $R$ and $S$ are not 
equivalent. Proposition \ref{prop:spmfld} applied to $M'$ shows that 
$C \subset M^{\prime o}$; in particular, $C$ is disjoint from $A'_1$. 

There exists a map $\hat q: \che M' \ra \che M_h$ extending the 
quotient map $q:M' \ra M_h$ identifying newly created split faces 
in $M'$. There exist hemispheric $n$-crescents $\hat q(R)$ and 
$\hat q(S)$ in $\che M_h$ with same interior as $R^o$ and $S^o$ 
included in $M_h - A_1$ by Proposition \ref{prop:hoocorr}.

Since $\nu_R \cap \nu_S \cap M'$ belongs to $M' - A'_1 = M_h - A_1$, 
it follows that $q(\nu_R)$ and $q(\nu_S)$ meet in $M_h - A_1$. 
Since obviously $q(\nu_R) \subset R'$ and $q(\nu_S) \subset S'$, 
we have that $R'$ and $S'$ meet in $M_h - A_1$. However, since 
$R'$ and $S'$ are hemispheric $n$-crescents in $\che M_h$, 
$q(C)$ is a subset of the pre-two-faced submanifold $A_1$. This is 
a contradiction. Therefore, we have either $R = S$ or 
$R \cap S = \emp$ for $n$-crescents $R$ and $S$ in $\che M'$. 
Finally, since $\che M_h^{\sf s}$ is a subset of $\che M'$ obviously, 
we also have $R = S$ or $R \cap S = \emp$ for $n$-crescents $R$ 
and $S$ in $\che M_h^{\sf s}$.
  
The above shows that $M^{\sf s}$ has no two-faced submanifold 
arising from hemispheric $n$-crescents. Let $\cal H$ denote the set 
of all hemispheric $n$-crescents in $\che M_h^{\sf s}$. As in Section 
\ref{sec:pfmain}, $p|\bigcup_{R\in \cal H} R \cap M_h^{\sf s}$ 
is a covering map to the concave affine manifold $K$ of type I. 
Since any two hemispheric $n$-crescents are equal or disjoint, 
it is easy to see that $p^{-1}(K^o) = \bigcup_{R \in \cal H} R^o$. 
Then $N$, $N = M^{\sf s} - K^o$, is a real projective $n$-manifold 
with totally geodesic boundary; in fact, $M^{\sf s}$ decomposes 
into $N$ and $K$ along totally geodesic $(n-1)$-dimensional submanifold 
\[p(\bigcup_{R \in \cal H} \nu_R \cap M_h^{\sf s}).\]

We see that $p^{-1}(N)$ equals $M_h^{\sf s} - p^{-1}(K^o)$, and 
so $M_h^{\sf s} - p^{-1}(K^o)$ covers $N$. As we saw in Section 
\ref{sec:crescres}, we may choose a component $L'_i$ of 
$M_h^{\sf s} - p^{-1}(K^o)$ for each component $L_i$ of $N$ 
covering $L_i$ as a holonomy cover with developing map $\dev|L'_i$ 
and holonomy homomorphism as described there; $\coprod_{i=1}^n L'_i$ 
becomes a holonomy cover $N_h$ of $N$.  

We will show that $\che N_h$ includes no hemispheric $n$-crescent,
which implies that $N$ contains no concave affine manifold 
of type I. This completes the construction of the second column 
of arrows of equation \ref{eqn:arrows}.

The projective completion of $L'_i$ is denoted by $\che L'_i$. The 
projective completion $\che N_h$ of $N_h$ equals the disjoint union 
of $\che L'_i$. Since the inclusion map $i: L'_i \ra M_h^{\sf s}$ is 
distance decreasing, it extends to $\cim: \che N_h \ra \che M_h^{\sf s}$. 
If $\che N_h$ includes any hemispheric $n$-crescent $R$, then 
$\cim(R)$ is a hemispheric $n$-crescent by Proposition 
\ref{prop:corrsub}. First, since $\cim(R) \cap N_h$ is a subset 
of $p^{-1}(K)$ by the construction of $K$, and 
$\cim(R) \cap M_h^{\sf s} \supset R \cap N_h$, it follows that 
$R \cap N_h \subset p^{-1}(K)$. Second, since $N_h \cap p^{-1}(K)$ 
equals $p^{-1}(\Bd K)$, it follows that $N_h \cap p^{-1}(K)$ 
includes no open subset of $M_h^{\sf s}$.  Finally, this is 
a contradiction while $R \cap N_h \supset R^o$. 
Therefore, $\che N_h$ includes no hemispheric $n$-crescent.
(We may see after this stage, the completions of the covers of 
the spaces never includes any hemispheric $n$-crescents as 
the splitting and taking submanifolds do not affect this fact 
which we showed in Section \ref{sec:crescres}.)

Now we go to the second stage of the construction. Suppose that 
$\che N_h$ includes bihedral $n$-crescents and $A_2$ is the 
two-faced $(n-1)$-submanifold arising from bihedral $n$-crescents. 
Then we obtain the splitting $N^{\sf s}$ of $N$ along $p(A_2)$. 

We split $N_h$ along $A_2$ to obtain $N^{{\sf s} \prime}$. Then the 
holonomy cover $N_h^{\sf s}$ of $N^{\sf s}$ is a disjoint 
union of components of $N^{{\sf s} \prime}$ chosen for each component 
of $N^{\sf s}$. Let $\che N_h^{\sf s}$ denote the completion. 

The reasoning using Proposition \ref{prop:boocorr} as in the seventh 
paragraph above shows that $\Lambda(R) = \Lambda(S)$ or 
$\Lambda(R) \cap \Lambda(S) = \emp$ for every pair of bihedral 
$n$-crescents $R$ and $S$ in $\che N_h^{\sf s}$. The proof of Theorem 
\ref{thm:main2} shows that $N^{\sf s}$ contains a concave affine manifold 
$T$ of type II with the covering map 
\[ p|\bigcup_{R \in \cal B} \Lambda(R) \cap N_h^{\sf s}:
\bigcup_{R \in \cal B} \Lambda(R) \cap N_h^{\sf s} \ra K\] where 
$\cal B$ denotes the set of representatives of the equivalence 
classes of bihedral $n$-crescents in $\che N_h^{\sf s}$. 
And we see that $N^{\sf s} - T^o$ is a real projective manifold 
with convex boundary. We let $S = N^{\sf s} - T^o$. Thus $N^{\sf s}$ 
decomposes into $S$ and $T$.


For each component $S_i$ of $S$, we choose a component $S'_i$ of 
$N_h^{\sf s} - p^{-1}(T^o)$. Then $\coprod S'_i$ is a holonomy 
cover $S_h$ of $S$. The projective completion $\che S_h$ equals 
the disjoint union $\coprod \che S'_i$. As in the fourth paragraph 
above, we can show that $\che S_h$ includes no hemispheric or 
bihedral $n$-crescent using Proposition \ref{prop:corrsub}. 
If $S$ is not $(n-1)$-convex, then $\che S_h$ includes 
an $n$-crescent since Theorem \ref{thm:n-1conv} easily generalizes 
to the case when the real projective manifold $M$ has convex boundary 
instead of totally geodesic one or empty one. Since $\che S_h$ does 
not include a bihedral $n$-crescent, it follows that $S$ is $(n-1)$-convex. 

Now, we will show that the decomposition of Corollary \ref{cor:spl1}
is canonical. First, the two-faced submanifolds $A_1$ and $A_2$ 
are canonically defined. Now let $M = M^{\sf s}$ and $N = N^{\sf s}$
for convenience.

First, suppose that $M$ decomposes into $N'$ and $K'$ where $K'$ is 
a submanifold whose components are concave affine manifolds of type I 
and $N'$ is the closure of $M -K'$ and $N'$ includes no concave 
affine manifold of type I. This means that $\che N'_h$ includes no
hemispheric $n$-crescent where $\che N'_h$ is the completion of 
the holonomy cover $N'_h$ of $N'$. We will show that $N' = N$ and 
$K' = K$. 

Let $K'_i$, $i=1, \dots, n$, be the components of $K'$, and let 
$K'_{i, h}$ their respective holonomy cover, $\che K'_{i, h}$ the
projective completions. Each $\che K'_{i, h}$ equals a hemispheric 
$n$-crescent $R_i$. We claim that $p^{-1}(K')$ equals 
$\bigcup_{R \in \cal H} R \cap M_h$ where $\cal H$ is the set of 
all $n$-crescents in $\che M_h$. This will prove our claim in 
the above paragraph.

Each component $K^j_i$ of $p^{-1}(K'_i)$ is a holonomy cover of $K'_i$
(see Section \ref{sec:crescres}). Let $l^j_i$ denote the lift of 
the covering map $K'_{h, i} \ra K'_i$ to $K^j_i$ which is a homeomorphism.
$\dev\circ l^j_i$ is a developing map for $K'_{h, i}$ as it is a 
real projective map (see Ratcliff \cite{RT:94}). We may put a metric 
$\bdd$ on $K'_{h, i}$ induced from $\bdd$ on $\SI^n$, a quasi-isometric 
to any such choice of metric, using developing maps. Thus, we may 
identify $K^j_i$ with $K'_{h, i}$ and their completions for a moment. 
From the definition of concave affine manifolds of type I, 
the completion of $K'_{h, i}$ equals a hemispheric $n$-crescent 
$R_i$, and $K'_{h, i} = R_i \cap K'_{h, i}$. Proposition 
\ref{prop:corrsub} shows that there exists a hemispheric 
$n$-crescent $R'_i$ in $\che M_h$ with identical interior 
as that of $R_i$, and clearly $R'_i$ includes $K^j_i$ in 
$\che M_h$ so that $K^j_i = R'_i \cap M_h$. Since this is true 
for any component $K^j_i$, we have that $p^{-1}(K')$ is a union 
of hemispheric $n$-crescents intersected with $M_h$ and a subset of 
$\bigcup_{R \in \cal H} R \cap M_h$. 

Suppose that there exists a hemispheric $n$-crescent $R$ in $\che M_h$ 
so that $R \cap M_h$ is not a subset of $p^{-1}(K')$. Suppose $R$ 
meet $p^{-1}(K')$. Then $R$ meets a hemispheric $n$-crescent
$S$ where $S \cap M_h \subset p^{-1}(K')$. If $R$ and $S$ overlap,
then $R = S$, which is absurd. Thus, $R$ and $S$ may meet only 
at $\nu_R \cap M_h$ and $\nu_S \cap M_h$. This gives us a component 
of a pre-two-faced submanifold. By our assumption on $M$, this does 
not happen. Thus we see that $R$ and $S$ are disjoint, and $R^o$ is 
a subset of $p^{-1}(M - K^{\prime o})$. Since each component of 
$p^{-1}(M- K^{\prime o})$ is a holonomy cover of a component of 
$M- K^{\prime o}$, it follows that the completion of a holonomy 
cover of a component of $M-K^{\prime o}$ includes
a hemispheric $n$-crescent by Proposition \ref{prop:corrsub}. 
Obviously this hemispheric $n$-crescent intersected with the 
holonomy cover covers a concave affine manifold of type I.
As this contradicts our assumption about $M- K^{\prime o}$, 
we have that $p^{-1}(K')$ equals $\bigcup_{R \in \cal H} R \cap M_h$.  

Second, if $N$ decomposes into $S'$ and $T'$ where $T'$ is a concave
affine manifold of type II, and $S'$ do not contain any concave affine
manifold of type II, then we claim that $S' = S$ and $T' = T$. As above,
this follows if the completion of the holonomy cover of each component 
of $S'$ does not contain any bihedral $n$-crescents. The proof of this 
claim is exactly analogous to the above three paragraphs. This completes 
the proof of Corollary \ref{cor:spl1}.

To end this section, we supply the proof of Corollary \ref{cor:cdcr}.
Corollary \ref{cor:spl1} shows that the real projective surface
$\Sigma$ decomposes into convex real projective surfaces and 
concave affine surfaces along simple closed geodesics. By Euler 
characteristic considerations, their Euler numbers are all zero. 


\begin{rem}
Finally, we remark that Corollary \ref{cor:spl1} also holds if 
we simply assume that $M_h$ is not real projectively diffeomorphic 
to an open $n$-hemisphere or $n$-bihedron. In this case 
the decomposition could be a trivial one. Given the assumption 
in Corollary \ref{cor:spl1}, the decomposition has to take place 
at least once. 
\end{rem}

\vfill
\break

\section{Left-invariant real projective structures on Lie groups}
\label{sec:Liegp}

Finally, we end with an application to affine Lie groups.
Let $G$ be a Lie group with a left-invariant real projective 
structure, which means that $G$ has a real projective structure 
and the group of left translations are projective automorphisms.
Consider the holonomy cover $G_h$ of $G$. Then $G_h$ is also 
a Lie group with the induced real projective structure, which is
clearly left-invariant.

As before if $G$ is not $(n-1)$-convex, then 
$G_h$ is not projectively diffeomorphic to an open 
$n$-bihedron or an open $n$-hemisphere.

\begin{thm}\label{thm:liegp}
If $G$ is not $(n-1)$-convex as a real projective manifold, then 
the projective completion $\che G_h$ of $G_h$ includes an 
$n$-crescent $B$.
\end{thm}
\begin{pf}
This is proved similarly to Theorem \ref{thm:n-1conv} by 
a pull-back argument. The reason is that the left-action of $G_h$ 
on $G_h$ is proper and hence given two compact sets $K$ and $K'$ 
of $G_h$, the set $\{g \in G_h| g(K) \cap K' \ne \emp\}$ is 
a compact subset of $G_h$. That is, all arguments of Section 
\ref{sec:n-1conv} go through by choosing an appropriate sequence 
$\{g_i\}$ of elements of $G_h$ instead of deck transformations.  
Obviously, if $G_h$ contains a cocompact discrete subgroup, then
this is a corollary of Theorem \ref{thm:n-1conv}. But if not, we 
have to do this parallel argument.
\end{pf}

Suppose that from now on $\che G_h$ includes an $n$-crescent $B$.
Since the action of $G_h$ on $G_h$ is transitive, $\che G_h$ equals 
the union of $g(B)$ for $g \in G_h$. We claim that $g(B) \sim g'(B)$ 
for every pair of $g$ and $g'$ in $G_h$: That is, there exists 
a chain of $n$-crescents $B_i$, $i=1, \dots, k$, of same type as $B$ 
so that $B_1 = g(B)$, $B_i$ overlaps with $B_{i+1}$ for each 
$i = 1, \dots, k-1$, and $B_k = g'(B)$: Let $p$ be a point of $B^o$ 
and $B(p)$ a tiny ball of $p$ in $B^o$. We can choose a sequence $g_0, 
\dots, g_n$ with $g_0 = g$ and $g_n = g'$ where $g_{i+1}^{-1} g_i(p)$ 
belongs to $B(p)^o$. In other words, we require $g_{i+1}^{-1} g_i$ to 
be sufficiently close to the identity element. Then $g_i(B^o) \cap 
g_{i+1}(B^o) \ne \emp$ for each $i$. Hence $g(B) \sim g'(B)$ for 
any pair $g, g' \in G_h$.

If $B$ is an $n$-hemisphere, then we claim that $g(B) = B$ for all $g$ 
and $G_h = B \cap G_h$. The proof of this fact is identical to that of 
Theorem \ref{thm:hcresint} but we have use the following lemma instead of 
Lemma \ref{lem:twocres}.

\begin{lem}\label{lem:Gtwocres}
Suppose that $\dev: \che G_h \ra \SI^n$ is an imbedding onto the 
union of two $n$-hemispheres $H_1$ and $H_2$ meeting each other on 
an $n$-bihedron or an $n$-hemisphere. Then $H_1 = H_2$, and 
$G_h$ is projectively diffeomorphic to an open $n$-hemisphere. 
\end{lem}
\begin{pf} 
If $H_1$ and $H_2$ are different, as in the proof of Lemma 
\ref{lem:twocres} we obtain two $(n-1)$-dimensional hemispheres 
$O_1$ and $O_2$ in $G_h$ where a subgroup of index one or two 
in $G_h$ acts on. Since the action of $G_h$ is transitive, this 
is clearly absurd.  
\end{pf}

So if $B$ is an $n$-hemisphere, then we obtain $G_h = B \cap G_h$.
Since $G_h$ is boundaryless, $\delta B$ must consist of ideal points. 
This contradicts the definition of $n$-crescents. Therefore, every 
$n$-crescent in $\che G_h$ is a bihedral $n$-crescent.

The above shows that $\che G_h = \Lambda(B)$ for a bihedral 
$n$-crescent $B$, $G_h$ is a concave affine manifold of type II and 
hence so is $G$: To prove this, we need to show that two overlapping 
$n$-crescents intersect transversally as the proof for the Lie group 
case is slightly different. This is proved as in the proof of Theorem 
\ref{thm:transversal} using the above Lemma \ref{lem:Gtwocres} instead 
of Lemma \ref{lem:twocres}.

Let $H$ be a Lie group acting transitively on a space $X$. It is 
well-known that for a Lie group $L$ with left-invariant 
$(H, X)$-structure, the developing map is a covering map onto 
its image, an open subset (see Proposition 2.2 in Kim \cite{KimH:96}).
Thus $\dev: G_h \ra \SI^n$ is a covering map onto its image. 

Recall that $\dev$ maps $\Lambda(B)^o$ into an open subset of an 
open hemisphere $H$, and $\delta \Lambda(B)$ is mapped into the 
boundary $\SI^{n-1}$ of $H$. Each point of $G_h$ belongs to $S^o$ 
for an $n$-crescent $S$ equivalent to $R$ since the action of 
$G_h$ on $G_h$ is transitive (see above). Since each point of 
$\dev(G_h)$ belongs to the interior of an $n$-bihedron $S$ with 
a side in $\SI^{n-1}$, the complement of $\dev(G_h)$ is a closed 
convex subset of $H^o$. Thus, $\dev| G_h$ is a covering map 
onto the complement of a convex closed subset of $\bR^n$.
As $\tilde G$ covers $G_h$, we see that this completes the proof 
of Theorem \ref{thm:leftinv}.

An affine $m$-convexity for $1 \leq m < n$ is defined as follows.
Let $M$ be an affine $n$-manifold, and let $T$ be 
an affine $(m+1)$-simplex in $\bR^n$ with sides $F_1, F_2,
\dots, F_{m+2}$. Then every nondegenerate affine map 
$f: T^o \cup F_2 \cup \dots \cup F_{m+2} \ra M$ extends to 
one $T \ra M$ (see \cite{uaf:97} for more details).

If $G$ has a left-invariant affine structure, then $G$ has a 
compatible left-invariant real projective structure. By Theorem 
\ref{thm:leftinv}, $G$ is either $(n-1)$-convex or $G_h$ is 
a concave affine manifold of type II. It is easy to see that the 
$(n-1)$-convexity of $G$ in the real projective sense implies 
the $(n-1)$-convexity of $G$ in the affine sense.

As before, if $G_h$ is a concave affine manifold of type II, 
the argument above shows that $G_h$ is mapped by $\dev$ to the 
complement of a closed convex set in $\bR^n$. This completes 
the proof of Corollary \ref{cor:leftinvco}.

Finally, we easily see that the following theorem holds with 
the same proof as the Lie group case:

\begin{thm}\label{thm:homg}
Let $M$ be a homogenous space on which a Lie group $G$ 
acts transitively and properly. Suppose $M$ have 
a $G$-invariant real projective structure. Then $M$ is either 
$(n-1)$-convex, or $M$ is concave affine. The same holds for 
$G$-invariant affine structures.
\end{thm}

\vfill
\break

\begin{appendix}

\section{Proof of Theorem \ref{thm:1conv}}

The proofs are a little sketchy here; however, they are elementary.

\begin{thm}
The following are equivalent\/{\rm :}
\begin{enumerate}
\item $M$ is $1$-convex.
\item $M$ is convex.
\item $M$ is real projectively isomorphic to a quotient of 
a convex domain in $\SI^n$.
\end{enumerate}
\end{thm}

(1)$\ra$(2): Since $M$ is $1$-convex, $\tilde M$ is $1$-convex. 
Any two points $x$ and $y$ in $\tilde M$ are connected by a chain of 
segments $s_i$, $i=1, \dots, n$, of $\bdd$-length $< \pi$ with 
endpoints $p_i$ and $p_{i+1}$ so that $s_i \cap s_{i+1} = \{p_{i+1}\}$ 
exactly. This follows since any path may be covered by tiny balls 
which are convex. We will show that $x$ and $y$ is connected by 
a segment of $\bdd$-length $\leq \pi$.

Assume that $x$ and $y$ are connected by such a chain with $n$ 
being a minimum. We can assume further that $s_i$ are in general 
position, i.e., $s_i$ and $s_{i+1}$ do not extend each other as 
an imbedded geodesic for each $i=1, \dots, n-1$, which may be 
achieved by perturbing the points $p_2, \dots, p_n$, unless $n =2$ 
and $s_1 \cup s_2$ form a segment of $\bdd$-length $\pi$; in which 
case, there is nothing to prove since $s_1\cup s_2$ is the segment 
we need. To show we can achieve this, we take a maximal sequence of 
segments which extend each other as geodesics. Suppose that 
$s_i, s_{i+1}, \dots, s_{j}$ form such a sequence for $j > i$. 
Then the total length of the segment will be less than $\pi|j -i|$. 
We divide the sequence into new segments of equal $\bdd$-length 
$s'_i, s'_{i+1}, \dots, s'_j$ where $s'_k$ has new endpoints 
$p'_k, p'_{k+1}$ for $i \leq k \leq j$ where $p'_i = p_i$ and 
$p'_{j+1} = p_{j+1}$. Then we may change $p'_k$ for 
$k = i+1, \dots, j$ toward one-side of the segments by a small 
amount generically. Then we see that new segments $s''_i, s''_{i+1}, 
\dots, s''_j$ are in general position together with $s_{i-1}$ and 
$s_{j+1}$. This would work unless $j-i =2$ and the total $\bdd$-length 
equals $\pi$ since changing $p'_{i+1}$ still preserves 
$s''_i \cup s''_{i+1}$ to be a segment of $\bdd$-length
$\pi$. However, since $n \geq 3$, we may move $p'_i$ or $p'_{j+1}$ 
in some direction to put the segments into the general position. 

Let us choose a chain $s_i$, $i=1, \dots, n$, with minimal number 
of segments in general position. We assume that we are not in case when 
$n=2$ and  $s_1 \cup s_2$ forming a segment of $\bdd$-length $\pi$,

We show that the number of the segments equals one, which shows that 
$\tilde M$ is convex. If the number of the segments is not one, then we 
take $s_1$ and $s_2$ and parameterize each of them 
by projective maps $f_i:[0, 1] \ra s_i$, $i=1,2$, so that $f_i(0) = p_2$. 
Then since $p_2$ is a point of its tiny ball $B(p_2)$, it follows that 
there exists a nondegenerate real projective map $f_t: \triangle_t \ra B(p_2)$ 
where $\triangle_t$ is a triangle in $\bR^2$ with vertices $(0,0), (t,0),$
and $(0, t)$ and $f_t(0, 0) = p_2$, $f_t(s, 0) = f_1(s)$ and 
$f_t(0, s) = f_2(s)$ for $0 \leq s \leq t$. $f_t:\triangle_t \ra \tilde M$ 
is always an imbedding since $\dev\circ f_t$ is a nondegenerate 
projective map $\triangle_t \subset \SI^n \ra \SI^n$.

We consider the subset $A$ of $[0,1]$ so that $f_t: \triangle_t \ra 
\tilde M$ is defined. Then $A$ is open in $[0,1]$ since as $f(\triangle_t)$ 
is compact, there exists a convex neighborhood of it in $\tilde M$
where $\dev$ restricts to an imbedding. 
$A$ is closed by $1$-convexity: we consider the union 
$K =\bigcup_{t \in A} f_t(A)$. Then the closure of $K$ in $\che M$ is 
a compact triangle in $\che M$ with two sides in $s_1$ and $s_2$. Thus, 
two sides of $K$ and $K^o$ are in $\tilde M$. By $1$-convexity, $K$ 
itself is in $\tilde M$. Thus, $\sup A$ also belongs to $A$ and $A$ 
is closed. Hence $A$ must equal $[0,1]$ and there exists a segment 
$s'_1$ of $\bdd$-length $< \pi$, namely $f_1(\ovl{(1,0)(0,1)})$, 
connecting $p_1$ and $p_3$. This contradicts the minimality, and 
$x$ and $y$ are connected by a segment of $\bdd$-length $\leq \pi$.

(2)$\ra$(3) We choose a point $x$ in $\tilde M$. For each point 
$y$ of $\tilde M$, there exists a segment of $\bdd$-length $\leq \pi$
connecting $x$ and $y$. This implies that $\dev(\tilde M)$ is 
a convex subset of $\SI^n$. As $\tilde M$ is convex, the closure of 
$\tilde M$ in $\che M$ is convex as we explain in Section 
\ref{sec:conv2}. Thus $\che M$ is a tame set, and 
$\dev:\che M \ra \SI^n$ is an imbedding onto a convex subset of 
$\SI^n$. Hence $\dev|\tilde M$ is an imbedding onto
a convex subset of $\SI^n$. Now, the equation $\dev\circ \vth =
h(\vth) \circ \dev$ holds for each deck transformation $\vth$ of $\tilde M$.
Therefore, it follows that $\dev$ induces a real projective diffeomorphism 
$\tilde M /\pi_1(M) \ra \dev(\tilde M)/h(\pi_1(M))$.

(3)$\ra$(1) This part is straightforward using the classification of 
convex sets in Theorem \ref{thm:class1} since $\tilde M$ can be
identified with a convex domain in $\SI^n$.

\vfill
\break

\section{The Proof of Proposition \ref{prop:blowup}}\label{app:B}

\begin{prop}
Suppose we have a sequence of $\eps$-$\bdd$-balls $B_i$ in a real 
projective sphere $\SI^n$ for some $n \geq 1$ and a fixed positive 
number $\eps$ and a sequence of projective maps $\vpi_i$. Assume 
the following\/{\rm :}
\begin{itemize}
\item The sequence of $\bdd$-diameters of $\vpi_i(B_i)$ goes to zero. 
\item $\vpi_i(B_i)$ converges to a point, say $p$. 
\item For a compact $n$-ball neighborhood $L$ of $p$, 
$\vpi_i^{-1}(L)$ converges to a compact set $L_\infty$.
\end{itemize}
Then $L_\infty$ is an $n$-hemisphere.
\end{prop}

Recall that $\bR^{n+1}$ has a standard euclidean metric and $\bdd$ on 
$\SI^n$ is obtained from it by considering $\SI^n$ as the standard unit 
sphere in $\bR^{n+1}$. 

The Cartan decomposition of Lie groups states that a real reductive 
Lie group $G$ can be written as $KTK$ where $K$ is a compact Lie group 
and $T$ is a maximal real tori. Since $\Aut(\SI^n)$ is isomorphic to 
$\SL_\pm(n+1, \bR)$, we see that $\Aut(\SI^n)$ can be written as 
$O(n+1)D(n+1)O(n+1)$ where $O(n+1)$ is the orthogonal group acting 
on $\SI^n$ as the group of isometries and $D(n+1)$ is the group of 
determinant $1$ diagonal matrices with positive entries listed 
in nonincreasing order where $D(n+1)$ acts in $\SI^n$ as a subgroup 
of $\GL(n+1, \bR)$ acting in the standard manner on $\SI^n$. In other 
words, each element $g$ of $\Aut(\SI^n)$ can be written
as $i(g)d(g)i'(g)$ where $i(g), i'(g)$ are isometries and 
$d(g) \in D(n+1)$ (see Carri\`ere \cite{Car:89} and Choi \cite{uaf:97},
and also \cite{uaf2:97}).

We may write $\vpi_i$ as $K_{1, i}\circ D_i \circ K_{2, i}$ 
where $K_{1, i}$ and $K_{2, i}$ are $\bdd$-isometries of $\SI^n$ and 
$D_i$ is a projective map in $\Aut(\SI^n)$ represented by a diagonal 
matrix of determinant $1$ with positive entries. More precisely, 
$D_i$ has $2n+2$ fixed points $[\pm e_0], \dots, [\pm e_n]$, 
the equivalence classes of standard basis vectors 
$\pm e_0, \dots, \pm e_n$ of $\bR^{n+1}$, and $D_i$ has a matrix 
diagonal with respect to this basis; the diagonal entries $\lambda_i$, 
$i= 0, 1, \dots, n$, are positive and in nonincreasing order. 
Let $O_{[e_0]}$ denote the open hemisphere containing $[e_0]$ whose 
boundary is the great sphere $\SI^{n-1}$ containing $[\pm e_j]$ for all 
$j$, $j \geq 1$, and $O_{[-e_0]}$ that containing $[-e_0]$ with the same 
boundary set. 

Recall that $\vpi_i = K_{1, i}\circ D_i \circ K_{2, i}$. Let us 
denote by $q_i = \vpi_i(p_i)$ for the $\bdd$-center $p_i$ of 
the ball $B_i$. Since $\vpi_i(B_i)$ converges to $p$, and 
$K_{1, i}$ is an isometry, the sequence of the $\bdd$-diameter 
of $D_i\circ K_{2, i}(B_i)$ goes to zero as $i \ra \infty$. 
We \awlg that $D_i(K_{2, i}(B_i))$ converges to a set consisting 
of a point by choosing a subsequence if necessary. By the following 
lemma \ref{lem:shrinkim}, $D_i(K_{2, i}(B_i))$ converges to one 
of the attractors $[e_0]$ and $[-e_0]$. We \awlg that $D_i(K_{2, i}(B_i))$ 
converges to $[e_0]$.  

Since $L$ is an $n$-dimensional ball neighborhood of $p$, $L$ 
includes a $\bdd$-ball $B_\delta(p)$ in $\SI^n$ with center $p$ with 
radius $\delta$ for some positive constant $\delta$. 
There exists a positive integer $N$ so that for $i > N$,  
we have \[ \vpi_i(p_i) \subset B_{\delta/2}(p)\] for the $\bdd$-ball 
$B_{\delta/2}(p)$ of radius $\delta/2$ in $\SI^n$.
Hence $B_{\delta/2}(q_i)$ is a subset of $L$ for $i > N$.

Since $K_{1, i}^{-1}(q_i) = D_i \circ K_{2, i}(p_i)$, 
the sequence $K_{1, i}^{-1}(q_i)$ converges to $[e_0]$ by the 
second paragraph above. There exists an integer $N_1 > N$ such that 
$K_{1, i}^{-1}(q_i)$ is of $\bdd$-distance less than $\delta/4$ 
from $[e_0]$ for $i > N_1$. Since $K_{1, i}^{-1}$ is a $\bdd$-isometry, 
$K_{1, i}^{-1}(B_{\delta/2}(q_i))$ includes the ball 
$B_{\delta/4}([e_0])$ for $i > N_1$. Hence $K_{1, i}^{-1}(L)$ 
includes $B_{\delta/4}([e_0])$ for $i > N_1$.

Since $[e_0]$ is an attractor under the action of the sequence 
$\{D_i\}$ by Lemma \ref{lem:shrinkim}, the images of 
$B_{\delta/4}([e_0])$ under $D_i^{-1}$ eventually include any compact 
subset of $O_{[e_0]}$. Thus, $D_i^{-1}(B_{\delta/4}([e_0]))$ converges 
to $\clo(O_1)$ geometrically, and up to a choice of subsequence 
$K_{2, i}^{-1}\circ D_i^{-1}(B_{\delta/4}([e_0]))$ converges to
an $n$-hemisphere. The equation
\begin{eqnarray}\label{eqn:a1}
\vpi_i^{-1}(L) &= K_{2, i}^{-1} \circ D_i^{-1} \circ K_{1, i}^{-1}(L) 
\nonumber \\ 
& \supset K_{2, i}^{-1}\circ D_i^{-1}(B_{\delta/4}([e_0])).
\end{eqnarray}
shows that 
$\vpi_i^{-1}(L)$ converges to an $n$-hemisphere.
\qed

The straightforward proof of the following lemma is left to the reader.
\begin{lem}\label{lem:shrinkim}
Let $K_i$ be a sequence of $\eps$-$\bdd$-balls in $\SI^n$ and $d_i$ 
a sequence of automorphisms of $\SI^n$ that are represented by diagonal 
matrices of determinant $1$ with positive entries for the standard basis
with the first entry $\lambda_i$ the maximum. Suppose $d_i(K_i)$ converges 
to the set consisting of a point $y$. Then there exists an integer $N$ 
so that for $i > N$, the following statements hold\/{\rm :} 
\begin{enumerate}
\item $[e_0]$ and $[-e_0]$ are attracting fixed points of $d_i$. 
\item $y$ equals $[e_0]$ or $[-e_0]$.   
\item The eigenvalue $\lambda_i$ of $d_i$ corresponding to $e_0$ and 
$-e_0$ is strictly larger than the eigenvalues corresponding to 
$\pm e_j$, $j=1, \dots, n$.
\item $\lambda_i/\lambda'_i \ra +\infty$ for the maximum eigenvalue 
$\lambda'_i$ of $d_i$ corresponding to $\pm e_j, j = 1 \dots n$.  
\end{enumerate}
\end{lem}

\end{appendix}

\vfill
\break

\bibliographystyle{plain}
\bibliography{proj}

\end{document}